\documentclass[11pt,a4paper]{article}
\usepackage[T1]{fontenc}
\usepackage[utf8]{inputenc}
\usepackage{lmodern}
\usepackage{amsmath,amssymb,amsthm,mathtools,bm,mathrsfs}
\usepackage[margin=28mm]{geometry}
\usepackage{microtype}
\usepackage{booktabs,longtable,array,enumitem}
\usepackage{xcolor}
\definecolor{linkblue}{RGB}{30,55,95}
\usepackage[colorlinks=true,linkcolor=linkblue,citecolor=linkblue,urlcolor=linkblue]{hyperref}
\hypersetup{pdftitle={Spectral gaps near one quarter for Weil--Petersson random surfaces with growing cusp number}}
\numberwithin{equation}{section}
\newtheorem{theorem}{Theorem}[section]
\newtheorem{lemma}[theorem]{Lemma}
\newtheorem{proposition}[theorem]{Proposition}
\newtheorem{corollary}[theorem]{Corollary}
\theoremstyle{definition}

\theoremstyle{remark}

\newcommand{\E}{\mathbb E}
\newcommand{\PP}{\mathbb P}
\newcommand{\R}{\mathbb R}

\newcommand{\Z}{\mathbb Z}
\newcommand{\ind}{\mathbf 1}
\DeclareMathOperator{\supp}{supp}
\DeclareMathOperator{\Spec}{Spec}
\DeclareMathOperator{\Tr}{Tr}

\DeclareMathOperator{\TV}{TV}
\DeclareMathOperator{\SpG}{SpG}
\title{Spectral gaps near one quarter for Weil--Petersson\\random surfaces with growing cusp number}
\author{Guan HUANG}
\date{}
\begin{document}
\maketitle
\begin{abstract}
We prove that the spectral gap of a Weil--Petersson random hyperbolic
surface of genus $g$ with $n=n(g)\ge1$ labelled cusps converges in
probability to $1/4$ as $g\to\infty$, for every sequence
$n=o(\sqrt g)$. More precisely, for each fixed
$0<\varepsilon<1/4$, there are constants $C_\varepsilon\ge2$ and
$g_\varepsilon\ge3$, independent of $n$, such that the expected
number of nonzero $L^2$ eigenvalues in $(0,1/4-\varepsilon]$,
counted with multiplicity and including residual eigenvalues, is
at most $C_\varepsilon(n+1)^2/g$ whenever $g\ge g_\varepsilon$
and $g\ge C_\varepsilon(n+1)^2$. This proves a conjecture of
Shen and Wu. A matching lower bound shows that the count estimate
is optimal up to constants for each fixed $\varepsilon$, and that
$n=o(\sqrt g)$ is necessary for convergence of the gap to $1/4$.
The proof combines volume estimates uniform in the cusp count with
polynomial approximation and a positive trace on truncated surfaces.
\end{abstract}
\section{Introduction}\label{sec:introduction}

For a complete finite-area hyperbolic surface with cusps, the
continuous spectrum begins at $1/4$. We prove that the spectral gap
of a Weil--Petersson random surface converges in probability to this
threshold as the genus $g$ tends to infinity, for every positive
integer cusp sequence $n=n(g)=o(\sqrt g)$. Our quantitative result
bounds the expected number of nonzero eigenvalues in
$(0,1/4-\varepsilon]$ for each fixed $0<\varepsilon<1/4$,
uniformly in the cusp count.

Let $\mathcal M_{g,n}$ be the moduli space of connected, complete,
finite-area oriented hyperbolic surfaces of genus $g$ with $n$
labelled cusps, without boundary or orbifold points. We equip it
with normalized Weil--Petersson measure $\PP_{g,n}$ and write
$\E_{g,n}$ for expectation. The spectral statements concern integers
$g\ge3$ and $n\ge1$; the geometric expansions also allow $n=0$.
The hyperbolic area of each such surface is
\begin{equation}\label{intro:area}
 \mathcal A_{g,n}=2\pi(2g-2+n).
\end{equation}
We reserve $V_{g,n}$ for the Weil--Petersson volume of moduli space.

For $X\in\mathcal M_{g,n}$ with $n\ge1$, the nonnegative self-adjoint Laplacian
$\Delta_X$ has a simple eigenvalue zero, continuous spectrum
$[1/4,\infty)$, and possibly nonzero discrete eigenvalues below $1/4$.
For $0<\varepsilon<1/4$, define
\begin{equation}\label{intro:count}
 N_\varepsilon(X)=
 \#\{\lambda\in\Spec_{L^2}(\Delta_X):0<\lambda\le1/4-\varepsilon\}.
\end{equation}
Here the count is taken with multiplicity and includes both cuspidal
and residual eigenvalues. Write
$\SpG(X)=\inf\Spec(\Delta_X|_{\mathbf1^\perp})$, where
$\mathbf1^\perp$ is the orthogonal complement of the constant
functions in $L^2(X)$. Thus $\SpG(X)=1/4$ if there is no nonzero
eigenvalue below $1/4$, and
$\{\SpG(X)\le1/4-\varepsilon\}=\{N_\varepsilon(X)\ge1\}$.

\begin{theorem}[Spectral gaps throughout the sub-square-root cusp range]
\label{thm:subcritical}
For every fixed $0<\varepsilon<1/4$ there are constants
$C_\varepsilon\ge2$ and $g_\varepsilon\ge3$, independent of $n$,
such that, for all integers $g\ge3$ and $n\ge1$ satisfying
\begin{equation}\label{intro:subcritical-domain}
 g\ge g_\varepsilon,\qquad g\ge C_\varepsilon(n+1)^2,
\end{equation}
one has
\begin{equation}\label{intro:subcritical-count}
 \PP_{g,n}\{\SpG(X)\le1/4-\varepsilon\}
 \le\E_{g,n}N_\varepsilon(X)
 \le C_\varepsilon\frac{(n+1)^2}{g}.
\end{equation}
In particular, for every positive integer sequence $n(g)=o(\sqrt g)$,
the expected count and this probability tend to zero. For every
envelope $n_{\max}(g)\ge1$ with $n_{\max}(g)=o(\sqrt g)$, the following bound holds
for all sufficiently large $g$, uniformly over integer cusp counts
$1\le n\le n_{\max}(g)$:
\begin{equation}\label{intro:subcritical-uniform}
 \sup_{1\le n\le n_{\max}(g)}\E_{g,n}N_\varepsilon
 \le C_\varepsilon\frac{(n_{\max}(g)+1)^2}{g}\longrightarrow0.
\end{equation}
\end{theorem}

For example, the expected count is
$O_{\varepsilon,\delta}(g^{2\delta-1})$ uniformly for
$1\le n\le g^\delta$, for each fixed $0\le\delta<1/2$.
It is also $O_\varepsilon((\log\log g)^{-2})$ uniformly for
$1\le n\le\sqrt g/\log\log g$. Thus the estimate covers cusp
sequences outside every range $n=O(g^\delta)$ with fixed $\delta<1/2$.
The same bounds hold for the probability in
\eqref{intro:subcritical-count}. Throughout, $\varepsilon$ is fixed:
the estimate gives neither a shrinking spectral window nor a bound
tending to zero when $n\asymp\sqrt g$. The proof is completed in
Section~\ref{sec:additivespectral}.

A matching lower bound shows that the order $(n+1)^2/g$ is optimal
for both the expected count and the probability in the quadratic
domain, for each fixed $\varepsilon$.
Proposition~\ref{sharp:lower} gives the lower bound
$c_\varepsilon\min\{(n+1)^2/g,1\}$ for all $n\ge1$ once $g$
exceeds a bound depending only on $\varepsilon$.
Corollary~\ref{sharp:consequences} gives the two-sided
estimate and shows that convergence of the spectral gap to $1/4$
holds if and only if $n=o(\sqrt g)$.

\subsection{Relation to previous results}\label{sec:related}
For closed Weil--Petersson random surfaces, Mirzakhani established a
uniform positive spectral gap with probability tending to one
\cite{Mir13}. Wu and Xue \cite[Theorem~1]{WX22} and, independently,
Lipnowski and Wright \cite[Theorem~1.1]{LW24} raised the lower bound to
$3/16-\varepsilon$. Anantharaman and Monk subsequently reached
$1/4-\varepsilon$ for every fixed $\varepsilon>0$ through their work on
Friedman--Ramanujan functions and geometric expansions
\cite[Theorem~1.6]{AMGap}. Hide, Macera and Thomas gave a different proof
using the polynomial method and allowed the spectral window to
shrink at a polynomial rate: for some fixed $\eta>0$, the gap is at
least $1/4-O(g^{-\eta})$ with probability tending to one
\cite[Theorem~1.1]{HMT25}.

For Weil--Petersson surfaces with cusps, the known lower bounds
depend on the growth of the cusp count. Hide proved that, for each
fixed $0\le\sigma<1/2$ and $n=O(g^\sigma)$, the spectral gap exceeds
\[
 \frac14-\left(\frac{2\sigma+1}{4}\right)^2-\varepsilon
\]
with probability tending to one \cite[Theorem~1.3]{Hide23}.
He, Wu and Xue improved this bound to
\[
 \frac14-\left(\frac{1}{6(1-\sigma)}\right)^2-\varepsilon
\]
in the same fixed-power range \cite[Theorem~1.2]{HWX26}.
When the number of cusps is bounded, the two bounds are $3/16-\varepsilon$ and
$2/9-\varepsilon$, respectively. The proof of the latter identifies a
second-order cancellation between geometric contributions.

Shen and Wu obtained a positive spectral lower bound throughout
$n=o(\sqrt g)$ and formulated the corresponding $1/4-\varepsilon$
statement as Conjecture~4 \cite[Theorem~3 and Conjecture~4]{SW25}.
The same conjecture is restated in \cite[Conjecture~1.1]{HWX26}.
Shen and Wu also showed that the spectral gap tends to zero in
probability when $n/\sqrt g\to\infty$ and $n/g\to0$
\cite[Theorem~1]{SW25}. At the critical scale
$n/\sqrt g\to c\in(0,\infty)$, their Theorem~2 gives an explicit
limit for the distribution function of the Cheeger constant on an
interval near zero. Together, these results identify $\sqrt g$ as
a natural scale for the change in spectral behavior.

Li and Wang obtained large-genus volume asymptotics uniform over
all cusp counts \cite[Theorem~1.1]{LiWang26}. They also extended the
vanishing-gap conclusion to every sequence with
$n/\sqrt g\to\infty$ as $g\to\infty$, without assuming $n/g\to0$
\cite[Corollary~1.6]{LiWang26}.

Together with Shen and Wu's critical-scale result and the inverse
Cheeger inequality, this already implies that $n=o(\sqrt g)$ is
necessary for convergence in probability of the spectral gap to $1/4$;
see \cite[Theorems~2 and~11]{SW25} and
\cite[Corollary~1.6]{LiWang26}.
Proposition~\ref{sharp:lower} gives a uniform quantitative form of
this obstruction.

Theorem~\ref{thm:subcritical} proves Shen--Wu's Conjecture~4 for every
positive integer sequence $n=o(\sqrt g)$, with the quantitative bound
$C_\varepsilon(n+1)^2/g$ on the expected count. Residual eigenvalues
are included, as they are in the spectral-gap statements cited above.

For comparison, random-cover results use a different probability
measure: one samples uniformly from homomorphisms of the base
surface group to the symmetric group $S_d$, obtaining degree-$d$
covers with labelled fibre. For a fixed finite-area cusped base
and every fixed $\varepsilon>0$,
Hide and Magee proved that, with probability tending to one, the
spectrum below $1/4-\varepsilon$ consists exactly of the eigenvalues
inherited from the base, with the same multiplicities \cite[Theorem~1.1]{HM23}. Magee, Puder and van Handel
proved the corresponding statement for closed bases, with threshold
$1/4-o(1)$ \cite[Theorem~1.2]{MPvH}. For closed bases, Hide, Macera and
Thomas refined the threshold to $1/4-cd^{-\beta}$, where $d$ is the
covering degree and $c,\beta>0$ \cite[Theorem~1.1]{HMTCover25}.

In the Brooks--Makover model, Shen and Wu obtained a spectral gap greater than
$1/4-C/\log M$ with probability tending to one for the closed Belyi
surfaces obtained by conformally compactifying the gluings of $2M$
ideal triangles
\cite[Theorem~1.1]{SWBelyi26}. These sampling measures differ from the Weil--Petersson measure
with prescribed genus and cusp number used here.

The geometric argument uses Mirzakhani's integration formula and
volume-polynomial identities \cite{Mir07,Mir07Intersection,Mir08}
and the large-genus expansions of Mirzakhani and Zograf
\cite{MZ15}. Anantharaman and Monk developed integration and expansion
formulas organized by the subsurface filled by a geodesic
\cite[\S4.1, Theorems~5.7, 5.15 and~6.1]{AMI}.
For a fixed number of marked points, Guillermit refines the expansion-order
bounds and obtains volume asymptotics for boundary lengths on the
$\sqrt g$ scale, while leaving the dependence on the marked-point
count implicit \cite[Theorems~2, 4 and~5, Remark~3]{Gui26}.
We use the joint estimates of Hide, Macera and Thomas, which track both
the marked-point count and the expansion order
\cite[Theorem~4.2]{HMT25}. On the analytic side, the polynomial method
for strong convergence introduced by Chen, Garza-Vargas, Tropp and
van Handel \cite{CGTV24} was extended to surface-group permutation
representations in \cite{MPvH} and adapted to the closed
Weil--Petersson model in \cite{HMT25}.

The uniform count bound uses an additive genus condition whose
quadratic cusp term is independent of the detector degree, and a
bound for the first orbital correction independent of that degree.
All polynomial tests of a given maximum degree are evaluated on
one cusp truncation.

\subsection{A finite-parameter estimate}\label{sec:quantitative}

The proof applies a degree-$Q$ polynomial to a bounded spectral
multiplier. For $\lambda\ge1/4$, the chosen detector test decays exponentially
in $Q$. The additive genus condition below permits
$Q\asymp\log g$, which yields Theorem~\ref{thm:subcritical}.

Fix an even nonnegative function $b_0\in C_c^\infty(\R)$ of integral one,
and set $h=b_0*b_0$ and $f=\widehat h$. Define
$L=\max\{1,\sup\{|u|:u\in\supp h\}\}$, which is determined by $f$.

Our Fourier transform conventions are
\begin{equation}\label{intro:fourier}
 \widehat H(r)=\int_{\R}H(u)e^{-iru}\,du,\qquad
 \check F(u)=\frac1{2\pi}\int_{\R}F(r)e^{iru}\,dr.
\end{equation}
The function $f$ is even and entire. On the real axis,
$f=\widehat b_0^{\,2}$ and $0\le f\le1$; on the positive imaginary
axis, $t\mapsto f(it)$ is strictly increasing. In particular,
\begin{equation}\label{intro:a}
 a=f(i/2)>1.
\end{equation}
In the spectral variable $\lambda=1/4+r^2$, the multiplier $f(r)$
maps $\lambda\ge1/4$ into $[0,1]$ and the zero eigenvalue to $a$.
We fix $f$ throughout. For $0<\varepsilon<1/4$, define
\begin{equation}\label{intro:detectorconstants}
 b_\varepsilon=f(i\sqrt\varepsilon),\qquad
 \kappa=\frac{(b_\varepsilon-1)^2}{16a^2}.
\end{equation}
Then $1<b_\varepsilon<a$ and $\kappa>0$.
Let $Q_\varepsilon$ be the integer threshold defined in
\eqref{poly:detector-constants} of Section~\ref{sec:detector}; it depends
only on $f$ and $\varepsilon$.
Eigenvalues counted by $N_\varepsilon$ are mapped into
$[b_\varepsilon,a)$. The detector constructed in
Section~\ref{sec:detector} is uniformly positive there and
exponentially small on $[0,1]$. The lower bound
$Q\ge Q_\varepsilon$ ensures the detector and retained-mass estimates
used below.

Put $s=n+1$. For each surface $X$ and degree bound $Q$, we use
one cusp truncation for all real polynomial tests $p$ of degree
at most $Q$ with $p(0)=0$. Its height is
\begin{equation}\label{intro:cutoff}
 Y_Q=3e^{LQ}.
\end{equation}
Writing $p(x)=x\widetilde p(x)$, the resulting traces are integrals
of $\widetilde p$ against a common finite positive measure on $[0,a]$
(Theorem~\ref{trace:exact}).

\begin{theorem}[Finite-parameter spectral count]\label{thm:sharpfinite}
There are fixed constants $C_{\rm add},D_0,D_3\ge2$, depending only
on the multiplier $f$, such that for every fixed
$0<\varepsilon<1/4$ there is $C_{f,\varepsilon}$ with the following
property. For all integers $g\ge3$, $n\ge1$ and $Q\ge Q_\varepsilon$
satisfying
\begin{equation}\label{intro:additive-domain}
 g\ge C_{\rm add}s^2+C_{\rm add}Q^{D_3}s,
\end{equation}
one has
\begin{equation}\label{intro:additive-count}
 \begin{split}
 \PP_{g,n}(N_\varepsilon\ge1)\le\E_{g,n}N_\varepsilon
 \le C_{f,\varepsilon}\bigg\{
 &(g+nQ^2)e^{-\kappa Q}+\frac{s^2}{g}
       +\frac{s^3Q^{3D_0}}{g^2}\\
 &+\frac{s^2}{g}2^{-Q}
       +\frac{n}{\mathcal A_{g,n}Y_Q}\bigg\}.
 \end{split}
\end{equation}
\end{theorem}

The constants $D_0$, $D_3$ and $C_{\rm add}$ are specified in
Lemma~\ref{sharppm:coefficients} and Proposition~\ref{addspec:pm};
Theorem~\ref{addspec:spectral} proves the estimate. To deduce
Theorem~\ref{thm:subcritical}, take
$Q=\lceil K_\varepsilon\log g\rceil$, where
$K_\varepsilon=\max\{1,4/\kappa,2/L\}$.
If $g\ge2C_{\rm add}s^2$, the quadratic term in
\eqref{intro:additive-domain} is at most $g/2$. The linear term is
also at most $g/2$ once $g$ exceeds a bound depending only on $f$
and $\varepsilon$, since $Q^{D_3}=o(\sqrt g)$.
Moreover, the higher-order term in \eqref{intro:additive-count},
relative to $s^2/g$, is
$O_{f,\varepsilon}(Q^{3D_0}/\sqrt g)=o(1)$ uniformly in this range.
The remaining terms are bounded by a constant multiple of $s^2/g$.
After fixing $f$, this yields the stated constants
$C_\varepsilon$ and $g_\varepsilon$, independent of $n$.

\subsection{Proof structure and conventions}\label{sec:strategy}

The proof has three main steps.

First, we obtain geometric expansions uniform in the cusp count and
the expansion order. We integrate over simple and nonsimple geodesics,
including repetitions of primitive geodesics. For a nonsimple
geodesic, the integration is organized by the subsurface it fills,
called its carrier. Complementary-volume estimates and labelled
carrier integration give an orbital expansion whose coefficients are
independent of genus and polynomial in $n$ (Sections~\ref{sec:volume},
\ref{sec:orbital} and~\ref{sec:cusppolynomiality}). Retaining the
full binomial decay in the complementary-volume tails gives an
additive genus condition: a quadratic cusp term with an absolute
coefficient, plus a linear cusp term whose coefficient depends on
the expansion order and length cutoff
(Section~\ref{sec:additivethreshold}).

Second, we use one cusp cutoff for all monomials of a given maximum
degree. This gives a common positive trace functional with an exact
formula that includes the parabolic contribution. The leading orbital
functional is represented by a unit atom at $a=f(i/2)$ plus a finite
signed measure on $[0,1]$. Subtracting the constant eigenfunction's
contribution to the truncated trace cancels this atom up to the
explicitly controlled mass removed by truncation
(Section~\ref{sec:trace}). Every eigenfunction with eigenvalue counted by
$N_\varepsilon$, including a residual eigenfunction, retains a
uniformly positive fraction of its mass, allowing the polynomial
detector to bound the eigenvalue count
(Section~\ref{sec:polynomial}).

Third, we bound the first correction---the coefficient of $g^{-1}$
in the orbital expansion---independently of the detector degree.
For each fixed test, this coefficient is a polynomial of degree at
most two in $n$. Interpolation over auxiliary genera bounds the
canonical coefficients at finitely many cusp counts; their polynomial
dependence on $n$ extends the bounds to all cusp counts
(Section~\ref{sec:sharppm}). These bounds transfer the expansion
to polynomial tests at the additive genus threshold
(Section~\ref{sec:additivepolynomial}). At each fixed $n$,
positivity of the trace after subtracting the constant mode's
contribution shows that the restriction of the first correction's
distributional extension to $(1,\infty)$ is a locally finite positive
measure supported in $(1,a]$.
Its local mass bound and the detector's exponential decay on a fixed
interval strictly below the transition give an estimate independent of degree
(Section~\ref{sec:exceptionalmeasure}). Combining these estimates
proves Theorems~\ref{thm:sharpfinite} and~\ref{thm:subcritical}
(Section~\ref{sec:additivespectral}).

The table summarizes the parameter ranges. Here $R$ is the geometric
expansion order, $T$ the length cutoff, and $j$ the index of the
coefficient multiplying $g^{1-j}$ in the averaged polynomial orbital
expansion. Canonical coefficient bounds require no genus condition;
estimates for geodesic sums and spectral traces require the stated
bounds. Constants are absolute in the first row and may depend on
$f$ in the other rows; detector and spectral count bounds may also
depend on $\varepsilon$.

\begin{center}
\small
\renewcommand{\arraystretch}{1.15}
\begin{tabular}{@{}>{\raggedright\arraybackslash}p{0.25\linewidth}>{\raggedright\arraybackslash}p{0.25\linewidth}>{\raggedright\arraybackslash}p{0.42\linewidth}@{}}
\toprule
Estimate & Parameters & Required domain \\
\midrule
Orbital remainder & $R,n\ge0$, $T\ge0$;
$\rho=R+\lceil T\rceil+1$ &
$g\ge Cs^2+[C\rho]^Cs$, $s=n+1$ \\
Canonical coefficient norms & $j\ge2$, $n,Q\ge1$ &
No genus condition \\
Exceptional detector response & $n\ge1$, $Q\ge Q_\varepsilon$ &
No genus condition \\
Corrected polynomial transfer & $n,Q\ge1$ &
$g\ge C_{\rm add}s^2+C_{\rm add}Q^{D_3}s$ \\
Fixed-window count & $n\ge1$, $0<\varepsilon<1/4$ &
$g\ge g_\varepsilon$ and $g\ge C_\varepsilon s^2$ \\
\bottomrule
\end{tabular}
\end{center}

Our orbital sums are indexed by primitive \emph{oriented} geodesics,
with their repetitions included separately; the two orientations
are distinct. The separating kernel sums both orders of the two
components, and each pair of boundary components glued to one
another has a single twist coordinate. Weil--Petersson volumes use the intersection-theoretic
orbifold convention
$V_{0,3}=1$ and $V_{1,1}(\ell)=(\ell^2+4\pi^2)/48$.
Together with the Fourier convention \eqref{intro:fourier}, these
normalizations fix the factors in the gluing formulas and spectral trace.

\section{Uniform complementary volumes and boundary kernels}
\label{sec:volume}

The geometric averages in Section~\ref{sec:orbital} involve a chosen
subsurface, called the carrier, and the components glued to its
boundary. Their Weil--Petersson volume ratios form the boundary
kernels constructed here. We expand these kernels in inverse powers
of genus, with explicit bounds in the cusp count, expansion order,
and boundary lengths.

An envelope valid for every stable component controls configurations
with large deficit
(Sections~\ref{sec:volumeenvelope}--\ref{sec:volumetail}). Joint
intersection-number estimates give boundary-volume expansions and
shifted ratios for the large complementary component
(Sections~\ref{sec:volumesource}--\ref{sec:volumeshifts}).
Section~\ref{sec:volumekernels} assembles these inputs into a single
coefficient sequence, independent of genus and truncation order.
Section~\ref{sec:remainders} sharpens these genus thresholds, giving
the additive condition in Section~\ref{sec:additivethreshold}.

For nonnegative integers $u,v$, a type $(u,v)$ is stable when its
complexity $c=2u+v-2$ is at least one; this includes $(0,3)$ and
$(1,1)$. For a stable type, let $V_{u,v}(\boldsymbol\ell)$ denote
the Weil--Petersson orbifold volume of genus-$u$ surfaces with $v$
labelled ends of lengths $\boldsymbol\ell\in[0,\infty)^v$.
A positive length specifies a geodesic boundary and a zero length
specifies a cusp. We write $V_{u,v}=V_{u,v}(0,\ldots,0)$.
All expansion orders are nonnegative integers, and all constants
in this section are absolute. We allow $C$ to increase between
estimates; one sufficiently large choice works throughout,
independently of genus, cusp count, expansion order and carrier type.

\subsection{A volume envelope valid for every stable component}
\label{sec:volumeenvelope}

The following bounds apply to the carrier and every complementary
component, including those of genus zero or one.

\begin{lemma}\label{vol:envelope}
Put $A_0=4\pi^2$.  There are $U\ge1$, $\beta_{\rm v}>0$ and
$\kappa_{\rm v}>0$ such that
\begin{align}
 V_{u,v}&\le U A_0^{c-1}\frac{(c-1)!}{\sqrt c}
 &&(c=2u+v-2\ge1),\label{vol:upper}\\
 V_{g,n}&\ge\beta_{\rm v} A_0^{N}\frac{N!}{\sqrt g}
       \exp\!\left(-\frac{\kappa_{\rm v} n(n+1)}{2g}\right)
 &&(g\ge2,\ n\ge0,\ N=2g+n-3).\label{vol:lower}
\end{align}
Also, for every stable type and every nonnegative length vector,
\begin{equation}\label{vol:lengthbound}
 V_{u,v}(\boldsymbol\ell)\le
 V_{u,v}\exp\!\left(\frac12\sum_j\ell_j\right).
\end{equation}
\end{lemma}

\begin{proof}
We use Mirzakhani's boundary estimate and monotonicity
$V_{u,v}\le V_{u+1,v-2}$ for $v\ge4$ and $2u+v-4>0$
\cite[Lemma 3.2(2)--(3)]{Mir13}, the positive fixed-$v$ large-genus
asymptotic \cite[Theorem 1.2]{MZ15}, and the global marked-point ratio
estimate of Mirzakhani--Zograf \cite[Lemma~5.1]{MZ15}, in the form
recalled in \cite[Lemma~B.6]{HWX26}:
\begin{equation}\label{vol:globalratio}
 \left|\frac{(2g-2+j)V_{g,j}}{V_{g,j+1}}-A_0^{-1}\right|
 \le C_B\frac{j+1}{2g-2+j}\qquad(g\ge2,\ j\ge0).
\end{equation}
The $j+1$ formulation in \cite[Lemma~B.6]{HWX26} includes the
zero-mark case $j=0$.
For $j=0$, the same estimate also follows from the fixed-marked-point
asymptotic, after enlarging $C_B$ to cover the finitely many remaining
genera.

For the upper bound, first absorb the type $(0,4)$ into $U$.
For every other stable type, iterate the monotonicity while $v\ge4$.
It preserves $c$ and
ends in one of the four families $v=0,1,2,3$.  In each family the
fixed-$v$ asymptotic bounds the normalized volume in
\eqref{vol:upper}; the finitely many initial stable terms are absorbed
in $U$.  This includes $(0,3)$ and $(1,1)$, without a step from
$v=3$ to $v=1$.

For the lower bound, set $1+\epsilon_j=A_0(2g-2+j)V_{g,j}/V_{g,j+1}$ and
$\kappa_{\rm v}=A_0C_B$.  Positivity and \eqref{vol:globalratio} give
$\log(1+\epsilon_j)\le\kappa_{\rm v}(j+1)/(2g-2+j)$.
Since $2g-2+j\ge g$, telescoping yields
\[
 V_{g,n}\ge V_{g,0}A_0^n\frac{(2g+n-3)!}{(2g-3)!}
             e^{-\kappa_{\rm v} n(n+1)/(2g)}.
\]
The $n=0$ asymptotic and positivity imply
$\inf_{g\ge2}V_{g,0}\sqrt g/[A_0^{2g-3}(2g-3)!]>0$, proving
\eqref{vol:lower}.  The first cited boundary estimate gives
\eqref{vol:lengthbound}, including zero lengths by continuity.
\end{proof}

\subsection{The labelled complementary tail}
\label{sec:volumetail}

Fix integers $g\ge2$, $n\ge0$ and a carrier type $(h_0,a_0,k)$,
where $h_0$ is the genus, $a_0$ the cusp count, and $k$ the number
of boundary components. Assume $h_0\ge0$, $0\le a_0\le n$,
$k\ge1$ and $m=2h_0+a_0+k-2\ge1$, and label the boundary slots
by $\{1,\ldots,k\}$.

A gluing pattern $P=(M,\Pi)$ consists of a matching $M$ of $2\sigma$
slots into $\sigma$ pairs and a partition $\Pi=(I_1,\ldots,I_b)$
of the remaining slots into $b\ge1$ nonempty blocks, ordered by
least label. A matched pair glues two carrier boundaries to one
another; a block $I_i$ specifies the boundaries attached to the
$i$th complementary component. Thus $k_i=|I_i|$ satisfies
$\sum_i k_i=k-2\sigma$. The complementary genera and cusp counts
are nonnegative integers $(g_i,a_i)$ satisfying
\begin{equation}\label{vol:admissible}
 \sum_i a_i=n-a_0,\qquad
 \sum_i g_i=g-h_0-k+\sigma+b,\qquad
 c_i=2g_i+a_i+k_i-2\ge1.
\end{equation}
An empty admissible sum is zero. Write $N$ for the ambient factorial
index and $R_i$ for the corresponding index of the $i$th complementary
component. The deficit $\delta$ is the loss of index on passing to a
largest complementary component. Define
\begin{gather}
 N=2g+n-3,\qquad R_i=c_i-1,\qquad
 \delta=N-\max_iR_i,\qquad d=m+b-1,\label{vol:deficit}\\
 W=\frac{n!}{a_0!\prod_i a_i!}
       \frac{V_{h_0,a_0+k}\prod_iV_{g_i,a_i+k_i}}{V_{g,n}},\qquad
 \mathcal T_{g,n}^{h_0,a_0,k}(R)
   =\sum_P\ \sum_{\substack{\eqref{vol:admissible}\\\delta\ge R+1}}W.
 \label{vol:tailsum}
\end{gather}
The multinomial in $W$ assigns the $n$ cusp labels to the carrier
and the complementary components, and is included only once. This
labelled sum provides an upper bound before the gluing symmetry
factors are inserted in Section~\ref{sec:orbital}. A matched pair
contributes a gluing edge, rather than an additional annular component.

\begin{theorem}[Complementary-volume tail]\label{vol:tail}
There is an absolute $K\ge4096$ such that, simultaneously for all the
above data and all integers $R\ge0$, with $t=R+m+1$ and $s=n+1$,
\begin{equation}\label{vol:tailestimate}
 g\ge[Kts]^8\quad\Longrightarrow\quad
 \mathcal T_{g,n}^{h_0,a_0,k}(R)
 \le g^{-R-1}[Kts]^{16t}.
\end{equation}
For a fixed pattern $P$, let $E_P$ consist of one edge for each
matched pair and one edge for each unmatched slot, so $|E_P|=k-\sigma$.
Assign a length $y_e\ge0$ to each edge and put
$\mathbf y=(y_e)_{e\in E_P}$.
Write $x_j(\mathbf y)$ for the length of the edge incident to slot $j$;
matched slots therefore have the same length. Put
\begin{equation}\label{vol:edges}
 \Lambda=\sum_{j=1}^k x_j(\mathbf y),\qquad
 J_P(\mathbf y)=\prod_{e\in E_P}y_e.
\end{equation}
In the tail for this pattern, replace the numerator volumes by
boundary volumes with these slot lengths, keeping the cusp lengths
zero. After multiplication by $J_P$, the resulting tail is bounded by
\begin{equation}\label{vol:boundarytail}
 g^{-R-1}[Kt(n+1)]^{16t}(1+\Lambda)^k e^\Lambda
 \qquad\bigl(g\ge[Kt(n+1)]^8\bigr).
\end{equation}
\end{theorem}

\begin{proof}
Euler characteristic and \eqref{vol:admissible} give
\begin{equation}\label{vol:euler}
 \sum_i c_i=N+1-m,\qquad \sum_iR_i=N-d,\qquad \delta\ge d\ge m.
\end{equation}
When $g\ge\max\{m+2,\kappa_{\rm v}s^2\}$, the exponential in
\eqref{vol:lower} is at least $e^{-1/2}$, and a largest component has
$c_*\ge(2g+n-2-m)/b\ge g/b$.  Its square-root denominator therefore
cancels the ambient $\sqrt g$ up to a factor $\sqrt b$.
The total exponent of $A_0$ is $-b$.
Discarding the other square-root denominators in
Lemma~\ref{vol:envelope} gives, with one absolute $B\ge2$,
\begin{equation}\label{vol:factorialenvelope}
 W\le P_b\frac{\prod_iR_i!}{N!}\frac{n!}{a_0!\prod_i a_i!},
 \qquad P_b=B^{b+2}\sqrt b\,(m-1)!.
\end{equation}
This includes $R_i=0$.

Choose a component $*$ maximizing $R_i$, breaking ties by the fixed
block order. The identity
\[
 \delta=m+\sum_{i\ne *}c_i
\]
shows that bounded deficit bounds the total complexity of the carrier
and the other complementary components. The cusps outside component
$*$ include those on the carrier. Since $a_i\le c_i+1$, the number of
cusps outside component $*$ satisfies
\begin{equation}\label{vol:cuspbudget}
 A_{\mathrm{out}}=a_0+\sum_{i\ne *}a_i
 \le a_0+\delta-m+b-1
 =\delta+1-2h_0-k+b\le\delta+1.
\end{equation}
Since $a_*=n-A_{\mathrm{out}}$,
\[
 \frac{n!}{a_0!\prod_i a_i!}
 =\frac{n!}{(n-A_{\mathrm{out}})!\,a_0!\prod_{i\ne *}a_i!}
 \le s^{A_{\mathrm{out}}},
\]
which gives the first bound below; the multinomial theorem gives the second:
\begin{equation}\label{vol:labelbound}
 \frac{n!}{a_0!\prod_i a_i!}
 \le\min\{s^{\delta+1},(b+1)^n\}.
\end{equation}
This bound includes a pair of pants attached along one boundary,
whose other two ends may be cusps.

For a fixed pattern and deficit there are at most $b$ choices of $*$.
The other $R_i$ sum to $\delta-d$ and hence have at most
$(\delta+1)^{b-1}$ choices. For each $i\ne *$, the cusp count $a_i$
has at most $\delta+3$ choices, since
$0\le a_i\le R_i+2\le\delta+2$; the identity
$2g_i=R_i+3-k_i-a_i$ then fixes the genus. The two sum
constraints in \eqref{vol:admissible} then determine the largest
component. Discarding the remaining admissibility and parity
conditions gives at most
\begin{equation}\label{vol:deficitcount}
 b(\delta+3)^{2(b-1)}
\end{equation}
complementary allocations. Also
$\prod_{i\ne *}R_i!\le(\delta-d)!\le\delta!$.
Thus the contribution at this deficit is at most $P_bb x_\delta$, where
\begin{equation}\label{vol:endpointsequence}
 x_\delta=\frac{s^{\delta+1}(\delta+3)^{2b}}{\binom N\delta}.
\end{equation}

We split the deficit sum into a near range, where successive terms
decay geometrically, and a deep range, where a binomial coefficient
absorbs the total number of allocations. Set
$D=\max\{m,R+1\}$ and $J=\lfloor N/(64s)\rfloor$.
By increasing the absolute $K$, the hypothesis in
\eqref{vol:tailestimate} ensures
\[
 g\le N\le3g,\quad D\le t,\quad b\le m+2\le t+1,
 \quad J\ge t,\quad N-D+1\ge g,\quad N\ge128s.
\]
For $D\le\delta<J$, use $b\le\delta+2$ to obtain
\[
 \frac{x_{\delta+1}}{x_\delta}
 =s\frac{\delta+1}{N-\delta}
       \left(1+\frac1{\delta+3}\right)^{2b}
 \le e^2/63<1/2.
\]
Consequently the near tail $D\le\delta\le J$ for this pattern is at most
\begin{equation}\label{vol:neartail}
 2P_bb\,s^{D+1}(D+3)^{2b}D!\,g^{-R-1},
\end{equation}
since $(N)_D\ge g^D$ and $D\ge R+1$, where
$(N)_D=N(N-1)\cdots(N-D+1)$ denotes the falling factorial.

For the deep tail we use the following factorial inequality.
If $N,J$ are integers with $1\le J\le N/2$, $b\ge2$, and
$z_1,\ldots,z_b$ are nonnegative integers satisfying
$\sum_i z_i\le N$ and $z_i\le N-J$, then
\begin{equation}\label{vol:cappedfactorials}
 \prod_i z_i!\le(N-J)!J!.
\end{equation}
Indeed, first increase the coordinates until their sum is $N$, keeping each
at most $N-J$; this is possible since $b(N-J)\ge N$.  Transferring one unit from a smaller positive coordinate
$y$ to a largest coordinate $x<N-J$ multiplies the product by
$(x+1)/y\ge1$.  Continue until that coordinate equals $N-J$; the
remaining mass is $J$ and its factorial product is at most $J!$.
If $\delta>J$, apply this to $R_i$.  For $b=1$ the deep sector is empty,
since its only possible deficit is $m\le J$.  Otherwise the termwise
bound is $P_b(b+1)^n/\binom NJ$.  Each $g_i\le g$, since
$b\le k-2\sigma$, so there are at most $(N+1)^{2b}$ complementary allocations.  The bound $N\ge128s$ gives $J\ge N/(128s)$. Hence,
with $\alpha_0=(\log2)/128$,
\[
 \binom NJ\ge(N/J)^J\ge2^J\ge e^{\alpha_0N/s}.
\]
The deep contribution is therefore at most
\begin{equation}\label{vol:deeptail}
 P_b(b+1)^n(N+1)^{2b}e^{-\alpha_0N/s}\le P_bg^{-R-1}.
\end{equation}
To verify the last inequality uniformly in all parameters, note that
the threshold gives $b+1\le g$, $N+1\le g^2$ and
$R+1+4b+n\le8ts$.  The logarithm to absorb is at most $8ts\log g$.
Since $\log g\le\sqrt g$ and
$\sqrt g\ge(Kts)^4\ge(8/\alpha_0)ts^2$, it is bounded by
$\alpha_0g/s\le\alpha_0N/s$.

It remains to sum the estimates over gluing patterns.
A partial matching is encoded by each slot's partner or zero, giving
at most $(k+1)^k$ choices; a partition is encoded by its least block
labels, giving at most $k^k$ choices.  Thus the number of patterns is
at most $(2k+2)^{2k+2}\le(4t)^{4t}$.  The bounds
\[
 P_b\le(Ct)^{3t},\qquad
 2b\,s^{D+1}(D+3)^{2b}D!+1\le(Ct)^{5t}s^{2t}
\]
follow from $D\le t$, $b\le t+1$ and the factorial bounds
$(m-1)!\le t^{t-2}$, $D!\le t^t$.
The sum of the bounds in \eqref{vol:neartail} and \eqref{vol:deeptail}
over all patterns is at most $g^{-R-1}(Cts)^{12t}$, which proves
\eqref{vol:tailestimate} after enlarging $K$.

For the boundary estimate, fix one pattern. Each gluing edge appears
twice among the boundary lists of the numerator volumes. Their total
boundary length is therefore
$2\sum_{e\in E_P}y_e=2\Lambda-2\sum_{e\in M}y_e$.
Equation~\eqref{vol:lengthbound} contributes at most $e^\Lambda$, while
$J_P\le(1+\Lambda)^k$.  The zero-length tail for this pattern is nonnegative and at most
the full scalar tail in \eqref{vol:tailestimate}. This proves
\eqref{vol:boundarytail} on the pattern's own edge-coordinate space,
including zero lengths by continuity.
\end{proof}

\subsection{Joint intersection-number input and boundary summation}
\label{sec:volumesource}

The tail estimate leaves configurations with bounded deficit and one
large complementary component. We now expand its boundary volume.
The source estimates control both the number of marked ends and the
expansion order. Summing their intersection-number coefficients gives
an expansion for every nonnegative length vector, without a
length-dependent genus threshold.

For $p\ge1$ and $\mathbf d\in\Z_{\ge0}^p$, write
$|\mathbf d|=\sum_i d_i$. We use normalized intersection numbers
$[\tau_{\mathbf d}]_{G,p}$: the coefficient of
$\prod_i z_i^{2d_i}$ in $V_{G,p}(2\mathbf z)$ is
$[\tau_{\mathbf d}]_{G,p}/\prod_i(2d_i+1)!$.
Thus $[\tau_{\mathbf0}]_{G,p}=V_{G,p}$, and the bracket is zero
when $|\mathbf d|>3G+p-3$. For an integer $D\ge0$ and a real
$v\ge0$, put $E_D(v)=\sum_{q=0}^Dv^q/q!$.
For $G\ge2$ and $p\ge0$, define the scalar ratios
\begin{equation}\label{vol:scalarratios}
 H_{G,p}=A_0(2G-2+p)\frac{V_{G,p}}{V_{G,p+1}},\qquad
 P_{G,p}=\frac{V_{G-1,p+2}}{V_{G,p}}.
\end{equation}

\begin{proposition}[Joint source estimates]\label{vol:hma}
There is an absolute integer $c_0>600$ with the following properties.
Put
\[
 D_j(p)=500(c_0j)^{c_0j}(p+j+1)^j\quad(j\ge1),\qquad
 G_J(p)=500(p+J+2)(c_0(J+1))^{c_0}.
\]
There are coefficients $b_{\mathbf d,p}^j,h_p^j,p_p^j$, independent
of genus and truncation order, with $b_{\mathbf d,p}^0=1$ and
constant term one in each scalar expansion, such that for every
integer $J\ge0$ and integer $G>G_J(p)$,
\begin{align}
 \left|\frac{[\tau_{\mathbf d}]_{G,p}}{V_{G,p}}
       -\sum_{j=0}^Jb_{\mathbf d,p}^jG^{-j}\right|
 &\le G^{-J-1}D_{J+1}(p)E_{2J+2}(|\mathbf d|),
       &&p\ge1,\label{vol:intersectioninput}\\
 \left|H_{G,p}-1-\sum_{j=1}^Jh_p^jG^{-j}\right|
 +\left|P_{G,p}-1-\sum_{j=1}^Jp_p^jG^{-j}\right|
 &\le2G^{-J-1}D_{J+1}(p),&&p\ge0.\label{vol:ratioinput}
\end{align}
The first estimate is asserted only for
$|\mathbf d|\le3G+p-3$. The coefficients $b_{\mathbf d,p}^j$ and
their bounds below are defined for every nonnegative multiindex,
including those beyond the dimension cutoff at a given genus.
For every integer $j\ge1$,
\begin{equation}\label{vol:sourcecoefficients}
 |b_{\mathbf d,p}^j|\le D_j(p)E_{2j}(|\mathbf d|),\qquad
 |h_p^j|,|p_p^j|\le D_j(p).
\end{equation}
\end{proposition}

These are slightly weakened forms of
\cite[Theorem~4.2, A1, A3 and A4]{HMT25}, which refine the
fixed-parameter expansions of Mirzakhani--Zograf
\cite[Theorem~4.1]{MZ15}. Choose $c_0$ larger than the constant
required in \cite[Theorem~4.2]{HMT25}. Then $G_J(p)$ dominates all three source
thresholds. The scalar majorants follow directly, and the
multinomial identity
$\sum_{|\alpha|=q}\mathbf d^\alpha/\alpha!=|\mathbf d|^q/q!$
sums the intersection-number majorants to the factors $E_D$ above.

Mirzakhani's volume-polynomial and intersection-number formula
\cite[Theorems~1.1 and~4.4]{Mir07Intersection}, in the normalization recalled in
\cite[Theorem~4.4]{HMT25}, gives
\begin{equation}\label{vol:volumepolynomial}
 \frac{V_{G,p}(\mathbf0_{p-\ell},\mathbf x)}{V_{G,p}}
 =\sum_{|\mathbf d|\le3G+p-3}
   \frac{[\tau_{(\mathbf0,\mathbf d)}]_{G,p}}{V_{G,p}}
   \prod_{i=1}^{\ell}\frac{(x_i/2)^{2d_i}}{(2d_i+1)!}.
\end{equation}
The formula extends to nonnegative lengths by continuity.  Throughout,
we use the orbifold convention
$V_{1,1}(x)=(x^2+4\pi^2)/48$ of \cite[\S1]{DN09}; in particular,
$V_{1,1}=\pi^2/12$. The convention in \cite[\S2.4]{HMT25}
instead gives $V_{1,1}=\pi^2/6$.
All applications of \cite[Theorem~4.2]{HMT25} above have genera
$G,G-1\ge3$, where the generic stabilizer is trivial and the two
conventions agree.  Thus the imported normalized intersection numbers
and scalar ratios are unchanged.  Terminal components below always
use our declared orbifold volumes.
The estimates above imply, with one absolute $C$, the coarser bounds
$D_j(p)\le[Cj(p+1)]^{Cj}$ and
$G_J(p)\le[C(J+1)(p+1)]^C$. These bounds hold uniformly in both
the expansion order and the number of marked ends.

\begin{lemma}[Boundary-volume expansion]\label{vol:boundaryexp}
For $0\le\ell\le p$, $p\ge1$, there are continuous functions
$B_{j,p}^{(\ell)}$ on $[0,\infty)^\ell$, independent of genus and
truncation, such that, writing $X=\sum_i x_i$,
\begin{align}
 B_{0,p}^{(\ell)}(\mathbf x)&=
 \Phi_\ell(\mathbf x):=
       \prod_{i=1}^{\ell}\frac{\sinh(x_i/2)}{x_i/2},\label{vol:boundaryleading}\\
 |B_{j,p}^{(\ell)}(\mathbf x)|&\le
 [Cj(p+1)]^{Cj}(1+X)^{2j}e^{X/2}\quad(j\ge1),\label{vol:boundarycoefficients}\\
 \left|\frac{V_{G,p}(\mathbf0_{p-\ell},\mathbf x)}{V_{G,p}}
       -\sum_{j=0}^RG^{-j}B_{j,p}^{(\ell)}(\mathbf x)\right|
 &\le G^{-R-1}[C(R+1)(p+1)]^{C(R+1)}
           (1+X)^{2R+2}e^{X/2},\label{vol:boundaryremainder}
\end{align}
for every $R\ge0$ and $G\ge[C(R+1)(p+1)]^C$, with no condition on $X$.
Each factor $\sinh(x_i/2)/(x_i/2)$ is interpreted as one at
$x_i=0$, and an empty product is one.
\end{lemma}

\begin{proof}
Write $z_i=x_i/2$, $Z=\sum_i z_i$ and
$w_{\mathbf d}(\mathbf z)=\prod_i z_i^{2d_i}/(2d_i+1)!$.
The dimension-free summation bound is
\begin{equation}\label{vol:poissonsum}
 \sum_{\mathbf d\ge0}w_{\mathbf d}(\mathbf z)E_D(|\mathbf d|)
 \le(D+1)e^Z(1+Z)^D.
\end{equation}
To prove it, replace $(2d_i+1)!$ by $(2d_i)!$ and enlarge the even
indices $2d_i$ to all nonnegative indices.  Since $|\mathbf d|$ is at
most their sum, the result is bounded by $e^Z\E E_D(Y)$ for a Poisson
variable $Y$ of mean $Z$.  Its factorial moments give
$\E Y^q=\sum_{j=0}^q S(q,j)Z^j$, where $S(q,j)$ is the Stirling
number of the second kind. It counts partitions of a $q$-element
set into $j$ blocks and is at most $q!$
(encode a partition by its canonically ordered cycles), so
$\E Y^q/q!\le\sum_{j=0}^qZ^j\le(1+Z)^q$ for $q\ge1$.
The case $q=0$ is immediate.  Summing proves \eqref{vol:poissonsum}.

Define for every $j\ge0$
\begin{equation}\label{vol:Bdefinition}
 B_{j,p}^{(\ell)}(\mathbf x)=
 \sum_{\mathbf d\ge0}b_{(\mathbf0,\mathbf d),p}^j
                          w_{\mathbf d}(\mathbf x/2).
\end{equation}
The source coefficient bound and \eqref{vol:poissonsum} prove local
uniform absolute convergence, continuity and
\eqref{vol:boundarycoefficients}; $b^0=1$ gives
\eqref{vol:boundaryleading}.
For the remainder, first sum \eqref{vol:intersectioninput} over the
finite volume polynomial \eqref{vol:volumepolynomial};
\eqref{vol:poissonsum} gives the required bound for this part.
The infinite series defining $B_{j,p}^{(\ell)}$ also contain terms
beyond that polynomial's dimension cutoff. In this second part,
$v=|\mathbf d|\ge3G+p-2\ge G$.
For every $0\le j\le R$,
\begin{equation}\label{vol:degreetail}
 G^{R+1-j}\le v^{R+1-j},\qquad
 v^{R+1-j}E_{2j}(v)\le(2R+2)!E_{2R+2}(v).
\end{equation}
The second inequality follows term by term: the polynomial on the left
has degree at most
$R+j+1\le2R+1$, and the factorial quotient is at most $(2R+2)!$.
For $j=0$ use $b^0=1$ and $E_0=1$.  Equation~\eqref{vol:poissonsum}
then bounds the entire degree tail, including its constant-coefficient
part.  Absorbing $(2R+2)!$ and the $R+1$ summands into the displayed
order factor proves \eqref{vol:boundaryremainder}.
\end{proof}

\subsection{Finite products and shifted ratios}
\label{sec:volumeshifts}

To compare a large complementary component with the ambient surface,
we must multiply finitely many expansions and change their genus
parameter. We record the elementary estimates used for both operations.

\begin{lemma}[Quantitative algebra]\label{vol:finiteproduct}
Let $J\ge0$ and $M\ge1$ be integers. Suppose $M$ factors have expansions
$f_\nu(g)=1+\sum_{j=1}^Ja_{\nu,j}g^{-j}+r_\nu(g)$ with
$|a_{\nu,j}|\le A^j$ and $|r_\nu(g)|\le A^{J+1}g^{-J-1}$,
where $A\ge1$. At values of $g$ for which these expansions are valid
and $g\ge4MA$, their product has coefficients bounded
in order $j$ by $(MA)^j$ and remainder bounded by
$C(MA)^{J+1}g^{-J-1}$. Separately, a single expansion in
$(g-H)^{-1}$, $H\ge0$, with coefficient bounds $A^j$ and
remainder bound $A^{J+1}(g-H)^{-J-1}$ can be rewritten in
$g^{-1}$ with $A$ replaced by $C(A+H+1)$, provided
$g\ge4(H+1)$ and the original expansion is valid at $g-H$.
If the input coefficients are consistent across
truncation orders, both operations preserve this consistency.
\end{lemma}

\begin{proof}
For each total degree $j\ge1$, the number of weak compositions into
$M$ parts is at most $M^j$, since every such composition is the
multiplicity vector of a word of length $j$ on $M$ symbols.
This proves the coefficient bound and bounds the discarded polynomial
product by $2(MA/g)^{J+1}$.  Each truncated factor has modulus at most $1+2A/g$, and each exact
factor has modulus at most $1+3A/g$.  Replacing factors one at a time bounds
the sum of remainder contributions by
$MA^{J+1}g^{-J-1}e^{3MA/g}$, as required.

For shifting, use, for $j\ge1$,
\begin{equation}\label{vol:shiftidentity}
 (g-H)^{-j}=g^{-j}\sum_{v\ge0}
                    \binom{j+v-1}{v}(H/g)^v.
\end{equation}
The coefficient of $g^{-q}$ is
$\sum_{j=1}^qa_j\binom{q-1}{j-1}H^{q-j}$, bounded by
$[2(A+H+1)]^q$.  With $H/g\le1/4$, use
$\binom{j+v-1}{v}\le2^{j+v}$ and sum the geometric tail for
$v\ge J+1-j$.  Summing over $j\le J$ gives
$g^{-J-1}[C(A+H+1)]^{J+1}$.  The original remainder costs at most
$2^{J+1}A^{J+1}g^{-J-1}$.  The coefficient formulas involve only
orders up to $q$, proving consistency.
\end{proof}

Let $n,H\ge0$ and $\delta\ge1$ be integers, and set $G=g-H$ and
$p=n+2H-\delta\ge1$. The shifted ratio can be written entirely in
terms of the scalar ratios \eqref{vol:scalarratios}. Telescoping gives
\begin{align}
 \frac{V_{g-H,p}}{V_{g,n}}
 ={}&(8\pi^2g)^{-\delta}
       \prod_{j=0}^{\delta-1}H_{G,p+j}
       \prod_{u=0}^{H-1}P_{G+u+1,n+2H-2u-2}
       \prod_{j=0}^{\delta-1}
       \left(1+\frac{n-\delta+j-2}{2g}\right)^{-1}.
 \label{vol:telescoping}
\end{align}
Indeed the prefactor together with the first and last products gives
$V_{G,p}/V_{G,p+\delta}$, and the middle product gives
$V_{G,n+2H}/V_{g,n}$.
All marked-point indices in the middle product are nonnegative;
all genera along this path are large under the thresholds below.

To construct the expansion through order $J$, put $L_0=J+H+\delta+1$.
Proposition~\ref{vol:hma}, genus shifts and the geometric expansions
of the last factors in \eqref{vol:telescoping} have a common
coefficient majorant $A=[C_1L_0(n+1)]^{C_1}$, valid when $g\ge A$,
for an absolute $C_1$. There are $M=H+2\delta\le2L_0$ normalized
factors. Choosing the absolute constant $C$ in the genus threshold
$g\ge[CL_0(n+1)]^C$ sufficiently large ensures $g\ge4MA$.
Lemma~\ref{vol:finiteproduct} gives fixed formal coefficients
$v_j(n,H,\delta)$, $v_0=1$, such that
\begin{equation}\label{vol:zeroshiftedratio}
 \left|\frac{V_{g-H,p}}{V_{g,n}}
 -(8\pi^2)^{-\delta}g^{-\delta}
          \sum_{j=0}^Jv_j(n,H,\delta)g^{-j}\right|
 \le g^{-\delta-J-1}[CL_0(n+1)]^{C(J+1)}.
\end{equation}
For the coefficient of order $j$, one can take the majorant
$[C(j+H+\delta+1)(n+1)]^{Cj}$.

\begin{lemma}[Shifted boundary ratio]\label{vol:shiftedboundary}
For the preceding $n,H,\delta,p$ and $0\le\ell\le p$ there are
continuous functions $Q_j^{n,H,\delta;\ell}$, $j\ge0$, on
$[0,\infty)^\ell$, independent of genus and truncation order.
With $X=\sum_i x_i$, they satisfy $Q_0=\Phi_\ell$ and
\begin{equation}\label{vol:Qbound}
 |Q_j^{n,H,\delta;\ell}(\mathbf x)|
 \le[C(j+H+\delta+1)(n+1)]^{C(j+1)}
              (1+X)^{C(j+1)}e^{X/2}.
\end{equation}
For every $R\ge\delta$, $L_0=R+H+\delta+1$ and
$g\ge[CL_0(n+1)]^C$,
\begin{align}
 \frac{V_{g-H,p}(\mathbf0_{p-\ell},\mathbf x)}{V_{g,n}}
 &=(8\pi^2)^{-\delta}g^{-\delta}
       \sum_{j=0}^{R-\delta}Q_j^{n,H,\delta;\ell}(\mathbf x)g^{-j}
       +\mathcal E_R,\label{vol:shiftedboundaryformula}\\
 |\mathcal E_R|&\le g^{-R-1}[CL_0(n+1)]^{C(R+1)}
                     (1+X)^{C(R+1)}e^{X/2}.
 \label{vol:shiftedboundaryerror}
\end{align}
\end{lemma}

\begin{proof}
Shift the boundary coefficients of Lemma~\ref{vol:boundaryexp}:
$\widetilde B_0=B_{0,p}^{(\ell)}$ and, for $q\ge1$,
\[
 \widetilde B_q=\sum_{j=1}^q
        \binom{q-1}{j-1}H^{q-j}B_{j,p}^{(\ell)}.
\]
Define $Q_j=\sum_{u+v=j}v_u(n,H,\delta)\widetilde B_v$.
These finite convolution formulas prove consistency and
\eqref{vol:Qbound}. For the remainder put $J=R-\delta$ and $G=g-H$.
Write
\[
 \begin{aligned}
 Z(g)&=(8\pi^2)^\delta g^\delta\frac{V_{G,p}}{V_{g,n}},
 & S_J(g)&=\sum_{u=0}^Jv_u(n,H,\delta)g^{-u},\\
 D(g,\mathbf x)&=\frac{V_{G,p}(\mathbf0_{p-\ell},\mathbf x)}{V_{G,p}},
 & D_J(g,\mathbf x)&=\sum_{v=0}^J\widetilde B_v(\mathbf x)g^{-v}.
 \end{aligned}
\]
Enlarging the constant in the genus threshold gives $|S_J(g)|\le2$;
\eqref{vol:lengthbound} gives $|D(g,\mathbf x)|\le e^{X/2}$.
Use the exact identity
\[
 ZD-S_JD_J=(Z-S_J)D+S_J(D-D_J).
\]
After multiplication by $(8\pi^2)^{-\delta}g^{-\delta}$, the first
term is controlled by \eqref{vol:zeroshiftedratio} and the second by
the shifted boundary remainder. Each bound contains at most one
factor $e^{X/2}$.
In the product $S_JD_J$, only the boundary factor depends
on $\mathbf x$, with length growth at most $e^{X/2}(1+X)^{2J}$.
Extracting this common factor from the discarded terms of total
degree greater than $J$ leaves scalar coefficient majorants, to which
the geometric estimate in Lemma~\ref{vol:finiteproduct} applies. This
bounds all discarded terms by \eqref{vol:shiftedboundaryerror},
without a length-dependent genus threshold.
\end{proof}

\subsection{Canonical retained kernels}
\label{sec:volumekernels}

We now combine the complementary-volume tail with the shifted
boundary expansion. For a fixed carrier type and pattern $P$, let
$\mathbf x=(x_j(\mathbf y))_{j=1}^k$ be the slot lengths induced by
the edge lengths $\mathbf y$. The full nonnegative kernel is
\begin{equation}\label{vol:kerneldefinition}
 F_{g,n}^{h_0,a_0,k;P}(\mathbf y)=
 J_P(\mathbf y)V_{h_0,a_0+k}(\mathbf0_{a_0},\mathbf x)
 \sum_{\eqref{vol:admissible}}
 \frac{n!}{a_0!\prod_i a_i!}
 \frac{\prod_iV_{g_i,a_i+k_i}(\mathbf0_{a_i},\mathbf x_{I_i})}{V_{g,n}}.
\end{equation}

\begin{theorem}[Retained boundary kernels]\label{vol:kernel}
There is an absolute $C\ge2$ and, for every fixed $n,h_0,a_0,k,P$,
a sequence of continuous functions $\mathcal F_{j,n}^{h_0,a_0,k;P}$, $j\ge1$,
independent of genus and truncation order, such that, simultaneously
for all $R\ge0$, $t=R+m+1$, and $\mathbf y\ge0$,
\begin{equation}\label{vol:kernelexpansion}
 \left|F_{g,n}^{h_0,a_0,k;P}(\mathbf y)
          -\sum_{j=1}^Rg^{-j}\mathcal F_{j,n}^{h_0,a_0,k;P}(\mathbf y)\right|
 \le g^{-R-1}[Ct(n+1)]^{Ct}(1+\Lambda)^{Ct}e^\Lambda
\end{equation}
when $g\ge[Ct(n+1)]^C$.  Moreover,
$\mathcal F_{j,n}^{h_0,a_0,k;P}=0$ for $j<m+b-1$.  For every $j\ge1$,
\begin{equation}\label{vol:kernelcoefficientbound}
 |\mathcal F_{j,n}^{h_0,a_0,k;P}(\mathbf y)|
 \le[C(j+m+1)(n+1)]^{C(j+m+1)}
             (1+\Lambda)^{C(j+m+1)}e^\Lambda.
\end{equation}
\end{theorem}

\begin{proof}
Fix the carrier type and pattern. To define coefficients independently
of genus, we record all complementary components except the one that
will become large. A record consists of a distinguished block $*$
and nonnegative integers $(g_i,a_i)_{i\ne *}$ for which the
corresponding components are stable. The distinguished component
will be reconstructed from the ambient genus and the sum constraints
in \eqref{vol:admissible}. For such a record set
\begin{align}
 \delta&=m+\sum_{i\ne *}c_i,&
 A_{\mathrm{out}}&=a_0+\sum_{i\ne *}a_i,\notag\\
 H&=h_0+k-\sigma-b+\sum_{i\ne *}g_i,&
 p&=n-A_{\mathrm{out}}+k_*.
 \label{vol:recordparameters}
\end{align}
For each integer $\delta\ge1$, let $\mathscr E_\delta$ be the set
of records with this deficit and $A_{\mathrm{out}}\le n$.
These sets are independent of genus and truncation order.

Now fix $R\ge0$. Theorem~\ref{vol:tail} controls the configurations
with deficit at least $R+1$. For a record of deficit $\delta\le R$,
direct substitution gives
\begin{equation}\label{vol:recordbounds}
 2H+A_{\mathrm{out}}-k_*=\delta,\quad
 0\le H,\quad A_{\mathrm{out}}\le\delta+1,\quad
 2H\le\delta+k_*\le R+m+2,\quad p\le n+2t.
\end{equation}
The distinguished component therefore has genus $G=g-H$ and
$p=n-A_{\mathrm{out}}+k_*$ marked ends: $n-A_{\mathrm{out}}$
cusps and $k_*$ gluing boundaries. Also,
$p=n+2H-\delta\ge k_*\ge1$, as required by
Lemma~\ref{vol:shiftedboundary}.
The common threshold ensures $g>H$ and $N>2R$, so each such record
reconstructs an admissible large component. Its factorial index is
$2(g-H)+p-3=N-\delta$. Since $\delta\le R$, this exceeds
$\delta-m$, the total complexity of the other complementary
components, and makes it the unique largest component.
Conversely, any admissible
configuration with deficit at most $R$ has a unique largest component
and determines one of these records. Moreover,
\begin{equation}\label{vol:recordcount}
 \mathscr E_\delta=\varnothing\quad(\delta<m+b-1),\qquad
 |\mathscr E_\delta|\le b(\delta+3)^{2(b-1)},
\end{equation}
by the same independent choices of small indices and cusp counts
used in \eqref{vol:deficitcount}.

For a record $e$, let $X_* =\sum_{j\in I_*}x_j$ and put
\begin{equation}\label{vol:recordprefactor}
 A_e(\mathbf y)=J_PV_{h_0,a_0+k}(\mathbf0_{a_0},\mathbf x)
 \frac{n!}{a_0!(n-A_{\mathrm{out}})!\prod_{i\ne *}a_i!}
 \prod_{i\ne *}V_{g_i,a_i+k_i}(\mathbf0_{a_i},\mathbf x_{I_i}).
\end{equation}
The carrier and small components have total complexity $\delta$;
their factorial indices sum to $\delta-b$.  Lemma~\ref{vol:envelope}
and the product-factorial inequality bound their zero-length product
by $U^bA_0^{\delta-b}(\delta-b)!$.
The label factor is at most $(n+1)^{\delta+1}$.
Their total boundary length is at most $2\Lambda-X_*$, including
self-pairs only in the carrier.  Hence
\begin{equation}\label{vol:prefactorbound}
 |A_e(\mathbf y)|\le[C(\delta+1)(n+1)]^{C(\delta+1)}
                (1+\Lambda)^k e^{\Lambda-X_*/2}.
\end{equation}

Using these record sets, define the canonical coefficient sequence:
\begin{equation}\label{vol:canonicalcoefficients}
 \mathcal F_{j,n}^{h_0,a_0,k;P}(\mathbf y)=
 \sum_{\delta=1}^{j}\ \sum_{e\in\mathscr E_\delta}
 (8\pi^2)^{-\delta}A_e(\mathbf y)
 Q_{j-\delta}^{n,H_e,\delta;k_{*,e}}
                   (\mathbf x_{I_{*,e}}).
\end{equation}
These sums are finite and continuous; \eqref{vol:recordcount} proves
the stated vanishing.  For every retained record apply
Lemma~\ref{vol:shiftedboundary} through order $R$.
Its $L_0=R+H_e+\delta+1\le3t$, so the same threshold
$g\ge[Ct(n+1)]^C$ suffices for all records.
Multiplication by \eqref{vol:prefactorbound} changes $e^{X_*/2}$
to $e^\Lambda$.  Since $k\le m+2\le t+1$, all length powers are
bounded by $(1+\Lambda)^{Ct}$.
The number of records is at most $(Ct)^{Ct}$ by
\eqref{vol:recordcount}.  Summing their errors and adding
\eqref{vol:boundarytail} proves \eqref{vol:kernelexpansion}.
If $R<m+b-1$ there are no retained records and the tail alone
proves it, including $R=0$.
The same finite sum, \eqref{vol:Qbound} and
\eqref{vol:prefactorbound} with $R=j$ give
\eqref{vol:kernelcoefficientbound}.
\end{proof}

At the first possible order $d=m+b-1$, the coefficient is
\begin{equation}\label{vol:kernelleading}
 \mathcal F_{d,n}^{h_0,a_0,k;P}(\mathbf y)
 =(8\pi^2)^{-d}\sum_{e\in\mathscr E_d}A_e(\mathbf y)
       \prod_{q\in I_{*,e}}\frac{\sinh(x_q/2)}{x_q/2}.
\end{equation}
Each small component contributing here has complexity one, hence
type $(0,3)$ or $(1,1)$, where the marked ends may be cusps or gluing
boundaries. In particular, terminal pairs of pants already contribute
at this leading order. With a single complementary component,
$d=m$. Section~\ref{sec:orbital} treats the unstable cylinder directly,
without assigning a volume to $(0,2)$.

\section{Orbital integration and cusp-polynomial coefficients}
\label{sec:orbital}

Throughout this section $g\ge3$ and $n\ge0$ are integers, and every cusp is labelled.
For a real bounded even Borel function $H$, supported in $[-T,T]$ with
$T\ge0$, define
\begin{equation}\label{orb:statistic}
 \mathcal O_H(X)=\sum_{\gamma\in\mathcal P_{\rm or}(X)}
 \sum_{\nu\ge1}\frac{\ell_\gamma H(\nu\ell_\gamma)}
 {2\sinh(\nu\ell_\gamma/2)}.
\end{equation}
Here $\mathcal P_{\rm or}$ denotes primitive oriented geodesics; thus both
orientations and all positive repetitions occur.  The primitive underlying
geodesic determines the simple and nonsimple parts of this statistic.
We expand the expectation in powers of $g^{-1}$ and determine the
dependence of its coefficients on $n$. We first treat simple
geodesics by cutting along them, then treat nonsimple geodesics by
integrating over the subsurfaces they fill. Combining these exact
formulas with Section~\ref{sec:volume} gives the canonical coefficients
in Theorem~\ref{orb:expansion} and the geometric moments in
Proposition~\ref{orb:cgm}. Sections~\ref{sec:firstcorrection}
and~\ref{sec:cusppolynomiality} compute the first correction and prove
that the coefficient of $g^{-j}$ has degree at most $j+1$ in $n$.
This polynomiality concerns the coefficients; the sharper uniform
remainder needed when $n$ grows with $g$ is proved in
Section~\ref{sec:remainders}.

\subsection{The two simple-curve kernels}
\label{sec:simplekernels}

Put $\phi(x)=\sinh(x/2)/(x/2)$, with $\phi(0)=1$, and define
\begin{align}
 U_{g,n}(\ell)&=\frac{V_{g-1,n+2}(\mathbf0_n,\ell,\ell)}{V_{g,n}},
 \label{orb:simplekernels}\\
 S_{g,n}(\ell)&=
 \sum_{\substack{0\le i\le g,\ 0\le a_0\le n\\
                   2i+a_0\ge2,\ 2(g-i)+n-a_0\ge2}}
 \binom n{a_0}
 \frac{V_{i,a_0+1}(\mathbf0_{a_0},\ell)
 V_{g-i,n-a_0+1}(\mathbf0_{n-a_0},\ell)}{V_{g,n}}.
 \label{orb:sepkernel}
\end{align}
The kernel $U_{g,n}$ corresponds to a nonseparating cut, and $S_{g,n}$
to a separating cut with a distinguished side. Thus the separating
sum is ordered; its factor $\binom n{a_0}$ chooses the cusp labels on
that side. All volume factors are stable. The orientation and
side-counting factors are verified in Lemma~\ref{orb:simpleintegration}.

\begin{lemma}[Simple-kernel expansion]\label{orb:simplekernel}
There is an absolute $C\ge2$ and, for each $n$, a fixed sequence of continuous
functions $u_{j,n},s_{j,n}$ such that, for every integer $R\ge0$ and
$r=R+1$,
\begin{align}
 &\left|U_{g,n}(\ell)-\phi(\ell)^2-\sum_{j=1}^R g^{-j}u_{j,n}(\ell)\right|
 +\left|S_{g,n}(\ell)-\sum_{j=1}^R g^{-j}s_{j,n}(\ell)\right|
 \notag\\[-2pt]
 &\hspace{15mm}\le
 g^{-R-1}[Cr(n+1)]^{Cr}(1+\ell)^{Cr}e^\ell,
 \qquad g\ge[Cr(n+1)]^C,\quad \ell\ge0.
 \label{orb:simplekernelbound}
\end{align}
The coefficients are independent of genus and truncation order, and satisfy
\begin{equation}\label{orb:simplecoefficientbound}
 |u_{j,n}(\ell)|+|s_{j,n}(\ell)|
 \le[C(j+1)(n+1)]^{C(j+1)}(1+\ell)^{C(j+1)}e^\ell.
\end{equation}
In particular,
\begin{equation}\label{orb:firstseparating}
 s_{1,n}(\ell)=\frac{\phi(\ell)}{4\pi^2}
 \left[\frac{\ell^2+4\pi^2}{48}+\binom n2\right].
\end{equation}
\end{lemma}

\begin{proof}
We expand the nonseparating kernel first. For the separating kernel,
we bound the large-deficit tail and then expand the finitely many
retained endpoint types. The handle-ratio coefficients of Proposition~\ref{vol:hma}
are denoted by $p_n^j$, with $p_n^0=1$.  Factor
\begin{equation}\label{orb:nonsepfactor}
 U_{g,n}(\ell)=\frac{V_{g-1,n+2}}{V_{g,n}}
 \frac{V_{g-1,n+2}(\mathbf0_n,\ell,\ell)}{V_{g-1,n+2}}.
\end{equation}
Let $B_{j,p}^{(b)}$ be the fixed boundary coefficients of
Lemma~\ref{vol:boundaryexp}.  The binomial expansion of $(g-1)^{-j}$ gives
fixed coefficients
\begin{equation}\label{orb:ucoefficient}
 \beta_{0,n}=\phi^2,\qquad
 \beta_{d,n}(\ell)=\sum_{j=1}^d\binom{d-1}{j-1}
 B_{j,n+2}^{(2)}(\ell,\ell),\qquad
 u_{d,n}=\sum_{a+b=d}p_n^a\beta_{b,n}\quad(d\ge1).
\end{equation}
The shifted boundary remainder is bounded by
$g^{-R-1}[Cr(n+1)]^{Cr}(1+2\ell)^{2R+2}e^\ell$.
To multiply the two expansions, write their error as
$(P-P^{[R]})D+P^{[R]}(D-D^{[R]})$, where $P,D$ are the two factors
in \eqref{orb:nonsepfactor} and the superscript $[R]$ denotes their
truncation through order $g^{-R}$. Discard the terms of degree greater
than $R$ in the product of these truncations.  Here $|D|\le e^\ell$, while the scalar
truncation $P^{[R]}$ is bounded at the stated genus threshold.  Every discarded
product has only one boundary coefficient, so its common length weight is
factored out once.  The geometric majorant of Lemma~\ref{vol:finiteproduct}
then proves the required bound.  This argument does not require
$\ell/g$ to be small and also applies when $R=0$.

For the separating tail write $N=2g+n-3$, $s=n+1$, and, in a summand of
\eqref{orb:sepkernel}, put
\begin{equation}\label{orb:sepdeficit}
 c_1=2i+a_0-1,\quad c_2=2(g-i)+n-a_0-1,\quad
 c_1+c_2=N+1,\quad \delta=\min(c_1,c_2)\ge1.
\end{equation}
The volume envelope of Lemma~\ref{vol:envelope}, at $g\ge Cs^2$,
gives a universal $B$ such that
\begin{equation}\label{orb:septerm}
 \binom n{a_0}\frac{V_{i,a_0+1}V_{g-i,n-a_0+1}}{V_{g,n}}
 \le B\binom n{a_0}\frac{(\delta-1)!(N-\delta)!}{N!}
 =\frac{B\binom n{a_0}}{\delta\binom N\delta}.
\end{equation}
Indeed, the powers of $4\pi^2$ leave a fixed factor, the exponential
denominator correction is bounded, and the larger complexity is at least
$g/2$, controlling the square-root factor.  A smaller side has genus $h_0$, $a_*$ cusps and one boundary,
with $2h_0+a_*-1=\delta$, so $a_*\le\delta+1$ and
$\binom n{a_0}=\binom n{a_*}\le s^{\delta+1}$.  There are at most
$2(\delta+2)$ ordered records at this deficit, including a possible tie.
Consequently their sum is at most $6B x_\delta$, where
$x_\delta=s^{\delta+1}/\binom N\delta$.

For $r=R+1$ and $g\ge(Krs)^8$, take $D=\lfloor N/(64s)\rfloor$.
An absolute $K$ ensures $r\le D$, $g\le N\le3g$, and $N-r+1\ge g$.
For $r\le\delta<D$,
\begin{equation}\label{orb:sepratio}
 \frac{x_{\delta+1}}{x_\delta}
 =\frac{s(\delta+1)}{N-\delta}\le\frac1{63},
 \qquad
 \sum_{\delta=r}^D6Bx_\delta
 \le12B s^{r+1}r!g^{-r}.
\end{equation}
For $\delta>D$, use binomial unimodality in
\eqref{orb:septerm}, the identity $\sum_{a_0}\binom n{a_0}=2^n$,
and the at most $g+1$ genus choices. Since
$D\ge N/(128s)$ and $\binom ND\ge2^D$, this part is bounded by
$B(g+1)2^n e^{-\alpha_0N/s}$, where $\alpha_0=\log2/128$.
The same threshold, with $K$ enlarged, gives
\begin{equation}\label{orb:sepdeep}
 \log\bigl(g^r(g+1)2^n\bigr)
 \le3rs\log g\le3rs\sqrt g\le\alpha_0N/s.
\end{equation}
Thus the deep tail is at most $Bg^{-r}$.  Boundary growth adds at most
$e^\ell$ for the two factors, proving the direct estimate
\begin{equation}\label{orb:septail}
 \sum_{\delta\ge R+1}\binom n{a_0}
 \frac{V_{i,a_0+1}(\mathbf0_{a_0},\ell)
 V_{g-i,n-a_0+1}(\mathbf0_{n-a_0},\ell)}{V_{g,n}}
 \le g^{-R-1}[K(R+1)(n+1)]^{8(R+1)}e^\ell.
\end{equation}
The estimate also covers a pair of pants on the smaller side,
equal-complexity sides, and the closed case $n=0$.

It remains to expand the retained endpoints.  Define the fixed finite set
\begin{equation}\label{orb:endpointset}
 \mathscr E_\delta(n)=\{(h_0,a_0)\in\mathbb Z_{\ge0}^2:a_0\le n,
                         \ 2h_0+a_0-1=\delta\}.
\end{equation}
When $N>2R$, a term with $\delta\le R$ has a unique smaller side.
Each record in \eqref{orb:endpointset} has exactly two placements in the
ordered sum.  Its other component has genus $g-h_0$ and
$p=n-a_0+1=n+2h_0-\delta\ge1$ marked points.  The shifted ratio of
Lemma~\ref{vol:shiftedboundary} has canonical coefficients
$Q_d^{n,h_0,\delta;1}(\ell)$, starting with $Q_0=\phi$, and leading scalar
$(8\pi^2)^{-\delta}g^{-\delta}$.  Therefore define
\begin{equation}\label{orb:scoefficient}
 s_{j,n}(\ell)=2\sum_{\delta=1}^j
 \sum_{(h_0,a_0)\in\mathscr E_\delta(n)}
 (8\pi^2)^{-\delta}\binom n{a_0}
 V_{h_0,a_0+1}(\mathbf0_{a_0},\ell)
 Q_{j-\delta}^{n,h_0,\delta;1}(\ell).
\end{equation}
Here $h_0\le(\delta+1)/2$, $a_0\le\delta+1$, and
$|\mathscr E_\delta(n)|\le\delta+2$.  The small-volume prefactor is at
most $[C(\delta+1)(n+1)]^{C(\delta+1)}e^{\ell/2}$ by
Lemma~\ref{vol:envelope}.  Multiplication by the shifted remainder
contributes $e^{\ell/2}$ once more.  Summing the $O(R^2)$ records and
adding \eqref{orb:septail} proves \eqref{orb:simplekernelbound}; the same
finite sums at order $j$ prove \eqref{orb:simplecoefficientbound}.
For $\delta=1$ the only records are $(1,0)$ and, if available, $(0,2)$.
Using $V_{1,1}(\ell)=(\ell^2+4\pi^2)/48$ and $V_{0,3}=1$ gives
\eqref{orb:firstseparating}.
\end{proof}

\subsection{Exact simple integration and repetitions}

Because \eqref{orb:statistic} counts both orientations, we convert
the unoriented cutting formula to an oriented density.

\begin{lemma}[Oriented simple integration]\label{orb:simpleintegration}
Writing $K_{g,n}=U_{g,n}+S_{g,n}$, one has
\begin{equation}\label{orb:simpledensity}
 \E_{g,n}\sum_{\gamma\in\mathcal P_{\rm or}^{\rm simp}(X)}F(\ell_\gamma)
 =\int_0^\infty F(\ell)\ell K_{g,n}(\ell)\,d\ell
\end{equation}
for every nonnegative Borel $F$, with infinity allowed.  Consequently
$\mathcal O_H^{\rm simp}$ is absolutely integrable and
\begin{equation}\label{orb:simpleexact}
 \E_{g,n}\mathcal O_H^{\rm simp}
 =\int_0^\infty\sum_{\nu\ge1}
 \frac{\ell^2H(\nu\ell)}{2\sinh(\nu\ell/2)}K_{g,n}(\ell)\,d\ell.
\end{equation}
\end{lemma}

\begin{proof}
For an unoriented simple orbit, the marked integration formula
\cite[\S4, (4.1), Theorem 4.1 and (4.3)]{Mir08} gives the density
$2^{-M}N(\gamma)^{-1}\ell\prod_Z\widetilde V_Z$ before division by
$V_{g,n}$.  Here $N(\gamma)$ is the index of the subgroup preserving the curve's
orientation in its full stabilizer, and $M$ counts intact one-holed
tori separated off by the curve.  In genus at least three there is at
most one such component.  Its geometric volume is $2V_{1,1}$
\cite[(4.5)]{Mir08}, so $2^{-M}\prod\widetilde V_Z=\prod V_Z$ in our
orbifold convention \cite[\S1]{DN09}. The corresponding oriented,
ordered convention for closed surfaces is also recorded in \cite[Example~2.6]{AMI}.

For a nonseparating curve $N(\gamma)=2$, since the two cut boundaries can
be exchanged.  Counting both orientations therefore gives $\ell U_{g,n}$.
For a separating curve first mark an actual side $(i,A)$, $A\subseteq[n]$.
Forgetting the side identifies $(i,A)$ with $(g-i,A^c)$.
Usually these are two distinct records and $N(\gamma)=1$.  The sole
one-record case is $n=0$, $i=g/2$, when $N(\gamma)=2$.
Even if numerical types coincide when $n>0$, their disjoint labelled cusp
sets remain distinct, so no side-exchanging stabilizer fixes all labels.
Thus the unoriented density is one half of the ordered side sum in all
cases.  Summing subsets introduces exactly $\binom n{a_0}$; orientation
doubles the result to $\ell S_{g,n}$.  This proves
\eqref{orb:simpledensity} for continuous tests of compact support away
from zero.  Equality of the resulting locally finite measures extends it
to nonnegative Borel tests.

Tonelli and $2\sinh(\nu\ell/2)\ge\nu\ell$ give
\begin{equation}\label{orb:repetitionmajorant}
 \sum_{\nu\ge1}\int_0^{T/\nu}
 \frac{\ell^2|H(\nu\ell)|K_{g,n}(\ell)}{2\sinh(\nu\ell/2)}\,d\ell
 \le\frac{\zeta(3)}2T^2\|K_{g,n}\|_{C[0,T]}\|H\|_\infty<\infty.
\end{equation}
The kernel is a finite sum of volume polynomials at fixed $g,n$.
Signed Fubini now gives \eqref{orb:simpleexact}.  If $T=0$, both sides
vanish.  Thus repetitions cause no singularity at zero length.
\end{proof}

For a continuous function $k:[0,\infty)\to\mathbb R$, define the linear operation
\begin{equation}\label{orb:Loperation}
 \mathcal L_k(H)=\int_0^\infty\sum_{\nu\ge1}
 \frac{\ell^2k(\ell)H(\nu\ell)}{2\sinh(\nu\ell/2)}\,d\ell.
\end{equation}
If $|k(\ell)|\le A(1+\ell)^de^\ell$, the preceding argument gives
$|\mathcal L_k(H)|\le\zeta(3)AT^2(1+T)^de^T\|H\|_\infty/2$.
Set
\begin{equation}\label{orb:etadefinition}
 \eta_{j,n}=\mathcal L_{u_{j,n}+s_{j,n}},\qquad
 \Omega(H)=\mathcal L_{\phi^2}(H)
 =\int_0^\infty\sum_{\nu\ge1}
 \frac{2\sinh^2(\ell/2)}{\sinh(\nu\ell/2)}H(\nu\ell)\,d\ell.
\end{equation}
These functionals are defined on all compactly supported bounded even
Borel tests, independently of a support cutoff.  Applying
Lemma~\ref{orb:simplekernel} inside \eqref{orb:simpleexact} proves
\begin{equation}\label{orb:simpleexpansion}
 \begin{aligned}
 &\left|\E\mathcal O_H^{\rm simp}-\Omega(H)
 -\sum_{j=1}^Rg^{-j}\eta_{j,n}(H)\right|\\
 &\qquad\le g^{-R-1}[C(R+1)(n+1)]^{C(R+1)}
 (1+T)^{C(R+1)}e^T\|H\|_\infty.
 \end{aligned}
\end{equation}
This holds under $g\ge[C(R+1)(n+1)]^C$.  The same estimate without the genus factor,
with $R$ replaced by $j$, bounds $\eta_{j,n}$.

\subsection{Cusp-relative carriers and filling fibers}

The use of the subsurface filled by a geodesic follows the local-type
integration framework of Anantharaman and Monk
\cite[\S4.1, Theorems~5.7 and~6.1]{AMI} for closed ambient surfaces.
We use the same construction with inherited cusps, keeping their
labels and the boundary-slot symmetries explicit. A boundary slot is
a boundary component of the normalized carrier; two slots may map to
the same geodesic on the ambient surface.

The normalized carrier of a nonsimple geodesic is constructed from a
small regular neighborhood of its geodesic graph.  Fill disk regions,
complete peripheral punctured-disk regions to their inherited cusp ends,
and tighten the remaining boundary curves to positive geodesics.  If two
boundaries tighten to the same ambient cuff, retain two distinct boundary
slots in the normalization.  This is the construction preceding
\cite[Proposition 3.1]{HWX26}, with peripheral completion made explicit.
The carrier interior contains the geodesic, and its inclusion in the
ambient surface is an injective local isometry.  Its type $(h_0,a_0,k)$ records the genus, the number of inherited
cusps and the number of positive boundary slots, respectively.
Its complexity is $m=2h_0+a_0+k-2\ge1$.
An interior geodesic fills a carrier if the complement of a small regular
neighborhood of its graph consists of disks, cusp neighborhoods, and
annular collars of individual positive boundary components.

\begin{lemma}[Geometry and reconstruction]\label{orb:carriergeometry}
For an interior geodesic $\gamma$ filling a normalized carrier $Y$,
\begin{equation}\label{orb:carrierbounds}
 \Lambda:=\ell(\partial Y)\le2\ell_\gamma,
 \qquad 2\pi m=\operatorname{Area}(Y)\le4\ell_\gamma,
 \qquad \ell_\gamma\ge\pi m/2.
\end{equation}
Gluing a filled carrier to stable complementary components, and possibly
self-pairing some slots, reconstructs its canonical normalized carrier.
Intrinsic and ambient primitivity are equivalent.
\end{lemma}

\begin{proof}
Cut $Y$ along the graph of $\gamma$.  The resulting regions are disks,
peripheral punctured disks, and annuli adjacent to individual positive
boundary components.  Their graph-side lengths sum to $2\ell_\gamma$.
Each positive boundary geodesic minimizes its free homotopy class,
proving the perimeter estimate, with both slots counted in a self-pair.

The area bound also holds for carriers with cusps.  On a disk lifted to the
hyperbolic plane use a unit Busemann vector field of divergence one.
On a peripheral region's parabolic cyclic cover use $-y\partial_y$,
also of norm one and divergence one; truncate high in the cusp and let
the extra horocycle flux tend to zero.  On an annular region's hyperbolic
cyclic cover use $\tanh(r)\partial_r$, where $r$ is signed distance from
the core geodesic.  This field has divergence one, norm at most one,
and zero core flux.  The divergence theorem, with corners smoothed,
bounds each region's area by its graph-side length.  This gives the stronger estimate
$\operatorname{Area}(Y)\le2\ell_\gamma$; the weaker bound in
\eqref{orb:carrierbounds} suffices below.  Gauss--Bonnet gives
$2\pi m=\operatorname{Area}(Y)$, and the lower bound for $\ell_\gamma$ follows.

For reconstruction let $\iota:Y^\circ\to X$ be the injective local
isometry after gluing.  If a unit-speed closed geodesic $\alpha$ has an
ambient period $d$, then
$\iota(\alpha(t+d))=\iota(\alpha(t))$ for every $t$, so injectivity makes
$d$ an intrinsic period.  The converse is immediate.  This proves
primitivity in both directions and also preserves oriented geodesics
and their self-intersections.

Outside a small graph neighborhood, intrinsic disks and cusp regions
are absorbed by the carrier construction.  The two collars of a
self-pair form an annulus, whose two boundary germs tighten to the same
seam.  Each stable exterior component, with its attached collars,
retains negative Euler characteristic and cannot be absorbed.
Hence the construction recovers precisely the initial carrier and its
slot pairing.  Conversely, cutting an ambient canonical carrier gives
these records uniquely: equality of cuff images determines self-pairs,
and connected exterior regions determine complementary blocks.
A cuff has only two sides.  All remaining seams are essential and
nonparallel, as can also be seen by extending the cut components to
pants decompositions.  Finally a simple geodesic has an annular graph
neighborhood; adjoining only disks, cusp regions and boundary collars
cannot produce negative Euler characteristic.  Thus every filling
geodesic of a carrier with $m\ge1$ is nonsimple.
\end{proof}

Fix a carrier type $(h_0,a_0,k)$ with $k\ge1$, $0\le a_0\le n$
and $m=2h_0+a_0+k-2\ge1$. Write
$\mathcal P_{\rm or}^{\rm fill}(Y)$ for the oriented primitive interior
geodesics filling $Y$.  On each positive boundary fiber define pointwise
\begin{align}
 A_H(Y)&=\sum_{\alpha\in\mathcal P_{\rm or}^{\rm fill}(Y)}\sum_{\nu\ge1}
 \frac{\ell_\alpha H(\nu\ell_\alpha)}{2\sinh(\nu\ell_\alpha/2)},
 \notag\\
 I_H^{h_0,a_0,k}(\mathbf x)
 &=\int_{\mathcal M_{h_0,a_0+k}(\mathbf0_{a_0},\mathbf x)}
 A_H(Y)\,d{\rm WP}_{\rm orb}(Y),\qquad
 \Psi_H^{h_0,a_0,k}(\mathbf x)=
 \frac{I_H^{h_0,a_0,k}(\mathbf x)}
 {V_{h_0,a_0+k}(\mathbf0_{a_0},\mathbf x)}.
 \label{orb:fiberdefinition}
\end{align}
The numerator and denominator use the same orbifold measure, including
its factor one half on an intact one-holed torus
\cite[\S2.3.4 and Remark 2.7]{AMI} and \cite[\S1]{DN09}.

The oriented primitive count in \cite[Lemma 3.7 and its preceding
convention]{HWX26} applies to surfaces with cusps and geodesic boundaries:
\begin{equation}\label{orb:primitivecount}
 \#\{\alpha\in\mathcal P_{\rm or}(Y):\ell_\alpha\le T\}
 \le\left(h_0-1+\frac{a_0+k}{2}\right)e^{T+6}
       +3h_0-3+a_0+k\le Cme^T.
\end{equation}
By \eqref{orb:carrierbounds}, a contributing repetition has
$\nu\le2T/(\pi m)$.  Since $\ell/(2\sinh(\nu\ell/2))\le1/\nu$,
\begin{equation}\label{orb:fiberbound}
 |A_H|,\ |\Psi_H|\le Cm(1+T)e^T\|H\|_\infty,
 \qquad \Psi_H=0\quad\hbox{if }\Lambda>2T\hbox{ or }m>2T/\pi.
\end{equation}
The pointwise bound and finite WP volume imply that every fiber
integral is finite, even when the fiber moduli space is noncompact.
We also need measurability in the boundary lengths, including on
self-pairing diagonals. The maps
$\mathbf x\mapsto I_H^{h_0,a_0,k}(\mathbf x)$ and
$\mathbf x\mapsto\Psi_H^{h_0,a_0,k}(\mathbf x)$ are Borel on
$(0,\infty)^k$: the countable set of marked primitive filling
classes has continuous length functions in Fenchel--Nielsen coordinates;
the sum defining $A_H$ is therefore Borel.  A countable orbifold atlas
and measurable partition express its integral as parameter integrals
against the internal length and twist measures, with finite stabilizer
weights.  Apply parameter integration first to $A_{|H|}$ and then to the
positive and negative parts.  Division by the positive continuous
volume preserves measurability.  This defines $I_H$ and $\Psi_H$ at every positive boundary-length
vector, rather than only almost everywhere.  It therefore
justifies restriction to the equal-length diagonals of self-paired slots.

\subsection{The exact carrier formula}

Use the patterns of Theorem~\ref{vol:kernel}.  Recall that a pattern
$P=(M,\Pi)$ on $[k]$ has a partial matching $M$ of $2\sigma$ slots and
a partition of the remaining slots into nonempty blocks
$I_1,\ldots,I_b$, $b\ge1$, ordered by their least labels.
Each self-pair gives one gluing edge, and each remaining slot gives
one edge joining the carrier to an exterior component. Hence there
are $e_P=k-\sigma$ actual edges. Their lengths $\mathbf y$ induce
the slot lengths $\mathbf x_P$ and
\begin{equation}\label{orb:patternlengths}
 \Lambda_P=2\sum_{\text{self edges }e}y_e+
                  \sum_{\text{external edges }e}y_e,
 \qquad J_P=\prod_e y_e.
\end{equation}
Writing $k_i=|I_i|$, the complementary genera and cusp counts satisfy
\begin{equation}\label{orb:complementdata}
 \sum_i a_i=n-a_0,\quad
 \sum_i g_i=g-h_0-k+\sigma+b,\quad
 g_i,a_i\in\mathbb Z_{\ge0},\quad 2g_i+a_i+k_i-2\ge1.
\end{equation}
The nonnegative kernel supplied there is
\begin{equation}\label{orb:carrierkernel}
 F_{g,n}^{h_0,a_0,k;P}(\mathbf y)=J_P
 V_{h_0,a_0+k}(\mathbf0_{a_0},\mathbf x_P)
 \sum_{\eqref{orb:complementdata}}
 \frac{n!}{a_0!\prod_i a_i!}
 \frac{\prod_iV_{g_i,a_i+k_i}(\mathbf0_{a_i},\mathbf x_{I_i})}{V_{g,n}}.
\end{equation}

A \emph{rooted} carrier here means one whose boundary slots are
bijectively labelled by $[k]$. The factor $1/k!$ below forgets these
labels after integration. Patterns in this formula have at least
one stable exterior component. Carriers with no such component are
handled separately in the proof of Theorem~\ref{orb:expansion}.

\begin{proposition}[Rooted carrier integration]\label{orb:carrierintegration}
Let $\mathcal O_H^{h_0,a_0,k}$ restrict \eqref{orb:statistic} to primitive
nonsimple geodesics with this normalized carrier type and at least one
stable complementary component.  Then it is absolutely WP integrable
and
\begin{equation}\label{orb:carrierexact}
 \E_{g,n}\mathcal O_H^{h_0,a_0,k}
 =\frac1{k!}\sum_{P\in\mathscr P_k}
 \int_{(0,\infty)^{e_P}}
 \Psi_H^{h_0,a_0,k}(\mathbf x_P(\mathbf y))
 F_{g,n}^{h_0,a_0,k;P}(\mathbf y)\,d\mathbf y.
\end{equation}
\end{proposition}

\begin{proof}
First fix actual disjoint cusp-label sets $A_0,A_1,\ldots,A_b$ of sizes
$a_0,a_1,\ldots,a_b$.  Label all carrier boundary germs by $[k]$ and all
complementary boundaries by their incident slots.  Orient an external
edge by its carrier-side boundary orientation and a self-edge by its
smaller-labelled slot; order edges by their least carrier label.
Lemma~\ref{orb:carriergeometry} gives inverse cutting and gluing for
these rooted records.  The dual graph has $b+1$ vertices and
$k-\sigma$ edges, giving precisely the genus equation in
\eqref{orb:complementdata}.

The marked Fenchel--Nielsen product measure has one twist circle of
full marked period $y_e$ per actual edge.  After twist integration it is
\begin{equation}\label{orb:markedmeasure}
 \left(\prod_e y_e\,dy_e\right)d{\rm WP}_{\rm orb}(Y)
                         \prod_i d{\rm WP}_{\rm orb}(Z_i).
\end{equation}
This is the coordinate formula underlying
\cite[Theorem 2.5 and Remark 2.7]{AMI}; inherited cusp labels are kept as
in \cite[\S4]{Mir08}.  To justify its use for a nonconstant carrier
observable, take a local cut-moduli chart $[U/G]$.  The gluing chart is
the corresponding twist bundle
$[U\times\prod_e(\mathbb R/y_e\mathbb Z)/G]$.  Orbifold integration
divides the cover integral by $|G|$.  Although $G$ may rotate twists or
permute filling geodesics, the sum $A_H$ is invariant and independent of
the twists.  Integrating in the twists therefore multiplies $A_H$ by
$\prod_e y_e$, leaving exactly the cut-component orbifold weights.  Product charts and a
measurable partition give \eqref{orb:markedmeasure} globally.
For an intact torus, write $B$ for its raw geometric chart and
$\iota$ for its elliptic involution.  On the diagonal quotient
$[B\times(\mathbb R/y\mathbb Z)/\langle\iota\rangle]$, integration
gives $\frac12\int_B\int_0^y A_H\,dt
=y\int_B A_H\,d{\rm WP}_{\rm orb}$, with raw geometric measure
understood on the left.  Thus the full marked twist period and the
orbifold volume account for the half-twist symmetry exactly once.
A torus created only after regluing a cut pair of pants introduces no
additional factor of $1/2$.

Integrating a nonnegative $A_H$ first in the carrier, then in the other
components and twists, gives $I_HJ_P\prod_i V_{g_i,a_i+k_i}$.
The number of cusp-label assignments is exactly
$n!/(a_0!\prod_i a_i!)$.  Summation and division by $V_{g,n}$ produce
the integrand $\Psi_HF_{g,n}^P$ in \eqref{orb:carrierexact}.
Finally every normalized carrier has $k!$ bijective labellings of its
boundary germs, even when two germs have the same ambient image.
The invariant filling sum descends through this finite symmetric-group
quotient.  Orbifold integration divides by $k!$, including at metrics
with automorphisms.  Self-pair reversals, permutations of self-pairs,
and exchanges of exterior blocks are already included in this quotient;
there is no additional $2^{-\sigma}$, $\sigma!^{-1}$ or $b!^{-1}$.
Automorphisms of the marked configuration preserve the distinguished carrier.
For example, the three patterns with one self-pair on a pants carrier
with three boundary slots contribute $3/3!=1/2$ in total.

For nonnegative $H$, Tonelli's theorem proves the identity with
$+\infty$ allowed.  At fixed
$g,n$ the complementary-data sum is finite and its kernel is continuous
on the closed length orthant.  Equation~\eqref{orb:fiberbound} restricts
the integrals to the compact simplex $\Lambda_P\le2T$ and bounds
$\Psi_H$ there.  Applying the identity to $|H|$ proves absolute
integrability, then signed Fubini proves the stated formula.
\end{proof}

Let $\mathcal F_{j,n}^{h_0,a_0,k;P}$ be the canonical coefficients of
Theorem~\ref{vol:kernel}.  Define, without choosing $R$ or $T$,
\begin{equation}\label{orb:xidefinition}
 \xi_{j,n}^{h_0,a_0,k}(H)=\frac1{k!}\sum_{P\in\mathscr P_k}
 \int\Psi_H^{h_0,a_0,k}(\mathbf x_P)
                    \mathcal F_{j,n}^{h_0,a_0,k;P}(\mathbf y)\,d\mathbf y.
\end{equation}
For each test $H$, the fiber bound restricts these integrals to a
bounded length simplex, where the continuous kernels are bounded.
Thus the integrals converge absolutely and define real linear
functionals independently of any chosen support cutoff.  Since $\mathcal F_j=0$ for
$j<m+b-1$, and $b\ge1$,
\begin{equation}\label{orb:xivanishing}
 \xi_{j,n}^{h_0,a_0,k}=0\quad(j<m),\qquad
 \xi_{j,n}^{h_0,a_0,k}(H)=0\quad(m>2T/\pi).
\end{equation}
The second vanishing follows from the support condition in
\eqref{orb:fiberbound}.

To integrate the kernel remainder, note that
$|\mathscr P_k|\le(k+1)^{2k}$: encode each matched partner, or zero,
and then the least label in each remaining block.  Also $e_P\le k\le m+2$
and
\begin{equation}\label{orb:simplex}
 \operatorname{Leb}\{\Lambda_P\le2T\}
 \le(2T)^{e_P}/e_P!\le(1+2T)^k.
\end{equation}
The twist factors are already in $F_{g,n}$.  With $\tau=R+m+1$,
Theorem~\ref{vol:kernel}, \eqref{orb:fiberbound}, and
\eqref{orb:simplex} bound the integrated remainder by
\begin{equation}\label{orb:fixedcarrierintermediate}
 g^{-R-1}[C\tau(n+1)]^{C\tau}
 Cm(k+1)^{2k}(1+T)(1+2T)^{C\tau+k}e^{3T}\|H\|_\infty.
\end{equation}
As $k\le m+2\le3\tau/2$, its extra factors can be absorbed by enlarging
one absolute $C$.  For $g\ge[C\tau(n+1)]^C$ we obtain
\begin{equation}\label{orb:fixedcarrierexpansion}
 \left|\E\mathcal O_H^{h_0,a_0,k}
       -\sum_{j=1}^Rg^{-j}\xi_{j,n}^{h_0,a_0,k}(H)\right|
 \le g^{-R-1}[C\tau(n+1)]^{C\tau}
 (1+T)^{C\tau}e^{3T}\|H\|_\infty.
\end{equation}
The same calculation bounds $\xi_j$ without the genus factor, using
$\tau=j+m+1$.  This includes $R<m$, when every retained coefficient
vanishes, and $R=0$.

\subsection{Exhaustion and the full expansion}

The support of $H$ bounds the complexity of every contributing
carrier. This makes the nonsimple part a finite sum. We first exclude
carriers without a stable exterior component at the stated genus
threshold, then combine the remaining carrier expansions with the
simple-curve expansion.

\begin{theorem}[Uniform orbital expansion]\label{orb:expansion}
The statistic \eqref{orb:statistic} is absolutely WP integrable for every
fixed $g\ge3,n,T,H$.  There are fixed real linear functionals
$\kappa_{j,n}$, $j\ge1$, and an absolute $C\ge2$ such that, for every
integer $R\ge0$ and
$\rho=R+\lceil T\rceil+1$, whenever $g\ge[C\rho(n+1)]^C$,
\begin{equation}\label{orb:fullbound}
 \left|\E_{g,n}\mathcal O_H-\Omega(H)
       -\sum_{j=1}^Rg^{-j}\kappa_{j,n}(H)\right|
 \le g^{-R-1}[C\rho(n+1)]^{C\rho}
 (1+T)^{C\rho}e^{3T}\|H\|_\infty.
\end{equation}
The sequence is independent of genus, order and support cutoff, and
\begin{equation}\label{orb:fullcoefficientbound}
 |\kappa_{j,n}(H)|
 \le[C(j+1)(n+1)]^{C(j+1)}(1+T)^{C(j+1)}e^{3T}\|H\|_\infty.
\end{equation}
\end{theorem}

\begin{proof}
First consider a nonsimple carrier with no stable exterior component.
Every ambient cusp has been absorbed, and every positive boundary slot
is self-paired.  For some $\sigma\ge0$ its parameters satisfy
\begin{equation}\label{orb:wholecarrier}
 a_0=n,\quad k=2\sigma,\quad g=h_0+\sigma,\quad m=2g+n-2.
\end{equation}
This includes $\sigma=0$, namely the whole ambient surface.  By
\eqref{orb:carrierbounds},
$\ell_\gamma\ge\pi(g-1+n/2)$.  The threshold in
\eqref{orb:fullbound}, with $C\ge2$, implies $g\ge4\rho^2$, whereas
$T\le\rho-1$.  Hence these geodesics have length strictly greater than
$T$ and contribute nothing, even at the support endpoints.

Absolute integrability does not require the expansion threshold.
Every nonsimple primitive has $m\ge1$ and length at least $\pi/2$.
The ambient version of \eqref{orb:primitivecount} bounds the number of
primitives up to $T$ by $C(g+n)e^T$.  Therefore
\begin{equation}\label{orb:nonsimpleabsolute}
 |\mathcal O_H^{\rm nonsimp}(X)|
 \le C(g+n)(1+T)e^T\|H\|_\infty.
\end{equation}
Countable length-function sums are measurable.  This proves absolute
integrability of the nonsimple part for every fixed genus; the simple
part was handled in Lemma~\ref{orb:simpleintegration}.

Under the expansion threshold all contributing nonsimple carriers
have a stable complement, $k\ge1$, and $1\le m\le M:=\lfloor2T/\pi\rfloor$.
For $d\ge0$ let
\begin{equation}\label{orb:types}
 \mathscr C_d(n)=\{(h_0,a_0,k)\in\mathbb Z_{\ge0}^2\times\mathbb Z_{\ge1}:
                         \ a_0\le n,
                         \ 1\le2h_0+a_0+k-2\le d\}.
\end{equation}
At fixed $m$ there are at most $(m+2)^2$ types: choose $h_0$, then
$a_0\le m+1-2h_0$, and $k$ is forced.  Thus
$|\mathscr C_d(n)|\le(d+2)^3$.  Reconstruction gives the exact finite
identity
\begin{equation}\label{orb:exhaustion}
 \mathcal O_H^{\rm nonsimp}(X)=
 \sum_{(h_0,a_0,k)\in\mathscr C_M(n)}\mathcal O_H^{h_0,a_0,k}(X).
\end{equation}
This sum introduces no further choices of cusp or slot labels.
Define the global coefficients by the finite sums
\begin{equation}\label{orb:kappadefinition}
 \zeta_{j,n}(H)=\sum_{(h_0,a_0,k)\in\mathscr C_j(n)}
                     \xi_{j,n}^{h_0,a_0,k}(H),
 \qquad \kappa_{j,n}=\eta_{j,n}+\zeta_{j,n}.
\end{equation}
If $M>j$, the extra types in $\mathscr C_M(n)$ have zero coefficient
at order $j$; if $j>M$, the extra types in $\mathscr C_j(n)$ have
$\Psi_H=0$. These are the two vanishing properties in
\eqref{orb:xivanishing}. Hence the coefficient sum over
$\mathscr C_M(n)$ always equals $\zeta_{j,n}(H)$.  Thus the same
coefficients work for every support cutoff.

Apply \eqref{orb:fixedcarrierexpansion} at order $R$ to every type in
\eqref{orb:exhaustion}.  Its parameter satisfies
$R+m+1\le R+\lceil T\rceil+1=\rho$.  The common threshold therefore
applies, including types with $m>R$, which remain in the remainder.
The at most $(M+2)^3\le(2\rho)^3$ errors are absorbed in the constant
in \eqref{orb:fullbound}.  If $M=0$ the nonsimple statistic and all its
coefficients on $H$ vanish.  Adding \eqref{orb:simpleexpansion} proves
\eqref{orb:fullbound}.  For \eqref{orb:fullcoefficientbound}, sum the
coefficient estimate following \eqref{orb:fixedcarrierexpansion} over
$m\le j$, so that $j+m+1\le2(j+1)$, and absorb the cubic type count.
\end{proof}

\subsection{Geometric moments with growing order}
\label{sec:geometricmoments}

We next apply the orbital expansion to convolution powers of the
fixed test from Section~\ref{sec:quantitative}. For its kernel
$h=b_0*b_0$ and Fourier multiplier $f=\widehat h$, with
$\supp h\subset[-L,L]$, define for $p(x)=\sum_{t\ge1}a_tx^t$
\begin{equation}\label{orb:polynomialtest}
 \begin{aligned}
 H_p&=\check{(p\circ f)}=\sum_{t\ge1}a_th^{*t},
 &\mathcal O_X(p)&=\mathcal O_{H_p}(X),\\
 \omega_{1,n}(p)&=\Omega(H_p),
 &\omega_{j,n}(p)&=\kappa_{j-1,n}(H_p)\quad(j\ge2).
 \end{aligned}
\end{equation}
The coefficient functionals $\omega_{j,n}$ are fixed linear functionals
on $x\mathbb R[x]$. Since $\omega_{1,n}$ is independent of $n$,
write it as $\omega_1$.

\begin{proposition}[Uniform geometric moments]\label{orb:cgm}
There is one integer $c=c_f\ge2$, depending only on the fixed multiplier,
such that, simultaneously for every $n\ge0$, integers $q>t\ge1$, and
$g\ge[cq(n+1)]^c$,
\begin{equation}\label{orb:cgmbound}
 \left|g^{-1}\E_{g,n}\mathcal O_X(x^t)
       -\sum_{j=1}^{q-1}g^{-j}\omega_{j,n}(x^t)\right|
 \le g^{-q}[cq(n+1)]^{cq}.
\end{equation}
\end{proposition}

\begin{proof}
The convolution support rule and Young's inequality give
$\supp(h^{*t})\subset[-tL,tL]$ and
$\|h^{*t}\|_\infty\le M_f^t$, where
$M_f=\max\{1,\|h\|_1,\|h\|_\infty\}$ is fixed.
In Theorem~\ref{orb:expansion} take $R=q-2$ and $T=tL$.  Then
\begin{equation}\label{orb:momentparameters}
 \rho=q-1+\lceil tL\rceil\le(1+L)q,\qquad
 1+tL\le(1+L)q.
\end{equation}
After division by $g$, the remainder is bounded by
$g^{-q}[C\rho(n+1)]^{C\rho}(1+tL)^{C\rho}e^{3tL}M_f^t$.
The factor multiplying $g^{-q}$ has logarithm
$O_f(q\log(q(n+1)))$. Choosing the integer $c$ sufficiently large
therefore gives \eqref{orb:cgmbound} and the stated genus threshold.
In particular, the case $q=2,t=1$ uses $R=0$ and gives
a first-moment bound.
\end{proof}

\subsection{The exact first orbital correction}
\label{sec:firstcorrection}

The coefficient $\kappa_{1,n}$ of $g^{-1}$ in
Theorem~\ref{orb:expansion} is polynomial in $n$ of degree at most two.
We compute its nonseparating part, combine it with the separating
coefficient \eqref{orb:firstseparating}, and enumerate the five
contributing carrier configurations. The resulting formula separates
the constant, linear and quadratic dependence on $n$. We retain the definitions
of $\mathcal L_k$ and $I_H^{h_0,a_0,k}$ from
\eqref{orb:Loperation} and \eqref{orb:fiberdefinition}, including both
orientations, all repetitions and the orbifold volume convention.
Write
\begin{equation}\label{first:notation}
 \begin{aligned}
 \phi(x)=\frac{\sinh(x/2)}{x/2},\qquad
 &c(x)=\cosh(x/2),\\
 w(x)=x\phi(x)=2\sinh(x/2),\qquad
 &v(x)=\frac{x^2+4\pi^2}{48}.
 \end{aligned}
\end{equation}
Thus $v=V_{1,1}$ in our orbifold convention, and $\phi(0)=1$.

\begin{lemma}[First nonseparating coefficient]\label{first:nonseparating}
For every integer $n\ge0$ and $x\ge0$,
\begin{equation}\label{first:uexplicit}
 u_{1,n}(x)=\frac{(2-n)\phi(x)(c(x)+1)-c(x)^2}{\pi^2}.
\end{equation}
\end{lemma}

\begin{proof}
Mirzakhani--Zograf's first handle-ratio coefficient
\cite[Theorem~1.4, (1.4)]{MZ15} is
$p_n^1=(3-2n)/\pi^2$, including $n=0,1$.
For fixed $p\ge1$ and $\mathbf d\in\Z_{\ge0}^p$, put
$q=|\mathbf d|$, $r=\#\{i:d_i>0\}$ and $S_2=\sum_i d_i^2$.
The normalization in \cite[(B.3)]{HWX26} agrees with
\eqref{vol:volumepolynomial}. The pair and single sums in
\cite[Proposition~B.9]{HWX26} simplify to
\[
 A=2(q^2-S_2)+2(p-1)q-pr+\frac{r^2+r}{2},\qquad
 B=\frac32q-S_2-r.
\]
By \cite[Theorem~1.4]{MZ15}, both
$V_{G,p-1}/V_{G,p}$ and $V_{G-1,p+1}/V_{G,p}$ equal
$(8\pi^2G)^{-1}+O_p(G^{-2})$.
Applying that proposition with $p,\mathbf d$ fixed as $G\to\infty$
therefore gives $\pi^2b_{\mathbf d,p}^1=-A/2+B$, hence
\begin{equation}\label{first:intersection}
 b_{\mathbf d,p}^1
 =-\frac{q^2+(p-\tfrac52)q-r(2p-r-5)/4}{\pi^2}.
\end{equation}
The pair sum is empty when $p=1$.
For the all-zero multi-index this gives zero, as the normalized
intersection number is identically one.  These coefficients agree
with those of Proposition~\ref{vol:hma} by uniqueness of asymptotic
coefficients at each fixed marked-point count and multi-index.
All genera in this identification are large, so no terminal-volume
convention enters.

Set $z=x/2$, $F(z)=\sinh z/z$ and
$D=(z/2)\,d/dz$.  In the series for
$B_{1,n+2}^{(2)}(x,x)$, the weights are
$a_d(z)=z^{2d}/(2d+1)!$ and
$q=d_1+d_2$, $r=\mathbf1_{d_1>0}+\mathbf1_{d_2>0}$.
All series below and their differentiated series converge locally
uniformly.  Their sums satisfy
\begin{align}
 \sum_{d_1,d_2\ge0}q\,a_{d_1}a_{d_2}&=D(F^2),
 &\sum_{d_1,d_2\ge0}q^2a_{d_1}a_{d_2}&=D^2(F^2),\notag\\
 \sum_{d_1,d_2\ge0}r\,a_{d_1}a_{d_2}&=2F(F-1),
 &\sum_{d_1,d_2\ge0}r^2a_{d_1}a_{d_2}
 &=2F(F-1)+2(F-1)^2. \label{first:series}
\end{align}
Substituting $p=n+2$ into \eqref{first:intersection} gives
\begin{equation}\label{first:Bsum}
 B_{1,n+2}^{(2)}(x,x)
 =-\frac{D^2(F^2)+(n-\tfrac12)D(F^2)
              -(n-1)F(F-1)+\tfrac12(F-1)^2}{\pi^2}.
\end{equation}
By \eqref{orb:ucoefficient},
$u_{1,n}=p_n^1F^2+B_{1,n+2}^{(2)}$; the shift from $g-1$ to $g$
has no additional first-order term.  The identities
$D(F^2)=F\cosh z-F^2$ and
$D^2(F^2)=\tfrac12\cosh^2z-\tfrac32F\cosh z+
(1+\tfrac12z^2)F^2$, together with $z^2F^2=\cosh^2z-1$,
yield \eqref{first:uexplicit}.  Continuity gives the value
$u_{1,n}(0)=(3-2n)/\pi^2$.
\end{proof}

Define three kernels independent of $n$:
\begin{align}
 k_0(x)&=\frac{2\phi(x)(c(x)+1)-c(x)^2}{\pi^2}
                       +\frac{\phi(x)v(x)}{4\pi^2},\notag\\
 k_1(x)&=-\frac{\phi(x)(c(x)+1)}{\pi^2},&
 k_2(x)&=\frac{\phi(x)}{4\pi^2}. \label{first:kernels}
\end{align}
For the nonsimple contribution, write $I_{\rm t}(x)=I_H^{1,0,1}(x)$ and
$I_a(\mathbf x)=I_H^{0,a,3-a}(\mathbf x)$ for $a=0,1,2$.
Their cusp labels are intrinsic labels on the fixed carrier; ambient
cusp choices will appear only in the binomial factors below.
Set
\begin{align}
 \mathcal Z_0(H)=\frac1{8\pi^2}\bigg[
 &\int_0^\infty w(x)I_{\rm t}(x)\,dx\notag\\
 &+\frac16\int_{(0,\infty)^3}
       w(x_1)w(x_2)w(x_3)I_0(x_1,x_2,x_3)\,d\mathbf x\notag\\
 &+\frac12\int_{(0,\infty)^2}u\,w(y)I_0(u,u,y)\,du\,dy
 \bigg], \label{first:Zzero}\\
 \mathcal Z_1(H)&=\frac1{16\pi^2}\int_{(0,\infty)^2}
                 w(x_1)w(x_2)I_1(x_1,x_2)\,dx_1\,dx_2,
 \label{first:Zone}\\
 \mathcal Z_2(H)&=\frac1{8\pi^2}\int_0^\infty w(x)I_2(x)\,dx.
 \label{first:Ztwo}
\end{align}
\begin{proposition}[Exact quadratic cusp dependence]\label{first:quadratic}
For every compactly supported bounded even Borel $H$, the integrals
\eqref{first:Zzero}--\eqref{first:Ztwo} converge absolutely.  The fixed
linear functionals
\begin{equation}\label{first:Lambdas}
 \Lambda_i(H)=\mathcal L_{k_i}(H)+\mathcal Z_i(H),
 \qquad 0\le i\le2,
\end{equation}
satisfy, for every integer $n\ge0$,
\begin{equation}\label{first:quadraticidentity}
 \kappa_{1,n}(H)=\Lambda_0(H)+n\Lambda_1(H)+\binom n2\Lambda_2(H).
\end{equation}
There is an absolute $C$ such that, if $\supp H\subset[-T,T]$,
\begin{equation}\label{first:supportbound}
 \sum_{i=0}^2|\Lambda_i(H)|
 \le C(1+T)^4e^{2T}\|H\|_\infty,\qquad
 |\kappa_{1,n}(H)|
 \le C(n+1)^2(1+T)^4e^{2T}\|H\|_\infty.
\end{equation}
If $H$ is additionally continuous, nonnegative and not identically zero,
then
\begin{equation}\label{first:positivequadratic}
 \Lambda_2(H)\ge\frac1{4\pi^2}\int_0^\infty xH(x)\,dx>0.
\end{equation}
Thus for such a test the degree in $n$ is exactly two.
\end{proposition}

\begin{proof}
Lemma~\ref{first:nonseparating} and
\eqref{orb:firstseparating} give
$u_{1,n}+s_{1,n}=k_0+nk_1+\binom n2 k_2$.
Linearity of \eqref{orb:Loperation} identifies the simple contribution.

For the nonsimple contribution, \eqref{vol:recordcount} at $j=1$
forces $m=b=1$, with no small complementary components.  The possible
carrier types are exactly $(1,0,1)$ and $(0,a,3-a)$ for $a=0,1,2$.
For a surviving pattern, \eqref{vol:kernelleading} reduces to
\begin{equation}\label{first:carrierleading}
 \mathcal F_{1,n}^{h_0,a_0,k;P}
 =\frac{\binom n{a_0}}{8\pi^2}
       J_P V_{h_0,a_0+k}(\mathbf0_{a_0},\mathbf x_P)
       \prod_{\text{unpaired slots }j}\phi(x_j).
\end{equation}
Multiplication by $\Psi_H=I_H/V$ cancels only the displayed carrier
volume; it does not alter its intrinsic orbifold measure.

After extracting the common factor $(8\pi^2)^{-1}$, the cusp-label
and slot-symmetry weights are
\begin{center}
\small
\begin{tabular}{@{}lcc@{}}
\toprule
Carrier type $(h_0,a_0,k)$ & Self-pairs & Weight \\
\midrule
One-holed torus $(1,0,1)$ & $0$ & $1$ \\
Three-boundary pants $(0,0,3)$ & $0$ & $1/6$ \\
Three-boundary pants $(0,0,3)$ & $1$ & $1/2$ \\
One-cusp pants $(0,1,2)$ & $0$ & $n/2$ \\
Two-cusp pants $(0,2,1)$ & $0$ & $\binom n2$ \\
\bottomrule
\end{tabular}
\end{center}
The self-paired pants row combines the three choices of a pair,
giving $3/3!=1/2$. Its slot lengths are $(u,u,y)$ and its actual-edge
weight is $uy$. The self-edge contributes its twist length $u$,
while only the external edge contributes $\phi(y)$, giving the length
weight $u\,w(y)$. The torus convention is already accounted for in
\eqref{orb:carrierexact}. Pairing both slots of the one-cusp pants
would leave $b=0$ and is excluded.
These are all the patterns, so
$\zeta_{1,n}=\mathcal Z_0+n\mathcal Z_1+\binom n2\mathcal Z_2$.
This proves \eqref{first:quadraticidentity}, including $n=0,1$,
where unavailable cusp-labelled carriers have zero binomial weight.

It remains to justify absolute convergence and the support bound
\eqref{first:supportbound}.
For complexity one, \eqref{orb:fiberbound} gives
$|I_a|\le C(1+T)e^T\|H\|_\infty$ for pants and
$|I_{\rm t}(x)|\le C(1+x^2)(1+T)e^T\|H\|_\infty$ for the torus.
The integrands vanish unless the sum of carrier boundary lengths is
at most $2T$.  Thus the unpaired domains have
$\sum x_i\le2T$, and the self-paired domain has $2u+y\le2T$.
Since $w(x)\le e^{x/2}$, the product of external weights is at
most $e^T$.  The unpaired domains have dimension at most three;
the torus has one integration variable and a degree-two volume
factor; the self-paired domain has dimension two and the extra
factor $u\le T$.  In every case the bound is
$C(1+T)^4e^{2T}\|H\|_\infty$.
It also covers $T=0$, when every integral is zero.

The explicit kernels satisfy
$|k_i(x)|\le C(1+x)^2e^x$.
The bound following \eqref{orb:Loperation} therefore gives
$|\mathcal L_{k_i}(H)|\le C(1+T)^4e^T\|H\|_\infty$,
including the absolute sum of all repetitions.  This proves
\eqref{first:supportbound}.  Finally, for nonnegative $H$ both
$\mathcal Z_2(H)$ and every repetition in $\mathcal L_{k_2}(H)$
are nonnegative.  Its $\nu=1$ contribution is exactly
$(4\pi^2)^{-1}\int_0^\infty xH(x)\,dx$.
Since $H$ is continuous, even, and not identically zero, this integral
is positive, proving \eqref{first:positivequadratic}.
\end{proof}

\subsection{Cusp polynomiality at every fixed order}
\label{sec:cusppolynomiality}

The preceding quadratic formula is the first case of a degree bound
for every orbital coefficient. In a retained record of deficit $\delta$,
the cusp-label factor has degree at most $\delta+1$, while the shifted
coefficient at order $j-\delta$ has degree at most $j-\delta$.
Their product therefore has degree at most $j+1$. We establish these
bounds and justify integration of the coefficient polynomials below.
Finite interpolation then gives the binomial-basis formula for the
canonical orbital coefficients.
For a sequence $f(n)$ write
$\Delta f(n)=f(n+1)-f(n)$; $\Delta_n^r\kappa_{j,0}$ below denotes
the $r$th forward difference of $n\mapsto\kappa_{j,n}$ evaluated
at $n=0$. Put
$(z)_a=z(z-1)\cdots(z-a+1)$, with $(z)_0=1$.
Here polynomial extension always means extension of a coefficient,
not a volume with an inadmissible number of marked points.

\begin{lemma}[Polynomial boundary and shifted coefficients]
\label{poly:boundary}
Fix integers $s,\ell\ge0$ and put $p_0=\max\{1,\ell\}$.
The functions $B_{s,p}^{(\ell)}$, for integers $p\ge p_0$, are the
values of a polynomial in $p$ of degree at most $s$, with continuous
boundary-length functions as coefficients.
With $h_p^0=p_p^0=1$, the scalar coefficients $h_p^s,p_p^s$ are
polynomials in $p$ of degree at most $s$ on all integers $p\ge0$.
For fixed integers $H_0\ge0$, $\delta\ge1$ and $\ell\ge0$,
$Q_s^{n,H_0,\delta;\ell}$ likewise has a unique polynomial extension
in $n$ of degree at most $s$, agreeing with the canonical function
at every integer $n$ satisfying
\[
 n\ge0,\qquad p=n+2H_0-\delta\ge\max\{1,\ell\}.
\]
\end{lemma}

\begin{proof}
Lemma~4.7(i) of \cite{MZ15} fixes any
$\ell\ge1$ nonnegative descendant indices and pads with zeros.
It states that the coefficient of order $s$ is polynomial in the
total marked count $p\ge\ell$, of degree at most $s$.
By uniqueness of fixed-parameter asymptotic coefficients, these are
the $b_{(\mathbf0,\mathbf d),p}^s$ in
Proposition~\ref{vol:hma}; symmetry permits placing the padded zeros
before the displayed indices.

For $0\le i\le s$ define the Lagrange polynomials
\[
 L_i(p)=\prod_{\substack{0\le r\le s\\r\ne i}}
             \frac{p-(p_0+r)}{i-r}.
\]
For each fixed multi-index the source statement gives its exact
interpolation identity at $p_0,\ldots,p_0+s$.
Multiply by the weights in \eqref{vol:Bdefinition} and sum.
All $s+1$ node series converge locally uniformly and absolutely by
\eqref{vol:sourcecoefficients} and \eqref{vol:poissonsum}.
At every admissible integer $p$ the resulting series is the canonical
one.  Thus
\begin{equation}\label{poly:Binterpolation}
 B_{s,p}^{(\ell)}(\mathbf x)
 =\sum_{i=0}^s L_i(p)B_{s,p_0+i}^{(\ell)}(\mathbf x)
 \qquad(p\ge p_0).
\end{equation}
The right side supplies the continuous function-valued polynomial
extension. This interchange involves only the $s+1$ absolutely
convergent node series. If $\ell=0$, the normalized volume equals
one: $B_{0,p}^{(0)}=1$ and $B_{s,p}^{(0)}=0$ for $s\ge1$,
which proves that case directly.

For the scalar ratios use
\cite[Theorem~4.1(ii), (4.2)--(4.3), and Lemma~4.7(ii)]{MZ15}.
Their domain explicitly includes $p=0$.
Equation (4.2) is our $H_{G,p}$, whereas (4.3) is the reciprocal
of our $P_{G,p}$.  If
$a(t)=1+\sum_{r\ge1}a_r(p)t^r$ has
$\deg_p a_r\le r$, its formal inverse has coefficients
\[
 c_0=1,\qquad c_s=-\sum_{r=1}^s a_r c_{s-r},
\]
and hence $\deg_p c_s\le s$ by induction.
This establishes the asserted degree bound for both scalar sequences,
including $p=0,1$, after identifying coefficients by uniqueness.

The same degree rule is preserved by finite products, affine
translations of $p$, and fixed genus shifts.  In fact, the order-$q$
coefficient after shifting $g$ to $g-H_0$ is
\[
 \sum_{r=1}^q\binom{q-1}{r-1}H_0^{q-r}a_r(p),
\]
whose degree is at most $q$.
In \eqref{vol:telescoping}, the order-$q$ coefficient of each
geometric denominator factor is
$[-(n-\delta+j-2)/2]^q$, of degree at most $q$ in $n$.
The numbers of factors depend only on $H_0,\delta$.
Consequently $\deg_n v_q(n,H_0,\delta)\le q$.
The convolution defining $Q_s$ in the proof of
Lemma~\ref{vol:shiftedboundary}, with
$p=n+2H_0-\delta$, has degree at most $s$.
It agrees with the canonical coefficients at every admissible $n$.
There are infinitely many such integers, proving uniqueness.
\end{proof}

\begin{lemma}[Polynomial retained kernels]\label{poly:kernels}
For every $j\ge1$, $u_{j,n}$ and $s_{j,n}$ are polynomials in $n$
with continuous length-dependent coefficients, of degrees at most
$j$ and $j+1$, respectively.
For a fixed carrier type $(h_0,a_0,k)$ and pattern $P$,
$\mathcal F_{j,n}^{h_0,a_0,k;P}$ has a continuous function-valued polynomial
extension of degree at most $j+1$.
On nonnegative integers $n<a_0$ this extension is zero.
\end{lemma}

\begin{proof}
The formula \eqref{orb:ucoefficient} and
Lemma~\ref{poly:boundary} give $\deg_n\beta_{r,n}\le r$ and
$\deg_n p_n^r\le r$, so $\deg_n u_{j,n}\le j$.

For separating curves remove the condition $a_0\le n$ from
\eqref{orb:endpointset}, obtaining the finite set
\[
 \mathscr E_\delta^{\rm univ}
 =\{(h_0,a_0)\in\mathbb Z_{\ge0}^2:2h_0+a_0-1=\delta\}.
\]
In \eqref{orb:scoefficient}, replace each $Q_{j-\delta}$ by
its polynomial extension and $\binom n{a_0}$ by
$(n)_{a_0}/a_0!$, and sum over this universal set.
Each summand has degree at most
$a_0+j-\delta\le j+1$.
For an integer $n\ge a_0$, its other component has
$p=n-a_0+1\ge1$, so this is precisely the original term.
For $0\le n<a_0$, the polynomial falling factorial is zero.
Thus the extended sum equals $s_{j,n}$ at every nonnegative integer.

For a carrier pattern, remove $A_{\rm out}\le n$ from the records
in \eqref{vol:recordparameters}, retaining the distinguished block,
the stable small-component data and the exact deficit $\delta$.
Call the resulting set $\widehat{\mathscr E}_\delta$.
It is finite independently of $n$, since
$\sum_{i\ne *}c_i=\delta-m$ bounds every small genus and cusp count.
For a record $e$, put $A_e^{\rm lab}=a_0+\sum_{i\ne *}a_i$ and
write the prefactor independent of the ambient cusp count as
\begin{equation}\label{poly:recordprefactor}
 C_e(\mathbf y)=
 \frac{J_P(\mathbf y)V_{h_0,a_0+k}(\mathbf0_{a_0},\mathbf x)}
 {a_0!\prod_{i\ne *}a_i!}
 \prod_{i\ne *}V_{g_i,a_i+k_i}(\mathbf0_{a_i},\mathbf x_{I_i}).
\end{equation}
This continuous function is independent of $n$.
Denote polynomial extensions of shifted coefficients by
$\widehat Q$.  The polynomial kernel is
\begin{equation}\label{poly:universalF}
 \widehat{\mathcal F}_j(n;\mathbf y)=
 \sum_{\delta=1}^j\ \sum_{e\in\widehat{\mathscr E}_\delta}
 (8\pi^2)^{-\delta}(n)_{A_e^{\rm lab}}C_e(\mathbf y)
 \widehat Q_{j-\delta}^{\,n,H_e,\delta;k_{*,e}}
                    (\mathbf x_{I_{*,e}}).
\end{equation}
The combinatorial cusp budget
\eqref{vol:cuspbudget} holds for these records without any
ambient-cusp restriction, giving $A_e^{\rm lab}\le\delta+1$.
Every summand therefore has degree at most $j+1$.

For $n\ge A_e^{\rm lab}$, the identity
$2H_e+A_e^{\rm lab}-k_{*,e}=\delta$ gives
$p=n-A_e^{\rm lab}+k_{*,e}\ge k_{*,e}\ge1$.
Hence the term is exactly the canonical term in
\eqref{vol:canonicalcoefficients}.
For $0\le n<A_e^{\rm lab}$, its falling factorial vanishes.
The extension of $Q$ has already been defined, so this statement
never evaluates a volume with a negative number of cusps.
In particular all terms vanish if $n<a_0$.
This proves the claimed identification and continuity.
\end{proof}

\begin{theorem}[Fixed-order cusp polynomiality]\label{poly:orbital}
For every integer $j\ge1$ there are unique fixed real linear
functionals $\Lambda_{j,r}$, $0\le r\le j+1$, on compactly supported
bounded even Borel tests, such that for every integer $n\ge0$,
\begin{equation}\label{poly:orbitalidentity}
 \kappa_{j,n}(H)=\sum_{r=0}^{j+1}\binom nr\Lambda_{j,r}(H).
\end{equation}
They are independent of genus, cusp count, truncation and support
cutoff, and are explicitly specified by
\begin{equation}\label{poly:Newton}
 \Lambda_{j,r}(H)=\Delta_n^r\kappa_{j,0}(H)
 =\sum_{q=0}^r(-1)^{r-q}\binom rq\kappa_{j,q}(H).
\end{equation}
One absolute $C\ge2$ satisfies, simultaneously for all $j\ge1$,
$T\ge0$ and $\supp H\subset[-T,T]$,
\begin{align}
 \sum_{r=0}^{j+1}|\Lambda_{j,r}(H)|
 &\le[C(j+1)]^{C(j+1)}(1+T)^{C(j+1)}e^{3T}\|H\|_\infty,
 \label{poly:functionalbound}\\
 |\kappa_{j,n}(H)|
 &\le[C(j+1)]^{C(j+1)}(n+1)^{j+1}
       (1+T)^{C(j+1)}e^{3T}\|H\|_\infty.
 \label{poly:sharpcoefficientbound}
\end{align}
\end{theorem}

\begin{proof}
Fix $j$.  By Lemma~\ref{poly:kernels} and the absolute repetition
bound following \eqref{orb:Loperation}, the simple functional
$\eta_{j,n}$ is polynomial of degree at most $j+1$.
To justify passing polynomial coefficients through the integral,
express those coefficient functions
as finite linear combinations of the kernels at $n=0,\ldots,j+1$;
each such kernel satisfies \eqref{orb:simplecoefficientbound}.
Thus polynomial extraction commutes with the absolutely convergent
repetition integrals for every test in the statement.

For nonsimple carriers replace $\mathscr C_j(n)$ in
\eqref{orb:kappadefinition} by the finite set obtained by dropping
$a_0\le n$.  Its size remains at most $(j+2)^3$.
The fiber average $\Psi_H^{h_0,a_0,k}$ depends on its carrier
type, not on the ambient cusp count.
Use the polynomial kernel \eqref{poly:universalF} for each type.
On the fiber's support $\Lambda_P\le2T$, every coefficient of that kernel
is bounded, while \eqref{orb:fiberbound} bounds $\Psi_H$.
The length domain has finite measure by \eqref{orb:simplex}.
Each polynomial coefficient is therefore absolutely integrable.
Finite pattern and type sums preserve polynomiality and degree.
For every nonnegative integer $n$, unavailable types and records vanish
by Lemma~\ref{poly:kernels}, and all available terms agree with
\eqref{orb:xidefinition}--\eqref{orb:kappadefinition}.
Consequently $\kappa_{j,n}(H)$ is polynomial of degree at most
$d=j+1$ on all nonnegative integers.
The universal carrier-type set depends only on $j$; the support
condition \eqref{orb:fiberbound} removes the types with $m>2T/\pi$
for each test.

Newton interpolation now gives \eqref{poly:orbitalidentity} and
\eqref{poly:Newton}.  Uniqueness follows because the binomial
polynomials of orders $0,\ldots,d$ form a basis.
The finite formula also proves linearity and independence of
the auxiliary parameters.

To retain quantitative control in $j$, apply the estimate
\eqref{orb:fullcoefficientbound} only at $q=0,\ldots,d$.
Since $q+1\le j+2\le2(j+1)$, its maximum is at most
\[
 M_j(H):=[C(j+1)]^{C(j+1)}
           (1+T)^{C(j+1)}e^{3T}\|H\|_\infty
\]
after one absolute enlargement of $C$.
Equation~\eqref{poly:Newton} gives
$|\Lambda_{j,r}(H)|\le2^r M_j(H)$.
Summing over $r\le d$ contributes a factor at most $2^{d+1}$,
which is absorbed by enlarging $C$.  This proves
\eqref{poly:functionalbound}.
Finally $0\le\binom nr\le(n+1)^d$ for nonnegative integers $n$
and $0\le r\le d$, proving \eqref{poly:sharpcoefficientbound}.
\end{proof}

For $j=1$, uniqueness identifies these functionals with
$\Lambda_0,\Lambda_1,\Lambda_2$ of
Proposition~\ref{first:quadratic}.
For $j\ge2$, the theorem gives only an upper bound on the degree.
In the geometric-moment notation \eqref{orb:polynomialtest},
$\omega_{j,n}$ consequently has degree at most $j$ in $n$ for
$j\ge2$, while $\omega_1$ is independent of $n$.
Section~\ref{sec:remainders} controls the remainder when $n$ varies
with $g$; its final estimate is Theorem~\ref{addgeo:orbital}.

\section{Uniform remainder estimates}\label{sec:remainders}

This section sharpens the remainder in Theorem~\ref{orb:expansion}
while preserving its canonical coefficients.
Section~\ref{sec:allorderremainder} tracks the cusp powers in
the retained-record errors and complementary tails, giving a remainder
proportional to $(n+1)^{R+2}$ in Theorem~\ref{allrem:orbital}.
Section~\ref{sec:additivethreshold} then improves the tail estimate.
Theorem~\ref{addgeo:orbital} gives the same remainder under
an additive genus threshold: a quadratic cusp term with an absolute
coefficient, plus an order-dependent term linear in $n+1$.

Throughout, $g\ge3$ and $n\ge0$ are integers and $s=n+1$.
For integers $G\ge3$, $p\ge1$ and $0\le\ell\le p$, write
\[
 \mathcal V_{G,p}^{(\ell)}(\mathbf x)
 =\frac{V_{G,p}(\mathbf0_{p-\ell},\mathbf x)}{V_{G,p}}
\]
for $\mathbf x\in[0,\infty)^\ell$. All constants in this section
are absolute. As before, $C$ may increase between estimates; one
sufficiently large choice is independent of genus, cusp count, expansion
order, support size and carrier type.

\subsection{An all-order remainder with explicit cusp powers}
\label{sec:allorderremainder}

We first keep the power of the marked-point count explicit at every
order. Lemmas~\ref{allrem:boundary} and~\ref{allrem:shift} control
the retained records under genus thresholds linear in $s$.
Lemma~\ref{allrem:tail} controls the omitted records and is the source
of the stronger eighth-power cusp threshold in this subsection.

\begin{lemma}[Graded boundary errors]\label{allrem:boundary}
There is an absolute $C\ge2$ such that, for all integers $J\ge0$,
$p\ge1$, $0\le\ell\le p$, all $\mathbf x\in[0,\infty)^\ell$,
and every integer $G\ge[C(J+1)]^C(p+1)$, with $X=\sum_i x_i$,
\begin{equation}\label{allrem:boundaryerror}
 \left|\mathcal V_{G,p}^{(\ell)}(\mathbf x)
       -\sum_{j=0}^JG^{-j}B_{j,p}^{(\ell)}(\mathbf x)\right|
 \le G^{-J-1}(p+1)^{J+1}
 [C(J+1)]^{C(J+1)}(1+X)^{2J+2}e^{X/2}.
\end{equation}
For every $j\ge1$,
\begin{equation}\label{allrem:boundarycoeff}
 |B_{j,p}^{(\ell)}(\mathbf x)|
 \le(p+1)^j[C(j+1)]^{Cj}(1+X)^{2j}e^{X/2}.
\end{equation}
The scalar expansions of $H_{G,p}$ and $P_{G,p}$ have the same
order-$J$ error with the length factors removed, also for $p=0$.
Their order-$j$ coefficients satisfy the corresponding bound
$(p+1)^j[C(j+1)]^{Cj}$.
All genus thresholds here are independent of $X$.
\end{lemma}

\begin{proof}
Using $p+j+1\le(j+1)(p+1)$ and
$p+J+2\le(J+2)(p+1)$, the source quantities in
Proposition~\ref{vol:hma} satisfy
\[
 D_j(p)\le(p+1)^j[C(j+1)]^{Cj},
 \qquad G_J(p)\le[C(J+1)]^C(p+1).
\]
Enlarge $C$ so that the non-strict genus threshold in the lemma
implies $G>G_J(p)$.
The scalar assertions follow immediately, including $p=0$.
For boundary coefficients, apply \eqref{vol:poissonsum} to
\eqref{vol:Bdefinition}, absorbing $2j+1$ into the order factor.

Within the volume polynomial's dimension cutoff, the order-$J$
error follows by the same summation with
$D_{J+1}(p)E_{2J+2}(v)$, where $v=|\mathbf d|$.
The coefficient series also contains terms beyond this cutoff;
there $v\ge3G+p-2\ge G$.
For $0\le j\le J$, use \eqref{vol:degreetail} to bound
the omitted order-$j$ series by
\[
 G^{-J-1}D_j(p)(2J+2)!\,
       (2J+3)(1+X)^{2J+2}e^{X/2},
\]
where $D_0(p)=1$ for this estimate.
Since $(p+1)^j\le(p+1)^{J+1}$, summing these $J+1$ terms
costs only another factor $[C(J+1)]^{C(J+1)}$.
This includes the omitted constant series at $j=0$ and proves
\eqref{allrem:boundaryerror}.
For $\ell=0$ the exact ratio equals one and the positive-order
boundary coefficients vanish, so the assertion holds directly.
\end{proof}

\begin{lemma}[Graded shifted-volume error]\label{allrem:shift}
Let $J,H_0\ge0$, $\delta\ge1$ and $n\ge0$ be integers, and set
\[
 p=n+2H_0-\delta\ge1,\qquad
 L=J+H_0+\delta+1,\qquad s=n+1.
\]
For every integer $0\le\ell\le p$, all
$\mathbf x\in[0,\infty)^\ell$, and $g\ge[CL]^Cs$,
with $X=\sum_i x_i$, one has
\begin{align}
 &\left|\frac{V_{g-H_0,p}(\mathbf0_{p-\ell},\mathbf x)}{V_{g,n}}
 -(8\pi^2)^{-\delta}g^{-\delta}
       \sum_{j=0}^Jg^{-j}Q_j^{n,H_0,\delta;\ell}(\mathbf x)\right|
 \notag\\
 &\hspace{12mm}\le
 g^{-\delta-J-1}s^{J+1}[CL]^{C(J+1)}
                 (1+X)^{2J+2}e^{X/2}.
 \label{allrem:shifterror}
\end{align}
Here the $Q_j$ are exactly the canonical functions of
Lemma~\ref{vol:shiftedboundary}.
For $0\le j\le J$ they satisfy
\begin{equation}\label{allrem:Qcoeff}
 |Q_j^{n,H_0,\delta;\ell}(\mathbf x)|
 \le s^j[CL]^{Cj}(1+X)^{2j}e^{X/2}.
\end{equation}
\end{lemma}

\begin{proof}
We estimate the scalar shifted ratio first, then the normalized
boundary factor, and finally their product. All marked counts along
\eqref{vol:telescoping} are at most
$n+2H_0+\delta$, and all genus shifts are at most $H_0$.
For orders $1\le j\le J+1$, the exact source majorants are therefore
bounded by $A^j$ with
\[
 A=[C_1L]^{C_1}s,
\]
where $C_1$ is an auxiliary absolute constant.
The source thresholds satisfy a bound of the same form as $A$.
The genus-shift part of Lemma~\ref{vol:finiteproduct} replaces $A$
by at most $C(A+H_0+1)$, still of this form.
The geometric denominator factors have expansion coefficients
$[-(n-\delta+q-2)/2]^j$ and a geometric remainder, so they
obey the same majorant, including when the numerator in brackets
is negative.

There are $H_0+2\delta\le2L$ scalar factors.
The product part of Lemma~\ref{vol:finiteproduct}, after enlarging
the absolute exponent in $[CL]^C$, gives a normalized scalar ratio
$Z(g)$, defined by
$V_{g-H_0,p}/V_{g,n}=(8\pi^2)^{-\delta}g^{-\delta}Z(g)$.
Its canonical truncation
$Z^{[J]}(g)=\sum_{j=0}^Jv_jg^{-j}$ satisfies
\begin{equation}\label{allrem:scalarcontrol}
 \begin{aligned}
 &|v_j|\le A^j\quad(0\le j\le J),\qquad
 |Z-Z^{[J]}|\le(A/g)^{J+1},\\
 &|Z^{[J]}|\le2,\qquad A/g\le\tfrac14.
 \end{aligned}
\end{equation}
All subsequent enlargements of $C_1$ are absolute, so $A$ remains
linear in $s$. Impose $g\ge4A$ and choose $C_1$ large enough that this
also ensures $g-H_0\ge g/2$ and every strict source threshold.
At the end choose the statement's constant $C$ so that
$[CL]^Cs\ge4[C_1L]^{C_1}s$ for every $L\ge1$.

Let $D(g,\mathbf x)=\mathcal V_{g-H_0,p}^{(\ell)}(\mathbf x)$.
Shift the boundary expansion of Lemma~\ref{allrem:boundary}.
Its coefficients are exactly
\[
 \widetilde B_0=\Phi_\ell,\qquad
 \widetilde B_q=\sum_{j=1}^q
        \binom{q-1}{j-1}H_0^{q-j}B_{j,p}^{(\ell)}.
\]
Consequently
$|\widetilde B_q|\le A^q(1+X)^{2q}e^{X/2}$ for $0\le q\le J$.
Put
$D^{[J]}(g,\mathbf x)=\sum_{q=0}^Jg^{-q}\widetilde B_q(\mathbf x)$.
For the shifted remainder, the original boundary error costs at
most $(2A/g)^{J+1}(1+X)^{2J+2}e^{X/2}$.
In the truncated coefficient sum, use
\eqref{vol:shiftidentity} and
$\binom{j+v-1}{v}\le2^{j+v}$ for the discarded powers
$v\ge J+1-j$.
Since $H_0/g\le1/4$, their sum is bounded by
$[C(A+H_0+1)/g]^{J+1}(1+X)^{2J}e^{X/2}$.
After another fixed enlargement of $A$, this proves
\begin{equation}\label{allrem:shiftboundary}
 |D-D^{[J]}|
 \le(A/g)^{J+1}(1+X)^{2J+2}e^{X/2}.
\end{equation}

Use the exact decomposition
\[
 ZD-Z^{[J]}D^{[J]}
 =(Z-Z^{[J]})D+Z^{[J]}(D-D^{[J]}).
\]
The exact bound $|D|\le e^{X/2}$ and
\eqref{allrem:scalarcontrol}--\eqref{allrem:shiftboundary}
control these two terms with one exponential length factor.
For the terms discarded from the product of the truncations, first
take out $(1+X)^{2J}e^{X/2}$.
The remaining sum is at most
\[
 \sum_{\substack{0\le u,v\le J\\u+v\ge J+1}}
       (A/g)^{u+v}
 \le (J+1)^2(A/g)^{J+1}.
\]
The polynomial factor in $J$ is absorbed into $[CL]^{C(J+1)}$.
Only the boundary factor carries length growth, so the scalar
smallness condition $A/g\le1/4$ remains independent of $X$.
The finite convolution $\sum_{u+v=j}v_u\widetilde B_v$ equals
$Q_j$, proving both assertions.
When $J=0$ the discarded sum is empty and the same argument
applies.
\end{proof}

\begin{lemma}[Positive tails with explicit cusp power]
\label{allrem:tail}
For every integer $R\ge0$, every admissible carrier type of complexity $m\ge1$, and
$t=R+m+1$,
\begin{equation}\label{allrem:zerotail}
 g\ge[Kt(n+1)]^8
 \quad\Longrightarrow\quad
 \mathcal T_{g,n}^{h_0,a_0,k}(R)
 \le g^{-R-1}(n+1)^{R+2}[Ct]^{Ct}.
\end{equation}
On a fixed pattern, its boundary tail has the additional factor
$(1+\Lambda)^k e^\Lambda$.
The separating kernel with deficit at least $R+1$ satisfies
\begin{equation}\label{allrem:septail}
 S_{g,n}^{\ge R+1}(x)
 \le g^{-R-1}(n+1)^{R+2}
       [C(R+1)]^{C(R+1)}e^x
 \qquad\bigl(g\ge[K(R+1)(n+1)]^8\bigr).
\end{equation}
\end{lemma}

\begin{proof}
Return to the near- and deep-deficit estimates in the proof of
Theorem~\ref{vol:tail}. Set $D=\max\{m,R+1\}$ and retain the
bound $(N)_D\ge g^D$ before weakening $g^{-D}$ to $g^{-R-1}$.
The near sector, including its endpoint, is
bounded by
$2P_bb\,s^{D+1}(D+3)^{2b}D!g^{-D}$.
Under the same eighth-power threshold $g\ge s$, and
\begin{equation}\label{allrem:tailidentity}
 g^{-D}s^{D+1}
 =g^{-R-1}s^{R+2}(s/g)^{D-R-1}
 \le g^{-R-1}s^{R+2}.
\end{equation}
All other factors, including the pattern count, are bounded by
$[Ct]^{Ct}$ because $D\le t$ and $b\le t+1$.
The deep sector in \eqref{vol:deeptail} is already
$P_bg^{-R-1}$ under the same threshold.
For $b=1$ it is empty; the endpoint of the near sector remains
included.  This proves \eqref{allrem:zerotail}, also for $R=0$.
The numerator-length calculation following
\eqref{vol:deeptail} supplies the stated boundary factor.

For the separating kernel put $r=R+1$ in
\eqref{orb:sepratio}.  The uncompressed near bound is
$12Bs^{R+2}(R+1)!g^{-R-1}$, and
\eqref{orb:sepdeep} bounds the deep part by $Bg^{-R-1}$.
Boundary growth adds at most $e^x$.
Absorb the factorial into the order factor to obtain
\eqref{allrem:septail}, including the central tie and $n=0$.
\end{proof}

\begin{proposition}[Kernel errors with explicit cusp powers]
\label{allrem:kernels}
Let $R,n\ge0$ be integers and put $r=R+1$, $s=n+1$. Under
$g\ge[Cr]^Cs^8$, the canonical simple expansions obey, for $x\ge0$,
\begin{align}
 \left|U_{g,n}(x)-\phi(x)^2-\sum_{j=1}^Rg^{-j}u_{j,n}(x)\right|
 &\le g^{-R-1}s^{R+1}[Cr]^{Cr}(1+x)^{2R+2}e^x,
 \label{allrem:nonsep}\\
 \left|S_{g,n}(x)-\sum_{j=1}^Rg^{-j}s_{j,n}(x)\right|
 &\le g^{-R-1}s^{R+2}[Cr]^{Cr}(1+x)^{2R+2}e^x.
 \label{allrem:sep}
\end{align}
For a fixed carrier and pattern, let $t=R+m+1$ and let
$\mathbf y\in[0,\infty)^{e_P}$ be its edge lengths.
Under $g\ge[Ct]^Cs^8$,
\begin{equation}\label{allrem:kernel}
 \left|F_{g,n}^{h_0,a_0,k;P}(\mathbf y)
       -\sum_{j=1}^Rg^{-j}\mathcal F_{j,n}^{h_0,a_0,k;P}(\mathbf y)\right|
 \le g^{-R-1}s^{R+2}[Ct]^{Ct}
               (1+\Lambda)^{2t}e^\Lambda.
\end{equation}
The coefficients are given by \eqref{orb:ucoefficient},
\eqref{orb:scoefficient}, and \eqref{vol:canonicalcoefficients},
respectively.
\end{proposition}

\begin{proof}
For the nonseparating ratio use
$U_{g,n}=P_{g,n}\mathcal V_{g-1,n+2}^{(2)}(x,x)$.
Lemma~\ref{allrem:boundary}, its genus shift by one, and the
scalar--boundary product argument just proved give
\eqref{allrem:nonsep}.  Replacing $n+3$ by $3s$ and $1+2x$
by $2(1+x)$ costs only the displayed order factor.
The formal convolution is precisely \eqref{orb:ucoefficient}.

For a retained separating endpoint of deficit $1\le\delta\le R$,
take $J=R-\delta$ and $H_0=h_0$ in
Lemma~\ref{allrem:shift}.
Its label factor is at most $s^{a_0}$, with $a_0\le\delta+1$.
The small volume is bounded by
$[C(\delta+1)]^{C(\delta+1)}e^{x/2}$.
Since $h_0\le(\delta+1)/2$, the shifted lemma has $L\le C(R+1)$.
The cusp factor in the resulting error is therefore bounded by
\[
 s^{a_0+R-\delta+1}\le s^{R+2}.
\]
The two boundary exponentials give $e^x$, and
$2(R-\delta)+2\le2R+2$.
Sum the at most $O((R+1)^2)$ endpoints, with their two ordered
placements, and add \eqref{allrem:septail}.
This proves \eqref{allrem:sep} with the canonical
coefficients from \eqref{orb:scoefficient}.

For a retained carrier record $e$ of deficit $\delta\le R$,
keep the label factor separate in
\eqref{vol:recordprefactor}.
The volume-envelope calculation there gives the sharper form
\begin{equation}\label{allrem:recordprefactor}
 |A_e(\mathbf y)|
 \le s^{A_{\rm out}}[C(\delta+1)]^{C(\delta+1)}
        (1+\Lambda)^k e^{\Lambda-X_*/2},
 \qquad A_{\rm out}\le\delta+1.
\end{equation}
Indeed, the product of the zero-length volumes of the carrier and
small components is at most
$U^bA_0^{\delta-b}(\delta-b)!$, and the remaining multinomial
is at most $s^{A_{\rm out}}$.
Take $J=R-\delta$ in Lemma~\ref{allrem:shift}.
The record bounds give $H_0\le(R+m+2)/2$ and $L\le Ct$.
Multiplication by \eqref{allrem:recordprefactor} gives cusp power
\begin{equation}\label{allrem:recordcusp}
 s^{A_{\rm out}+J+1}\le s^{R+2},
\end{equation}
exponential factor $e^\Lambda$, and length power at most
$k+2(R-\delta)+2\le2t$, since $k\le m+2$ and $\delta\ge m$.
The record count in \eqref{vol:recordcount} is at most $[Ct]^{Ct}$.
Summing and adding Lemma~\ref{allrem:tail} proves
\eqref{allrem:kernel}.
Choose $C\ge8$ large enough that $[Ct]^Cs^8$ dominates
both the tail threshold and every shifted-ratio threshold.

The record reconstruction is valid because this threshold implies
$N>2R$ and $g>H_0$.
Available records satisfy $p=n-A_{\rm out}+k_*\ge k_*\ge1$;
unavailable cusp allocations are omitted before any volume is
evaluated.
The retained sums are exactly the finite definitions
\eqref{vol:canonicalcoefficients}.
When there are no retained records, in particular if $R=0$ or
$m>R$, the positive tail proves the assertion directly.
\end{proof}

\begin{theorem}[All-order orbital remainder with explicit cusp power]
\label{allrem:orbital}
There is an absolute $C\ge8$ such that, for all integers $R,n\ge0$,
$T\ge0$, and bounded real even Borel $H$ supported in $[-T,T]$,
putting $\rho=R+\lceil T\rceil+1$, the condition
\begin{equation}\label{allrem:threshold}
 g\ge[C\rho]^C(n+1)^8
\end{equation}
implies
\begin{equation}\label{allrem:orbitalbound}
 \left|\E_{g,n}\mathcal O_H-\Omega(H)
          -\sum_{j=1}^Rg^{-j}\kappa_{j,n}(H)\right|
 \le g^{-R-1}(n+1)^{R+2}
       [C\rho]^{C\rho}e^{3T}\|H\|_\infty.
\end{equation}
The empty sum at $R=0$ is zero.
The coefficients are those of Theorem~\ref{orb:expansion},
independent of genus, truncation order and support cutoff.
\end{theorem}

\begin{proof}
Apply the absolute repetition estimate following
\eqref{orb:Loperation} to
\eqref{allrem:nonsep}--\eqref{allrem:sep}.
It bounds the simple error by
\[
 g^{-R-1}s^{R+2}[C(R+1)]^{C(R+1)}
       (1+T)^{2R+4}e^T\|H\|_\infty.
\]
Since $R+1\le\rho$ and $1+T\le\rho$, this has the form
on the right of \eqref{allrem:orbitalbound}.

Increasing $C$ in \eqref{allrem:threshold} ensures $g\ge4\rho^2$.
The whole-carrier lower bound
$\pi(g-1+n/2)>T$ excludes all carriers without a stable exterior,
including entirely self-paired carriers.
Thus the exact exhaustion \eqref{orb:exhaustion} sums over
$1\le m\le M=\lfloor2T/\pi\rfloor$, with $k\le m+2$.
For such a type, $t=R+m+1\le\rho$, so
Proposition~\ref{allrem:kernels} applies under the common threshold.
For one pattern, combine \eqref{allrem:kernel},
\eqref{orb:fiberbound}, and \eqref{orb:simplex}.
Its integrated error is at most
\[
 g^{-R-1}s^{R+2}[Ct]^{Ct}
       Cm(1+T)(1+2T)^{2t+k}e^{3T}\|H\|_\infty.
\]
There are at most $(k+1)^{2k}$ patterns and $(m+2)^2$ types
at each complexity.
The factor $1/k!$ can be discarded in this upper bound.
Since $k\le m+2\le2\rho$ and $M\le T$, all length powers and
type/pattern counts are absorbed into $[C\rho]^{C\rho}$.
No ambient cusp factor enters this last summation.

The uniform fiber bound \eqref{orb:fiberbound} and the bounded
length simplex give absolute convergence of all error integrals.
The two vanishing properties in \eqref{orb:xivanishing}
identify the subtracted coefficient sum exactly with
$\zeta_{j,n}$, even when $j>M$ or $M>R$.
Adding the simple part gives
$\kappa_{j,n}=\eta_{j,n}+\zeta_{j,n}$ and proves the result.
For $M=0$ the nonsimple part and its coefficients vanish on $H$;
for $T=0$ all terms vanish.
This also proves the result for the stated bounded Borel tests.
\end{proof}

The power $R+2$ agrees with the degree bound for the first omitted
coefficient. This does not assert optimality of the remainder.
The retained-record estimates already require genus only linear in
$n+1$; the next subsection lowers the overall threshold by sharpening
the positive tails.

\subsection{An additive cusp and order threshold}
\label{sec:additivethreshold}

The improvement comes from retaining the full binomial lower bound
in the deep tail. Lemma~\ref{addgeo:tails} gives an additive threshold;
Proposition~\ref{addgeo:kernels} combines it with the retained-record
estimates already proved. Theorem~\ref{addgeo:orbital} then integrates
these kernel bounds. We continue to write $s=n+1$.

\begin{lemma}[Positive tails at an additive threshold]
\label{addgeo:tails}
There are absolute $K,C\ge2$ such that, for every carrier type of
complexity $m\ge1$, integer $R\ge0$, and $t=R+m+1$,
\begin{equation}\label{addgeo:complementtail}
 g\ge K(s^2+t^2s)
 \quad\Longrightarrow\quad
 \mathcal T_{g,n}^{h_0,a_0,k}(R)
 \le g^{-R-1}s^{R+2}[Ct]^{Ct}.
\end{equation}
For a fixed pattern, replacing the numerator volumes by their boundary
volumes and inserting the actual-edge Jacobian gives the additional
factor $(1+\Lambda)^ke^\Lambda$.
Writing $r=R+1$, the ordered separating tail satisfies
\begin{equation}\label{addgeo:separatingtail}
 S_{g,n}^{\ge r}(x)
 \le g^{-r}s^{r+1}[Cr]^{Cr}e^x,
 \qquad g\ge K(s^2+r^2s),\quad x\ge0.
\end{equation}
The assertions include $R=0$, $n=0$, and empty admissible sums.
\end{lemma}

At $x=0$ and for fixed $r$, the powers of $g$ and $n$ in
\eqref{addgeo:separatingtail} agree with
\cite[Lemma~B.14]{HWX26} in the regime $n=o(\sqrt g)$, with
$k=r+1$ there. Here we control the dependence on $r$ and obtain
the explicit threshold $g\ge K(s^2+r^2s)$, together with the
corresponding complementary-carrier estimate.

\begin{proof}
We split the complementary tail into near and deep deficits. The
near sum is geometric; in the deep sum the full binomial denominator
must absorb both the cusp-label count and the number of allocations.
Choose an integer $K\ge\max\{4096,\kappa_{\rm v}\}$. Put
\begin{equation}\label{addgeo:tailparameters}
 \begin{aligned}
 &r=R+1,\quad N=2g+n-3,\quad D_*=\max\{m,r\},\\
 &J=\lfloor N/(64s)\rfloor,\quad L_s=\log(64s).
 \end{aligned}
\end{equation}
Here $D_*$ is a lower bound for every deficit in the tail under consideration and
$J$ is the splitting point. We have $t\ge2$, $D_*\le t-1$ and
$b\le k\le m+2\le t+1$.
Since $s+t^2\ge2t$, the new hypothesis gives
$g\ge128ts$ and $g\ge\kappa_{\rm v}s^2$. Consequently
\begin{equation}\label{addgeo:tailprerequisites}
 \begin{gathered}
 g\le N\le3g,\quad N+1\le4g,\quad N-D_*+1\ge g,\\
 N\ge128s,\quad J\ge t,\quad
 1\le J\le N/2,\quad J\ge N/(128s),\quad s/g\le1.
 \end{gathered}
\end{equation}
Indeed $g\ge128ts$ implies $g\ge m+2$, $g\ge D_*+2$,
$g\ge s$ and $N/(64s)\ge2t$; the floor is at least half
$N/(64s)$ because the latter is at least two.

Fix a pattern with $b\ge1$ exterior blocks.
Its indices $R_i=c_i-1$ satisfy
$\sum_iR_i=N-d$, $d=m+b-1$, and
$\delta=N-\max_iR_i\ge d$.
The volume envelope applies directly under
\eqref{addgeo:tailprerequisites}: the exponential factor arising from the
denominator estimate is at most $e^{1/2}$, and the largest component
complexity is at least $g/b$. Thus \eqref{vol:factorialenvelope} gives
\[
 W\le P_b\frac{\prod_iR_i!}{N!}
            \frac{n!}{a_0!\prod_i a_i!},
 \qquad P_b=B^{b+2}\sqrt b\,(m-1)!,
\]
with one absolute $B\ge2$.
The outside-cusp budget \eqref{vol:cuspbudget} is
$A_{\rm out}\le\delta+1$, so the single multinomial is bounded
both by $s^{\delta+1}$ and by $(b+1)^n$.
The fixed-deficit count is at most
$b(\delta+3)^{2(b-1)}$: choose the maximal block, the remaining
indices of total $\delta-d$, and their cusp counts; their genera
are then fixed. Ignoring the parity and availability restrictions only
enlarges these positive sums.

It follows that the contribution at deficit $\delta$ is at most
$P_bb x_\delta$, where
\[
 x_\delta=\frac{s^{\delta+1}(\delta+3)^{2b}}{\binom N\delta}.
\]
For $D_*\le\delta<J$, use $b\le\delta+2$ to obtain
\[
 \frac{x_{\delta+1}}{x_\delta}
 =\frac{s(\delta+1)}{N-\delta}
       \left(1+\frac1{\delta+3}\right)^{2b}
 \le e^2/63<1/2.
\]
The ratio includes the step from $J-1$ to $J$.
Keeping $(N)_{D_*}\ge g^{D_*}$, the inclusive near sum is at most
\begin{equation}\label{addgeo:neartail}
 2P_bb(D_*+3)^{2b}D_*!\,g^{-D_*}s^{D_*+1}
 \le2P_bb(D_*+3)^{2b}D_*!\,g^{-r}s^{r+1}.
\end{equation}
The last inequality uses
$g^{-D_*}s^{D_*+1}=g^{-r}s^{r+1}(s/g)^{D_*-r}$.
When $b=1$, the only possible deficit is $m\le J$:
the tail is empty if $m<r$, and otherwise is already
included in \eqref{addgeo:neartail}.

For $b\ge2$ and $\delta>J$, every $R_i<N-J$.
The capped-factorial inequality \eqref{vol:cappedfactorials}
therefore applies, since $\sum_iR_i=N-d\le N$.
The entire deep sector has at most $(N+1)^{2b}$ allocations:
each cusp count is at most $n\le N$, and each genus is at most
$g\le N$, because $b\le k-2\sigma$ in the genus equation.
Its contribution is bounded by
\begin{equation}\label{addgeo:deepraw}
 P_b(b+1)^n(N+1)^{2b}\binom NJ^{-1}.
\end{equation}
The product formula for the binomial coefficient now gives
\begin{equation}\label{addgeo:entropy}
 \binom NJ
 =\prod_{j=0}^{J-1}\frac{N-j}{J-j}
 \ge(N/J)^J\ge(64s)^J
 \ge\exp\!\left(\frac{NL_s}{128s}\right).
\end{equation}
The factor $\log(64s)$ in the exponent absorbs the cusp-label count
at the additive threshold.

To bound the deep sector by $P_bg^{-r}$, it remains to absorb the
logarithm of $g^r(b+1)^n(N+1)^{2b}$. The inequalities
$r+2b\le4t$, $b+1\le2t$ and $2b\le3t$ give
\begin{equation}\label{addgeo:deeplogarithm}
 r\log g+n\log(b+1)+2b\log(N+1)
 \le F_t(g):=4t\log g+s\log(2t)+3t\log4.
\end{equation}
We first check this at the threshold value
$g_0=Ks(s+t^2)$, then use monotonicity in $g$. The calculation holds
for all real $t,s\ge1$, including the case $t=1$ needed for the
separating tail. Since $s+t^2\le2st^2$, we have
\begin{align}
 \frac{sF_t(g_0)}{g_0L_s}
 &\le
 \frac{4t\log(2K)+8t\log s+8t\log t
                    +s\log(2t)+3t\log4}
      {K(s+t^2)L_s}\notag\\
 &\le\frac1K\left\{
 \frac{4\log(2K)+3\log4}{\log64}
       +8+\frac8{e\log64}+1\right\}
 <\frac1{128}. \label{addgeo:deepthreshold}
\end{align}
For the second inequality use $t/(s+t^2)\le1/t\le1$,
$L_s\ge\log64$, $\log s/L_s\le1$, and
$(\log t)/t\le1/e$.
For the remaining label term, if $t\le s$ then
$\log(2t)\le L_s$; if $t>s$, use
$s/(s+t^2)\le1/t$ and $\log(2t)\le t$.
In both cases $s\log(2t)/((s+t^2)L_s)\le1$.

The final strict bound holds for every $K\ge4096$.
At $K=4096=2^{12}$ the first fraction in braces is $29/3$,
and $8/(e\log64)<4/3$, so the braces are less than $20$ and
$20/4096<1/128$.
The function $\log(2K)/K$ and the other constant multiples of
$1/K$ decrease on this range.
For fixed $t,s$, also
\[
 \left(\frac{F_t(g)}g\right)'
 =\frac{4t-4t\log g-s\log(2t)-3t\log4}{g^2}<0
 \qquad(g\ge e).
\]
Since $g_0>e$, \eqref{addgeo:deepthreshold} holds for every
$g\ge g_0$. It gives
$F_t(g)<gL_s/(128s)\le NL_s/(128s)$.
Equations \eqref{addgeo:deepraw}--\eqref{addgeo:deeplogarithm}
therefore bound the whole deep sector by $P_bg^{-r}$.

The remaining factors are independent of the ambient cusp count.
For $t\ge2$, one has
\[
 P_b\le(Ct)^{3t},\qquad
 2b(D_*+3)^{2b}D_*!\le(2t)^{5t},\qquad
 |\mathscr P_k|\le(4t)^{4t}.
\]
The last count encodes each slot's matched partner or zero and
the least label of its exterior block; there is no extra $b!$.
Summing the near and deep bounds proves
\eqref{addgeo:complementtail} after enlarging $C$.
For a fixed pattern, the numerator volumes have total boundary
length $2\Lambda-2\sum_{\rm self}y_e$.
The boundary estimate costs at most $e^\Lambda$, and the
actual-edge Jacobian is at most $(1+\Lambda)^k$.
This proves the boundary assertion on the fixed pattern's coordinate space,
including zero lengths by continuity.

For the ordered separating tail use the same $N,J,L_s$, with
$g\ge K(s^2+r^2s)$ and $r\ge1$.
The analogues of \eqref{addgeo:tailprerequisites} give
$J\ge r$ and $N-r+1\ge g$.
The volume envelope \eqref{orb:septerm} applies under this
threshold and bounds each zero-length summand by
$B\binom na/(\delta\binom N\delta)$.
A smaller side has at most $\delta+1$ cusp labels, and the two
ordered placements give the fixed-deficit bound
$6Bs^{\delta+1}/\binom N\delta$.
Its successive ratio is at most $1/63$ on $r\le\delta<J$.
Thus the inclusive near sum is bounded by
$12Br!g^{-r}s^{r+1}$.

For $\delta>J$, binomial unimodality and symmetry give
$\binom N\delta\ge\binom NJ$, including the equal central
binomial coefficients when $N$ is odd.
The original ordered cusp-label sum is at most $(g+1)2^n$, so the deep
part is at most
\[
 B(g+1)2^n\exp\!\left(-\frac{NL_s}{128s}\right).
\]
There is no further factor two in this ordered sum.
The logarithm of $g^r(g+1)2^n$ is at most
\[
 (r+1)\log g+s\log2
 \le2r\log g+s\log2\le F_r(g).
\]
Apply \eqref{addgeo:deepthreshold} with $t=r\ge1$ to bound this
deep part by $Bg^{-r}$. Adding the near sum, absorbing $r!$
into $[Cr]^{Cr}$, and inserting the two boundary factors
$e^{x/2}$ proves \eqref{addgeo:separatingtail}.
\end{proof}

\begin{proposition}[Kernel errors at an additive threshold]
\label{addgeo:kernels}
The estimates \eqref{allrem:nonsep}--\eqref{allrem:sep} hold under
$g\ge Cs^2+[Cr]^Cs$, where $r=R+1$.
For a fixed carrier and pattern, \eqref{allrem:kernel} holds under
$g\ge Cs^2+[Ct]^Cs$, where $t=R+m+1$.
Here $C$ is a new absolute constant; all coefficient functions,
error powers and length factors are unchanged.
\end{proposition}

\begin{proof}
We reuse the retained-record calculations from the proof of
Proposition~\ref{allrem:kernels}, whose local thresholds are linear
in $s$, and replace its tail input by Lemma~\ref{addgeo:tails}.

For the nonseparating kernel, Lemma~\ref{allrem:boundary}, a genus
shift by one and the scalar--boundary product argument apply under
$g\ge[Cr]^Cs$. They give \eqref{allrem:nonsep} for all $R\ge0$.
For a separating endpoint of deficit $1\le\delta\le R$, use
Lemma~\ref{allrem:shift} with $J=R-\delta$ and
$H_0=h_0\le(\delta+1)/2$, so $L=R+h_0+1\le2r$.
Its threshold is again linear in $s$. Multiplication by the same
small-component prefactor gives the error in \eqref{allrem:sep};
the omitted endpoints are covered by \eqref{addgeo:separatingtail}.

For a retained carrier record $\delta\le R$, put $J=R-\delta$.
The exact record identities imply
\begin{equation}\label{addgeo:recordthresholds}
 \begin{aligned}
 &H_0\le(t+1)/2\le t,\qquad
 L=J+H_0+\delta+1=R+H_0+1\le2t,\\
 &p=n-A_{\rm out}+k_*\ge k_*\ge1.
 \end{aligned}
\end{equation}
All scalar marked counts $p_\nu$ are nonnegative and at most
$n+2H_0+\delta\le n+3t$.
For every source order $J\le t$, the exact threshold in
Proposition~\ref{vol:hma} is consequently bounded by
\[
 G_J(p_\nu)
 \le3000(2c_0)^{c_0}t^{c_0+1}s.
\]
The graded majorants, genus shifts and at most $2L$ scalar
factors in \eqref{vol:telescoping} therefore have a common base
$A=[C_1L]^{C_1}s$, with one auxiliary absolute $C_1$.
This is the base used in Lemma~\ref{allrem:shift};
neither finite multiplication nor shifting raises the power of $s$.
The geometric denominator coefficients are covered by the same
base, including negative values of $n-\delta+q-2$.

Choose $C_1$ after the fixed enlargements in that lemma.
Requiring $g\ge4A$ ensures that the scalar partial sums are bounded,
that $g-H_0\ge g/2$, that every strict source threshold
$g-H_0>G_J(p_\nu)$ holds, and that $g-H_0\ge4$.
The last condition ensures that the predecessor genera in the source
estimates are at least three.
The exact boundary ratio is bounded by $e^{X_*/2}$, so the
scalar smallness condition is independent of boundary length.
Thus \eqref{allrem:shifterror} applies directly.

Multiply it by \eqref{allrem:recordprefactor}.
The single label factor gives
$s^{A_{\rm out}+J+1}\le s^{R+2}$;
the boundary exponentials give $e^\Lambda$;
and $k+2(R-\delta)+2\le2t$ because $\delta\ge m$.
The record count and exact small-volume factors cost $[Ct]^{Ct}$.
Adding the fixed-pattern boundary tail of
Lemma~\ref{addgeo:tails} proves \eqref{allrem:kernel}.

To verify all thresholds simultaneously, choose the final $C$
after $K,C_1$ and the constants of the local graded lemmas.
For every $u\ge1$ it may be chosen to satisfy
\begin{equation}\label{addgeo:domination}
 C\ge K,\qquad [Cu]^C\ge Ku^2,\qquad
 [Cu]^C\ge4[2C_1u]^{C_1},\qquad
 [Cu]^C\ge8(u+2),
\end{equation}
and to dominate the other finitely many fixed polynomial
requirements in $u$.
With $u=t$ or $r$, the additive threshold supplies both the
new tail threshold and every local condition above.
It also gives $N>2R$ and $g>H_0$, so retained records reconstruct
a unique largest stable complementary component.
Unavailable cusp allocations are omitted before any volume is
evaluated, and the available boundary ratio has $p\ge1$.
The scalar paths alone may use zero marked points.

The finite retained sums are exactly
\eqref{vol:canonicalcoefficients}.
When $R=0$ or $m>R$ they are empty and the new positive tail
supplies the full error.
For $R=0$, the separating estimate likewise follows entirely from
the positive tail bound, so no shifted estimate of negative order is needed.
\end{proof}

\begin{theorem}[All-order orbital remainder at an additive threshold]
\label{addgeo:orbital}
There is an absolute $C\ge8$ such that, for all integers $R,n\ge0$,
$T\ge0$, and bounded real even Borel $H$ supported in $[-T,T]$,
putting $\rho=R+\lceil T\rceil+1$, the condition
\begin{equation}\label{addgeo:orbitalthreshold}
 g\ge C(n+1)^2+[C\rho]^C(n+1)
\end{equation}
implies
\begin{equation}\label{addgeo:orbitalbound}
 \left|\E_{g,n}\mathcal O_H-\Omega(H)
       -\sum_{j=1}^Rg^{-j}\kappa_{j,n}(H)\right|
 \le g^{-R-1}(n+1)^{R+2}
       [C\rho]^{C\rho}e^{3T}\|H\|_\infty.
\end{equation}
The sum is empty for $R=0$.
The coefficients are those of Theorem~\ref{orb:expansion}.
\end{theorem}

\begin{proof}
Enlarge $C$ beyond the constant of Proposition~\ref{addgeo:kernels}
so that \eqref{addgeo:orbitalthreshold}, with $s,\rho\ge1$,
also implies $g\ge4\rho^2$. Since $T\le\rho-1$,
$\pi(g-1+n/2)>T$; thus \eqref{orb:wholecarrier} excludes every
carrier without a stable exterior, including entirely self-paired carriers.
Put $M=\lfloor2T/\pi\rfloor$. The exhaustion \eqref{orb:exhaustion}
has $1\le m\le M\le T$, so $r=R+1\le\rho$ and
$t=R+m+1\le\rho$. Proposition~\ref{addgeo:kernels} therefore supplies
\eqref{allrem:nonsep}--\eqref{allrem:kernel} under the common additive
threshold, with unchanged coefficients, length factors and cusp powers.
Repeat the integration and finite summation in the proof of
Theorem~\ref{allrem:orbital}, using these kernel bounds.
The absolute repetition estimate following \eqref{orb:Loperation},
together with \eqref{orb:fiberbound} and \eqref{orb:simplex}, gives the
same error bound and absolute integrability for bounded Borel $H$;
hence signed Fubini applies without another genus condition.
The two vanishing properties in \eqref{orb:xivanishing} identify the nonsimple
coefficient sum with $\zeta_{j,n}$, even when $j>M$ or $M>R$.
Adding the simple part gives $\kappa_{j,n}=\eta_{j,n}+\zeta_{j,n}$
and proves \eqref{addgeo:orbitalbound}, including $R=0$ with empty sums.
If $M=0$, the nonsimple contribution and all its coefficients on $H$
vanish; if $T=0$, all terms vanish, regardless of $H(0)$.
\end{proof}

The coefficient of $(n+1)^2$ in the genus threshold is now
independent of order and support.
The remaining order-dependent term is linear in $n+1$.
Section~\ref{sec:additivepolynomial} transfers this estimate to
polynomials of growing degree; Section~\ref{sec:additivespectral}
uses that transfer to bound the spectral count.

\section{The leading measure and a common truncated trace}
\label{sec:trace}

Positive pre-trace estimates on domains truncated at the cusps
were used by Hide \cite[Sections~2.4--2.5]{Hide23}.
We derive an exact trace formula at one cusp height for all
polynomial tests of bounded degree. The leading orbital functional
is represented by a unit atom at $a$ and a finite signed measure on
$[0,1]$. Subtracting the constant eigenfunction cancels the atom up
to the explicit mass removed by cusp truncation.

Throughout, $h=b_0*b_0$, $f=\widehat h$, $\supp h\subset[-L,L]$,
and $a=f(i/2)>1$ are fixed as in Section~\ref{sec:quantitative}.
The degree bound $Q\ge1$ is an integer, and polynomial tests are real
and vanish at zero. We use
\begin{equation}\label{trace:fourier}
 \widehat H(r)=\int_{\R}H(u)e^{-iru}\,du,
 \qquad
 \check F(u)=\frac1{2\pi}\int_{\R}F(r)e^{iru}\,dr.
\end{equation}
In particular, $0\le f(r)\le1$ on $\R$, and
$1\le f(is)\le a$ for $0\le s\le1/2$.  The latter assertion follows
from $f(is)=\int h(u)\cosh(su)\,du$ and $h\ge0$.
For a polynomial $p$ with $p(0)=0$, put
\begin{equation}\label{trace:polynomial-test}
 H_p=\check{(p\circ f)},\qquad
 p(x)=\sum_{t=1}^{Q}c_tx^t
 \quad\Longrightarrow\quad
 H_p=\sum_{t=1}^{Q}c_th^{*t},\qquad
 \supp H_p\subset[-LQ,LQ].
\end{equation}

\subsection{A finite signed measure for the leading functional}
\label{sec:leadingmeasure}

The separation of the term $p(a)$ also occurs in the proof of
\cite[Lemma~3.3]{HMT25}. The argument below represents the remainder
by a finite signed measure on $[0,1]$ and bounds its total variation.

\begin{proposition}[Leading orbital measure]\label{trace:leading}
There is a finite signed real measure $\mu_0$, independent of $n$, such that
for every polynomial $p$ vanishing at zero,
\begin{equation}\label{trace:leading-representation}
 \begin{gathered}
 \omega_{1,n}(p)=\Omega(H_p)=p(a)+\int_{[0,1]}p(x)\,d\mu_0(x),\\
 \supp\mu_0\subset[0,1],\qquad
 \mu_0(\R)=-1,\qquad \|\mu_0\|_{\TV}<5.
 \end{gathered}
\end{equation}
Consequently $\omega_{1,n}=\omega_1$ is independent of the cusp number.
The measure $\mu_f=\delta_a+\mu_0$ has mass zero and variation less than
six; its atom at $a$ has coefficient exactly one.
\end{proposition}

\begin{proof}
The integral in \eqref{orb:etadefinition} is
\[
 \Omega(H)=\int_0^\infty\sum_{k\ge1}
 \frac{2\sinh^2(x/2)}{\sinh(kx/2)}H(kx)\,dx.
\]
Changing variables $u=kx$ gives
\begin{align}
 \Omega(H)&=\int_0^\infty A(u)H(u)\,du,\notag\\
 A(u)&=\sum_{k\ge1}\frac2k
             \frac{\sinh^2(u/(2k))}{\sinh(u/2)}
       =2\cosh(u/2)+R(u),\label{trace:leading-weight}\\
 R(u)&=-2e^{-u/2}
       +\sum_{k\ge2}\frac2k
             \frac{\sinh^2(u/(2k))}{\sinh(u/2)},
 \qquad R(0)=-2.\label{trace:remainder-weight}
\end{align}
Indeed, $\sinh^2v\le\frac12\sinh(2v)$ and convexity of $\sinh$
imply $\sinh^2(u/(2k))\le k^{-1}\sinh(u/2)$ for $k\ge2$.
Thus these summands are bounded by $2/k^2$; their individual limits at
zero vanish.  This proves continuity of $R$ at zero and justifies the
change of variables for signed compactly supported $H$ by an integrable
absolute majorant.

Let $\mathsf P_s$ be the Cauchy probability measure
\begin{equation}\label{trace:cauchy}
 \mathsf P_s(dr)=\frac{s}{\pi(s^2+r^2)}\,dr\quad(s>0),\qquad
 \mathsf P_0=\delta_0,\qquad
 \int_{\R}e^{-iru}\,d\mathsf P_s(r)=e^{-s|u|}.
\end{equation}
The expansion of the hyperbolic denominator leads to symmetric
second differences of the Cauchy measures. Set $s_j=j+1/2$,
$d_k=1/k$, and consider
\begin{equation}\label{trace:rho-series}
 \rho=-2\mathsf P_{1/2}
 +\sum_{k=2}^{\infty}\frac1k\sum_{j=0}^{\infty}
       \bigl(\mathsf P_{s_j-d_k}+\mathsf P_{s_j+d_k}-2\mathsf P_{s_j}\bigr).
\end{equation}
We first prove convergence in total variation, including the endpoint
$s_0=d_2$.  Differentiation of the density gives
\[
 \partial_s^2\mathsf P_s(dr)=\frac{2s(s^2-3r^2)}{\pi(s^2+r^2)^3}\,dr,
 \qquad \|\partial_s^2\mathsf P_s\|_{\TV}\le\frac{3}{2s^2}.
\]
For the inequality, replace $|s^2-3r^2|$ by $s^2+3r^2$, scale $r=sv$,
and use $\int_{\R}(1+v^2)^{-3}dv=3\pi/8$ and
$\int_{\R}v^2(1+v^2)^{-3}dv=\pi/8$.  The integral remainder for a
symmetric second difference in $L^1(\R)$ therefore gives, for $0<d<s$,
\begin{equation}\label{trace:second-difference}
 \|\mathsf P_{s-d}+\mathsf P_{s+d}-2\mathsf P_s\|_{\TV}
 \le\int_{-d}^{d}(d-|v|)\frac{3}{2(s+v)^2}\,dv
 \le\frac{3d^2}{2(s-d)^2}.
\end{equation}
When $k\ge3$, $s_j-d_k\ge s_j/3$; when $k=2,j\ge1$,
$s_j-1/2\ge2s_j/3$.  Since
$\sum_{j\ge0}(j+1/2)^{-2}=\pi^2/2$, their respective contributions to
the variation sum are at most
\[
 \frac{27\pi^2}{4}\left(\zeta(3)-\frac98\right),
 \qquad
 \frac{27}{64}\left(\frac{\pi^2}{2}-4\right).
\]
The term $(k,j)=(2,0)$ has variation at most $2$, and the initial
term $-2\mathsf P_{1/2}$ has variation $2$.  Hence
\begin{equation}\label{trace:rho-variation}
 \frac12\|\rho\|_{\TV}
 \le 2+\frac{27\pi^2}{8}\left(\zeta(3)-\frac98\right)
       +\frac{27}{128}\left(\frac{\pi^2}{2}-4\right)<5.
\end{equation}
To verify the last inequality using rational bounds, the integral test gives
$\zeta(3)\le\sum_{k=1}^{10}k^{-3}+1/200<1203/1000$; together with
$\pi^2<10$, the right side is less than $15499/3200<5$.

For $u>0$, expansion of $1/\sinh(u/2)$ gives
\begin{equation}\label{trace:rho-transform-computation}
 \frac2k\frac{\sinh^2(u/(2k))}{\sinh(u/2)}
 =\frac1k\sum_{j\ge0}
 \bigl(e^{-(s_j-d_k)u}+e^{-(s_j+d_k)u}-2e^{-s_ju}\bigr).
\end{equation}
The variation convergence just proved permits termwise Fourier
transformation of \eqref{trace:rho-series}.  Equations
\eqref{trace:cauchy} and \eqref{trace:rho-transform-computation} imply
\begin{equation}\label{trace:rho-transform}
 \int_{\R}e^{-iru}\,d\rho(r)=R(|u|),\qquad \rho(\R)=-2.
\end{equation}
The mass also follows directly because every second difference has
mass zero.  Define
\begin{equation}\label{trace:pushforward}
 \mu_0=\tfrac12 f_*\rho,\qquad \mu_f=\delta_a+\mu_0.
\end{equation}
The real map $f$ takes values in $[0,1]$, and pushforward does not
increase variation.  Thus the asserted support, mass and norm follow.

Finally, evenness and compact support of $H_p$ give
\begin{align}
 \Omega(H_p)
 &=\int_0^\infty2\cosh(u/2)H_p(u)\,du
       +\frac12\int_{\R}R(|u|)H_p(u)\,du\notag\\
 &=p(f(i/2))+\frac12\int_{\R}p(f(r))\,d\rho(r)
  =p(a)+\mu_0(p).\label{trace:leading-inversion}
\end{align}
The first equality in the second line uses the entire Fourier transform
of $H_p$; the second term uses Fubini, justified by $H_p\in L^1$ and
$\rho$ finite.  Because $a>1$, $f_*\rho$ has no mass at $a$.  Thus $\mu_f(\{a\})=1$.
\end{proof}

The endpoint term $(k,j)=(2,0)$ in \eqref{trace:rho-series}
contributes $\frac14\delta_{f(0)}=\frac14\delta_1$ to $\mu_0$.
Since $h(0)=\int_{\R}b_0(u)^2\,du>0$, the function $h$ is positive
near zero, so $f(r)<1$ for $r\ne0$. The measure
$\rho-\frac12\delta_0$ is absolutely continuous by total-variation
convergence. Hence $\mu_0(\{1\})=1/4$, although its total mass is $-1$.

The functional
$q\mapsto\int q\,d\mu_f$ extends $\omega_1$ continuously to
$C([0,a])$. Two finite signed measures on $[0,a]$ agreeing on
polynomials vanishing at zero differ by a multiple of $\delta_0$,
by polynomial density. Requiring total mass zero fixes this ambiguity.

\subsection{Abel kernels and localized spectral multipliers}
\label{sec:abelkernels}

Let $X=\Gamma\backslash\mathbb H$ be a connected complete hyperbolic
surface of finite area, with $n\ge1$ cusps and no orbifold points,
where $\Gamma\subset\mathrm{PSL}_2(\R)$ is torsion-free.
We use $d\mu=dx\,dy/y^2$ and
$\Delta=-y^2(\partial_x^2+\partial_y^2)$.  Its area is
$\mathcal A_{g,n}=2\pi(2g-2+n)$. For each cusp $b$, let $\Gamma_b$
be its full stabilizer and choose a scaling matrix $\sigma_b$ such
that $\sigma_b^{-1}\Gamma_b\sigma_b=\langle U\rangle$, where
$U:z\mapsto z+1$.

We use the Abel transform to pass from the one-dimensional test $H$
to a radial kernel on $\mathbb H$. For $T>0$ and even
$H\in C_c^\infty(\R)$ supported in $[-T,T]$, define,
for $\rho\ge0$,
\begin{equation}\label{trace:abel-inverse}
 K_H(\rho)=-\frac1{\sqrt2\pi}\int_\rho^\infty
          \frac{H'(u)}{\sqrt{\cosh u-\cosh\rho}}\,du.
\end{equation}
This agrees with the normalization in
\cite[\S5, (5.3)--(5.4)]{HWX26}. The kernel may change sign;
positivity of $A_X$ below will come from the fixed multiplier $f$.
To check smoothness at the origin, the function
$G(v)=H'(\operatorname{arcosh}v)/\sqrt{v^2-1}$ extends smoothly to
$v=1$: $H'(u)/\sinh u$ is smooth and even near zero.  Changing variables
in \eqref{trace:abel-inverse} gives
\[
 K_H(\rho)=-\frac1{\sqrt2\pi}
       \int_0^\infty G(\cosh\rho+v)v^{-1/2}\,dv.
\]
Differentiation with respect to $\cosh\rho$ is justified by the
smoothness and compact support of $G$ and the integrability of
$v^{-1/2}$ near zero.  Thus $K_H$ is a smooth
radial kernel supported on $[0,T]$, flat at the upper endpoint.
Its forward Abel transform is
\begin{equation}\label{trace:abel-forward}
 H(u)=\sqrt2\int_{|u|}^{\infty}
       \frac{K_H(\rho)\sinh\rho}
            {\sqrt{\cosh\rho-\cosh u}}\,d\rho.
\end{equation}
Indeed, inserting \eqref{trace:abel-inverse} and interchanging integrals
produces the inner integral
\[
 \int_u^v\frac{\sinh\rho\,d\rho}
 {\sqrt{\cosh\rho-\cosh u}\sqrt{\cosh v-\cosh\rho}}=\pi
 \quad(0\le u<v);
\]
the result is $-\int_u^\infty H'(v)dv=H(u)$.  Absolute integrability
follows by the same computation with $|H'|$.

The spectral multiplier of this kernel is
$\widehat H(r)$.  At the point $i$, write $y=e^u$ and $x=e^{u/2}v$.
Then $\cosh d(i,x+iy)=\cosh u+v^2/2$, and convolution of $K_H$ with
$y^{1/2+ir}$ is
\[
 \int_{\R} e^{iru}\left[\int_{\R}
 K_H\!\left(\operatorname{arcosh}(\cosh u+v^2/2)\right)dv\right]du
 =\int_{\R}e^{iru}H(u)du=\widehat H(r).
\]
The bracket is exactly \eqref{trace:abel-forward}.
Also
\begin{equation}\label{trace:identity-density}
 K_H(0)=\frac1{4\pi}\int_{\R}r\tanh(\pi r)\widehat H(r)\,dr.
\end{equation}
This follows from
$K_H(0)=-(2\pi)^{-1}\int_0^\infty H'(u)/\sinh(u/2)\,du$,
Fourier inversion, and
$\int_0^\infty\sin(ru)/\sinh(u/2)du=\pi\tanh(\pi r)$.
The last identity follows by expanding $1/\sinh(u/2)$ into exponentials
and using the partial-fraction expansion of $\tanh$; the resulting
series is $2\sum_{j\ge0}r/((j+1/2)^2+r^2)$.
All Fourier interchanges are justified by rapid decay and
$|\sin(ru)|\le\min\{1,|ru|\}$.
This also checks the spherical inversion normalization recorded in
\cite[\S5, immediately before (5.8)]{HWX26}.

The periodization
\begin{equation}\label{trace:periodization}
 \mathcal K_H(z,w)=\sum_{\gamma\in\Gamma}K_H(d(z,\gamma w))
\end{equation}
is locally finite on compact coordinate lifts, by proper discontinuity
and compact radial support. Choose an orthonormal basis $\{u_j\}$ of
the discrete spectral subspace of $L^2(X)$, consisting of the constant,
residual and cuspidal eigenfunctions. Write $\Delta u_j=\lambda_j u_j$ and
$\lambda_j=1/4+r_j^2$, with $r_j\ge0$ or $r_j\in i[0,1/2]$.
The Eisenstein series $E_b$ use the width-one cusp coordinates.
The local spectral expansion is
\begin{equation}\label{trace:local-spectral}
 \mathcal K_H(z,w)=
 \sum_j\widehat H(r_j)u_j(z)\overline{u_j(w)}
 +\frac1{4\pi}\sum_{b=1}^{n}\int_{\R}
  \widehat H(r)E_b(z,\tfrac12+ir)
                 \overline{E_b(w,\tfrac12+ir)}\,dr.
\end{equation}
The precise input is
\cite[\S3.3, (3.11)]{Bruin13}, with transform hypotheses in \S3.1
and the convergence bound in \S3.5, Lemma~3.2.  The transform is even and
entire, with arbitrarily rapid decay in every fixed horizontal strip;
thus those hypotheses hold. Our nonnegative Laplacian is Bruin's
$-\Delta$. On compact coordinate lifts, Lemma~3.2 gives a uniform
linear bound in the eigenvalue cutoff for the local spectral mass.
Rapid decay of the multiplier and Cauchy--Schwarz therefore give
locally uniform absolute convergence, permitting diagonal restriction
and integration over compact sets.

We specify the group convention in applying this input.  Since $n\ge1$,
$\Gamma$ is free of rank $2g+n-1$.  Lift free generators to
$\mathrm{SL}_2(\R)$.  The resulting representation projects to the
original faithful representation, so its image is an injective lift
omitting $-I$.  It is discrete since its projected image is discrete.
Bruin's integration divisor $\#(\Gamma\cap\{\pm I\})$ is therefore
one.  Weight-zero kernels and cusp coordinates agree with those in
$\mathrm{PSL}_2(\R)$; the sign of a lifted parabolic does not affect
its action. Bruin's transform is $\widehat H((s-1/2)/i)$.
With $s=1/2+ir$, its symmetry under $s\mapsto1-s$ makes the
choice of branch for an exceptional eigenvalue immaterial.
No factor two is introduced in \eqref{trace:local-spectral}.

We now return to the fixed multiplier $f$. Define the bounded
self-adjoint operator $A_X$ by its Laplace multiplier
\begin{equation}\label{trace:operator}
 f(r(\lambda)),\qquad
 r(\lambda)=
 \begin{cases}
  \sqrt{\lambda-1/4},&\lambda\ge1/4,\\
  i\sqrt{1/4-\lambda},&0\le\lambda<1/4.
 \end{cases}
 \qquad 0\le A_X\le a.
\end{equation}
This is the meaning of $f(\sqrt{\Delta_X-1/4})$. The constant
eigenfunction has $A_X$-eigenvalue $a$; the part of the spectrum of
$\Delta_X$ in $[1/4,\infty)$ maps into $[0,1]$.
For $p=x\widetilde p$, we use the two norms
\begin{equation}\label{trace:polynomial-norms}
 N_*(p)=\|\widetilde p\|_{C([0,a])},\qquad
 N_{\R}(p)=\|\widetilde p\|_{C([0,1])}.
\end{equation}
The first norm controls the full spectral trace. The second will
control the identity and parabolic terms using only real spectral
parameters.

\subsection{Uniform geometry of all cusp lifts}
\label{sec:cusplifts}

The following cusp neighborhoods have the same height normalization
on every surface. Control of all their lifts is needed when unfolding
the parabolic contribution: removing other cusp neighborhoods must
not remove any part of the translation kernel's support.

\begin{lemma}[Uniform cusp truncation]\label{trace:cusp-lifts}
For $Y>1$ the width-one cusp regions $C_b(Y)$, represented by
$0\le x<1$, $y>Y$, are embedded and pairwise disjoint, each with area
$1/Y$.  In the chart of a distinguished cusp, every other lift of every
$C_b(Y)$ lies below height $1/Y$.
\end{lemma}

\begin{proof}
We first prove the normalized Shimizu estimate.  In a discrete subgroup
containing $U=\left(\begin{smallmatrix}1&1\\0&1\end{smallmatrix}\right)$,
an element $S=\left(\begin{smallmatrix}a&b\\c&d\end{smallmatrix}\right)$
has either $c=0$ or $|c|\ge1$. Otherwise set $S_0=S$ and put
$S_{m+1}=S_mUS_m^{-1}$.
Direct multiplication gives
\begin{equation}\label{trace:shimizu}
 c_{m+1}=-c_m^2,\quad a_{m+1}=1-a_mc_m,
 \quad b_{m+1}=a_m^2,\quad d_{m+1}=1+a_mc_m.
\end{equation}
If $0<|c_0|<1$, then $|c_m|=|c_0|^{2^m}\to0$.
Once $|c_m|\le1/2$, the recurrence for $a_m$ bounds it uniformly;
it then implies $a_{m+1}\to1$.  Hence $S_m\to U$ through distinct
elements, contradicting discreteness.  The same proof applies to
matrix representatives of a subgroup of $\mathrm{PSL}_2(\R)$.

Let $\sigma_b,\sigma_d$ be two width-one scalings and
$M=\sigma_b^{-1}\gamma\sigma_d$
with lower-left entry $c$.  If $c=0$, the corresponding cusp lift is
the distinguished lift: equality of its fixed point with infinity
means that the cusp representatives agree and the element belongs
to the stabilizer.  Otherwise $MUM^{-1}$ belongs to
$\sigma_b^{-1}\Gamma\sigma_b$ and has lower-left entry $-c^2$.
The just-proved estimate gives $|c|\ge1$.  For $w=u+iv$, $v>Y$,
\begin{equation}\label{trace:cross-cusp-height}
 \operatorname{Im}(Mw)=\frac{v}{|cw+d|^2}
       \le\frac1{c^2v}<\frac1Y.
\end{equation}
This proves the assertion for different lifts of the same cusp and
for lifts of different cusps.  Since $Y>1$, they cannot meet the
distinguished region $y>Y$.  The area of its quotient strip is
$\int_0^1\int_Y^\infty y^{-2}dy\,dx=1/Y$.
\end{proof}

\subsection{The exact trace on a common truncation}
\label{sec:commontrace}

The complement $X(Y)=X\setminus\bigcup_bC_b(Y)$ is compact. Write
$\chi_Y=\ind_{X(Y)}$, as a multiplication operator, and fix
\begin{equation}\label{trace:common-cutoff}
 Y_Q=3e^{LQ}\qquad(Q\ge1).
\end{equation}
Every test $H_p$ with $\deg p\le Q$ is supported in $[-LQ,LQ]$.
The common cutoff therefore gives one geometric trace formula and
one positive weighted spectral measure for all these tests.
Individual geometric terms may still have either sign.

\begin{theorem}[Exact trace on a common truncation]\label{trace:exact}
For every such surface $X$, integer $Q\ge1$ and real polynomial
$p$ of degree at most $Q$ with $p(0)=0$, the compression
$\chi_{Y_Q}p(A_X)\chi_{Y_Q}$ is trace class. For $Y>1$, define
\begin{align}
 \mathcal I(p)&=K_{H_p}(0)
   =\frac1{4\pi}\int_{\R}r\tanh(\pi r)p(f(r))\,dr,
       \label{trace:I-definition}\\
 \mathcal P_Y(p)&=2\int_0^{LQ}K_{H_p}(\rho)\cosh(\rho/2)
       \mathsf H_{\lfloor2Y\sinh(\rho/2)\rfloor}\,d\rho,
       \label{trace:P-definition}
\end{align}
where $\mathsf H_m=\sum_{k=1}^mk^{-1}$ for integers $m\ge1$,
and $\mathsf H_0=0$.
Then the exact identity is
\begin{equation}\label{trace:exact-identity}
 \tau_{X,Q}(p):=\Tr\bigl(\chi_{Y_Q}p(A_X)\chi_{Y_Q}\bigr)
 =\left(\mathcal A_{g,n}-\frac n{Y_Q}\right)\mathcal I(p)
       +n\mathcal P_{Y_Q}(p)+\mathcal O_X(p).
\end{equation}
If $p$ is nonnegative on the spectrum of $A_X$, then
$\tau_{X,Q}(p)\ge0$.  There is a finite positive measure
$\nu_{X,Q}$ supported in $[0,a]$ such that
\begin{equation}\label{trace:weighted-measure}
 \begin{gathered}
 \tau_{X,Q}(x\widetilde p)=\int_{[0,a]}\widetilde p(s)\,d\nu_{X,Q}(s),\\
 \nu_{X,Q}([0,a])=\tau_{X,Q}(x),\qquad
 |\tau_{X,Q}(p)|\le N_*(p)\tau_{X,Q}(x).
 \end{gathered}
\end{equation}
\end{theorem}

\begin{proof}
\emph{Trace class and positivity.}
Apply \eqref{trace:local-spectral} to $H=h$.  All diagonal terms are
nonnegative, since the multiplier $f$ is nonnegative on the full
spectrum.  For any $Y>1$, the diagonal kernel is smooth on the compact
set $X(Y)$.  With $S=A_X^{1/2}\chi_Y$, spectral Parseval and Tonelli give
\begin{align}
 \|S\|_{\mathrm{HS}}^2
 &=\sum_j f(r_j)\|\chi_Yu_j\|_2^2
  +\frac1{4\pi}\sum_{b=1}^n\int_{\R}f(r)
       \|\chi_YE_b(\cdot,\tfrac12+ir)\|_2^2\,dr\notag\\
 &=\int_{X(Y)}\mathcal K_h(z,z)\,d\mu(z)<\infty.
       \label{trace:HS}
\end{align}
Explicitly, for an orthonormal basis $\{e_m\}$ of $L^2(X)$, Parseval
turns $\sum_m|\langle\chi_Ye_m,u_j\rangle|^2$ into
$\|\chi_Yu_j\|_2^2$, and does the same for each Eisenstein function
restricted to $X(Y)$.  Nonnegativity justifies interchanging this
basis sum and the spectral integrals.
Since $p=x\widetilde p$, functional calculus now gives
\begin{equation}\label{trace:factorization}
 \chi_Yp(A_X)\chi_Y=S^*\widetilde p(A_X)S,
 \qquad
 \|\chi_Yp(A_X)\chi_Y\|_1\le N_*(p)\|S\|_{\mathrm{HS}}^2.
\end{equation}
Thus $\chi_Yp(A_X)\chi_Y$ is trace class, also for complex polynomials.
It is positive when $p$ is nonnegative on the spectrum of $A_X$.
The inequality $|p(s)|\le N_*(p)s$ on $[0,a]$ dominates the signed
spectral diagonal by \eqref{trace:HS}.  Thus at $Y=Y_Q$,
\begin{align}
 \tau_{X,Q}(p)
 &=\int_{X(Y_Q)}\mathcal K_{H_p}(z,z)\,d\mu(z)\notag\\
 &=\sum_jp(f(r_j))\|\chi_{Y_Q}u_j\|_2^2
 +\frac1{4\pi}\sum_{b=1}^{n}\int_{\R}p(f(r))
       \|\chi_{Y_Q}E_b(\cdot,\tfrac12+ir)\|_2^2\,dr,
       \label{trace:spectral-trace}
\end{align}
with absolute convergence. With $S=A_X^{1/2}\chi_{Y_Q}$, define
$\nu_{X,Q}(B)=\|\ind_B(A_X)S\|_{\mathrm{HS}}^2$ for Borel
sets $B\subset[0,a]$. This proves
\eqref{trace:weighted-measure}; countable additivity follows from
orthogonality of the spectral projections.

\emph{Hyperbolic contribution.}
We compute the geometric integral for any even
$H\in C_c^\infty(\R)$ supported in $[-T,T]$, with $T>0$, provided
\begin{equation}\label{trace:cutoff-conditions}
 Y>\max\{1,e^{T/2},2\sinh(T/2)\}.
\end{equation}
The condition $Y>e^{T/2}$ ensures that hyperbolic terms vanish in
the removed cusp regions. The other condition, $Y>2\sinh(T/2)$,
will make the parabolic unfolding exact.
A hyperbolic element does not stabilize the distinguished cusp:
the cusp stabilizer is the primitive parabolic cyclic group.  Indeed, an
additional dilation fixing that cusp would conjugate its translations
to translations tending to the identity, contradicting discreteness.
Thus its lower-left entry satisfies $|c|\ge1$ and, for $z=x+iy$,
\begin{equation}\label{trace:hyperbolic-displacement}
 \operatorname{Im}(\gamma z)\le\frac1y,
 \qquad d(z,\gamma z)\ge
 \log\frac{y}{\operatorname{Im}(\gamma z)}\ge2\log y.
\end{equation}
Every hyperbolic kernel term vanishes on every deleted cusp tail.
Its truncated integral therefore equals its full integral.  The
periodization of $|K_H|$ is locally bounded on the compact truncation,
so these hyperbolic integrals are absolutely summable.

For a primitive oriented hyperbolic class $\gamma$ of length $\ell$,
the centralizer of every nonzero power is $\langle\gamma\rangle$.
Fix a positive repetition index $k\ge1$ and put $v=k\ell$.
Fermi coordinates $(s,r)$ on its cyclic quotient satisfy
$0\le s<\ell$, $d\mu=\cosh r\,ds\,dr$, and
\[
 \cosh d(z,\gamma^kz)
       =\cosh v+2\sinh^2(v/2)\sinh^2r.
\]
Changing variables from $r$ to $\rho=d(z,\gamma^kz)$, including both
signs of $r$, and applying \eqref{trace:abel-forward}, yields
\begin{align}
 \int_{\langle\gamma\rangle\backslash\mathbb H}
       K_H(d(z,\gamma^kz))\,d\mu(z)
 &=\frac{\ell}{\sqrt2\sinh(v/2)}\int_v^\infty
       \frac{K_H(\rho)\sinh\rho}
            {\sqrt{\cosh\rho-\cosh v}}\,d\rho\notag\\
 &=\frac{\ell}{2\sinh(k\ell/2)}H(k\ell).
       \label{trace:hyperbolic-orbit}
\end{align}
The inverse class is distinct: an orientation-preserving isometry
interchanging the endpoints of the axis has order two in
$\mathrm{PSL}_2(\R)$ and cannot belong to the torsion-free group.
Consequently conjugacy-class unfolding gives precisely
$\mathcal O_H(X)$, with both primitive orientations and every positive
repetition.

\emph{Parabolic contribution.}
Fix a cusp chart. The displacement of
its translation by $k\in\Z\setminus\{0\}$ is
\begin{equation}\label{trace:translation-displacement}
 d(z,z+k)=2\operatorname{arsinh}\frac{|k|}{2y}.
\end{equation}
Unfolding the conjugates of its stabilizer gives the width-one strip
with the lifted truncation indicator.  Apart from its own deleted
region $y>Y$, this strip has holes from other removed cusp lifts.
By Lemma~\ref{trace:cusp-lifts}, every such hole lies in $y<1/Y$.
On the support of the translation kernel, however,
\begin{equation}\label{trace:no-holes}
 y\ge\frac{|k|}{2\sinh(T/2)}
       \ge\frac1{2\sinh(T/2)}>\frac1Y.
\end{equation}
The translation kernel vanishes on every other removed cusp lift.
Hence replacing the unfolded domain by the whole strip
$0\le x<1$, $0<y\le Y$ leaves the integral unchanged. This is a
support argument and does not require the kernel to be nonnegative.
Retaining both signs of $k$ and substituting
$y=k/(2\sinh(\rho/2))$ gives the contribution of this one cusp:
\begin{align}
 \mathcal P_Y(H)
 &=\sum_{k\ne0}\int_0^Y
       K_H\!\left(2\operatorname{arsinh}\frac{|k|}{2y}\right)
       \frac{dy}{y^2}\notag\\
 &=\sum_{k\ge1}\frac2k
       \int_{2\operatorname{arsinh}(k/(2Y))}^{T}
           K_H(\rho)\cosh(\rho/2)\,d\rho\notag\\
 &=2\int_0^TK_H(\rho)\cosh(\rho/2)
           \mathsf H_{\lfloor2Y\sinh(\rho/2)\rfloor}\,d\rho.
       \label{trace:parabolic-unfolding}
\end{align}
An integral in the middle line is zero if its lower limit exceeds
$T$.  Only $k\le2Y\sinh(T/2)$ contribute, so the sum is finite and
the integrals are absolutely convergent.  This argument covers other lifts of the
same cusp and lifts of different cusps, and gives the identical term
for each of the $n$ cusps.

\emph{Assembly of the trace formula.}
The identity element contributes
$(\mathcal A_{g,n}-n/Y)K_H(0)$.  Nonidentity elements are hyperbolic
or parabolic, since there are no orbifold points.  Thus
\begin{equation}\label{trace:geometric-identity}
 \int_{X(Y)}\mathcal K_H(z,z)\,d\mu(z)
  =\left(\mathcal A_{g,n}-\frac nY\right)K_H(0)
        +n\mathcal P_Y(H)+\mathcal O_H(X).
\end{equation}
For $T=LQ$ and $Y=3e^{LQ}$, all inequalities in
\eqref{trace:cutoff-conditions} hold.  Taking $H=H_p$ proves
\eqref{trace:exact-identity}.
\end{proof}

The measure in \eqref{trace:weighted-measure} represents the
\emph{weighted} functional $\widetilde p\mapsto\tau_{X,Q}(x\widetilde p)$.
At any fixed $Y>1$, the factorization \eqref{trace:factorization}
also gives a trace-class operator $S^*\psi(A_X)S$ for every bounded
Borel function $\psi$ on $[0,a]$, with $S=A_X^{1/2}\chi_Y$.
Thus the measure construction requires no degree bound; the bound
$\deg p\le Q$ ensures that the geometric identity holds at $Y_Q$.
The compression of $1$ is $\chi_{Y_Q}$, an infinite-rank projection,
so it is not trace class.  Formula
\eqref{trace:P-definition} is independent of the chosen support radius,
provided the corresponding interval contains $\supp H_p$; at a fixed $Y_Q$ it is one linear
functional for all degrees at most $Q$.

\subsection{Deterministic norms and the constant eigenfunction}
\label{sec:tracenorms}

We first bound the identity and parabolic terms using only the
multiplier on the real spectral axis, then subtract the constant
mode. The parabolic estimate uses the following absolute Abel norm.
For an even $H\in C_c^\infty(\R)$ supported
in $[-T,T]$, with $T>0$, put
\begin{equation}\label{trace:abel-norm-definition}
 B(H)=2\int_0^T|K_H(\rho)|\cosh(\rho/2)\,d\rho,
 \qquad \log_+v=\log\max\{1,v\}\quad(v\ge0).
\end{equation}
The value of $B(H)$ is independent of the chosen support radius.

\begin{proposition}[Deterministic trace bounds]\label{trace:parabolic}
There are finite nonnegative constants $I_f,J_f,b_f$, depending only on the fixed
multiplier, such that for every integer $Q\ge1$ and real polynomial
$p=x\widetilde p$ of degree
at most $Q$,
\begin{equation}\label{trace:deterministic-bounds}
 |\mathcal I(p)|\le I_fN_{\R}(p),\qquad
 |\mathcal P_{Y_Q}(p)|\le J_fQ^2N_{\R}(p),\qquad
 |\mathcal P_{Y_Q}(x)|\le b_fQ.
\end{equation}
One may take
\begin{align}
 I_f&=\frac1{4\pi}\int_{\R}r\tanh(\pi r)f(r)\,dr,
 &C_{1,f}&=\frac1{2\pi}\int_{\R}|r|f(r)\,dr,
       \label{trace:constants-one}\\
 J_f&=LC_{1,f}(1+\log3+3L/2),
 &b_f&=B(h)\bigl[1+L+\log_+(6\sinh(L/2))\bigr].
       \label{trace:constants-two}
\end{align}
\end{proposition}

\begin{proof}
The absolute Abel norm obeys
\begin{equation}\label{trace:absolute-abel}
 B(H)\le\int_0^T|H'(u)|\,du.
\end{equation}
Indeed, substitute \eqref{trace:abel-inverse}, apply Tonelli, and use
\[
 \int_0^u\frac{\cosh(\rho/2)}
                   {\sqrt{\cosh u-\cosh\rho}}\,d\rho
      =\frac\pi{\sqrt2}\qquad(u>0).
\]
The substitution $v=\sinh(\rho/2)$ gives this elementary integral.
Taking absolute values in \eqref{trace:parabolic-unfolding} therefore
gives, without any positivity assumption on $K_H$,
\begin{equation}\label{trace:parabolic-absolute}
 |\mathcal P_Y(H)|\le B(H)
       \mathsf H_{\lfloor2Y\sinh(T/2)\rfloor}
 \le B(H)\bigl[1+\log_+(2Y\sinh(T/2))\bigr].
\end{equation}
Since $|p(f(r))|\le N_{\R}(p)f(r)$ on the real axis,
\eqref{trace:I-definition} proves the identity bound, and differentiated
Fourier inversion gives
$\|H_p'\|_\infty\le C_{1,f}N_{\R}(p)$.
At $T=LQ$, \eqref{trace:absolute-abel} yields
$B(H_p)\le LQC_{1,f}N_{\R}(p)$, while
\[
 2Y_Q\sinh(LQ/2)\le3e^{3LQ/2},\qquad
 \mathsf H_{\lfloor2Y_Q\sinh(LQ/2)\rfloor}
       \le1+\log3+3LQ/2.
\]
This proves the second bound in \eqref{trace:deterministic-bounds}.
For the base test $p=x$, the test function $H_p=h$ has support
in $[-L,L]$, independently of $Q$, so only the cutoff
$Y_Q$ contributes growth in \eqref{trace:parabolic-absolute}:
\[
 |\mathcal P_{Y_Q}(x)|
 \le B(h)\bigl[1+LQ+\log_+(6\sinh(L/2))\bigr]\le b_fQ.
\]
\end{proof}

Let $P_0$ denote the orthogonal projection onto the constant
eigenfunction.  The normalized eigenfunction is
$u_0=\mathcal A_{g,n}^{-1/2}$, so its truncated squared mass equals
$1-n/(\mathcal A_{g,n}Y_Q)$.  In particular
\begin{equation}\label{trace:remove-constant}
 \tau_{X,Q}^{\circ}(p)
 :=\Tr\bigl(\chi_{Y_Q}p(A_X)(1-P_0)\chi_{Y_Q}\bigr)
 =\tau_{X,Q}(p)
       -\left(1-\frac{n}{\mathcal A_{g,n}Y_Q}\right)p(a).
\end{equation}
Subtracting the displayed rank-one compression preserves trace class.
Since $P_0$ is a spectral projection, it commutes with $p(A_X)$.
The compressed operator is therefore positive whenever $p$ is
nonnegative on the remaining spectral image, and
$0\le\tau_{X,Q}^{\circ}(x)\le\tau_{X,Q}(x)$.

Combining \eqref{trace:exact-identity}, \eqref{trace:remove-constant}
and Proposition~\ref{trace:leading} gives the exact decomposition
\begin{equation}\label{trace:corrected-identity}
 \begin{split}
 \tau_{X,Q}^{\circ}(p)
 ={}&\left(\mathcal A_{g,n}-\frac n{Y_Q}\right)\mathcal I(p)
       +n\mathcal P_{Y_Q}(p)+\mu_0(p)\\
    &+\bigl(\mathcal O_X(p)-\omega_1(p)\bigr)
       +\frac{np(a)}{\mathcal A_{g,n}Y_Q}.
 \end{split}
\end{equation}
The remaining constant-mode term, $np(a)/(\mathcal A_{g,n}Y_Q)$,
comes from the mass removed by truncation.
The signed measure $\mu_0$ is supported on $[0,1]$, and the identity
and parabolic terms are controlled by $N_{\R}(p)$ above.
Section~\ref{sec:polynomial} estimates the remaining expectation
$\E_{g,n}[\mathcal O_X(p)-\omega_1(p)]$ through the orbital expansions.

\section{Polynomial approximation and spectral detection}
\label{sec:polynomial}

We prove the spectral count in Theorem~\ref{addspec:spectral} by
applying a positive polynomial detector to the truncated trace.
The detector is small on the real spectral image and close to one
on the interval corresponding to low eigenvalues. A uniform bound
on the mass lost in the cusps then converts its trace into an
eigenvalue count.

Proposition~\ref{poly:probability} illustrates the interpolation
method with a preliminary count. The final estimate uses the
coefficient bounds of Section~\ref{sec:sharppm}, the additive
remainder of Section~\ref{sec:additivepolynomial}, and the
degree-independent first-correction bound obtained from positivity
in Section~\ref{sec:exceptionalmeasure}.

Throughout this section, $g\ge3$, $n\ge1$ and the degree bound $Q\ge1$
are integers. Polynomial tests are real and have zero constant term.
The integer $c=c_f\ge2$ is fixed by Proposition~\ref{orb:cgm}; the
constants introduced for the additive transfer below are distinct.
For the preliminary approximation set
\begin{equation}\label{poly:parameters}
 \alpha=2c+2,\qquad \beta=2c+6,\qquad
 G_{n,Q}=\left\lceil4[c(2cQ+2)(n+1)]^{c+1}\right\rceil,
 \qquad Y_Q=3e^{LQ}.
\end{equation}
For $p=x\widetilde p$ we use the norms
$N_*(p)=\|\widetilde p\|_{[0,a]}$ and
$N_{\mathbb R}(p)=\|\widetilde p\|_{[0,1]}$ from the preceding section.
The interpolation method has a compact-surface precedent in
\cite[Lemma~3.1 and equations~(3.1)--(3.9)]{HMT25}; the argument below retains the cusp terms and gives an explicit genus threshold.

\subsection{A first moment and reciprocal-genus interpolation}
\label{sec:polynomialtools}

We begin with a first moment that will control all polynomial traces
through positivity. The exact identity of Theorem~\ref{trace:exact}, with
$\tau_{X,Q}(p)=\Tr(\chi_{Y_Q}p(A_X)\chi_{Y_Q})$, reads
\begin{equation}\label{poly:trace-interface}
 \tau_{X,Q}(p)
 = (\mathcal A_{g,n}-n/Y_Q)\mathcal I(p)
   +n\mathcal P_{Y_Q}(p)+\mathcal O_X(p),
 \qquad |\tau_{X,Q}(p)|\le N_*(p)\tau_{X,Q}(x).
\end{equation}
Proposition~\ref{trace:parabolic} supplies fixed finite constants
$I_f,J_f,b_f\ge0$ such that
\begin{equation}\label{poly:deterministic-bounds}
 |\mathcal I(p)|\le I_fN_{\mathbb R}(p),\qquad
 |\mathcal P_{Y_Q}(p)|\le J_fQ^2N_{\mathbb R}(p),\qquad
 |\mathcal P_{Y_Q}(x)|\le b_fQ.
\end{equation}
In particular, these bounds use only the multiplier on the real spectral axis; they assign no sign to the parabolic term.

Apply Proposition~\ref{orb:cgm} with $(q,t)=(2,1)$. Since
Proposition~\ref{trace:leading} gives $|\omega_1(x)|\le6a$, we have
\begin{equation}\label{poly:base-moment}
 \big|\E_{g,n}\mathcal O_X(x)-\omega_1(x)\big|
 \le\frac{[2c(n+1)]^{2c}}g\le g,
 \qquad g\ge[2c(n+1)]^c.
\end{equation}
The base orbital statistic is nonnegative because $h\ge0$, and its expectation is finite. Consequently the positive trace in
\eqref{poly:trace-interface} is integrable. Since
$0\le\mathcal A_{g,n}-n/Y_Q\le4\pi g+2\pi n$, equations
\eqref{poly:deterministic-bounds}--\eqref{poly:base-moment} imply
\begin{equation}\label{poly:first-moment}
 \E_{g,n}\tau_{X,Q}(x)\le D_f(g+nQ),
 \qquad D_f=1+6a+6\pi I_f+b_f.
\end{equation}
This estimate holds simultaneously for every $Q\ge1$ at the explicit genus threshold in \eqref{poly:base-moment}. Trace domination and the exact identity then give absolute integrability and
\begin{equation}\label{poly:apriori}
 \E_{g,n}|\mathcal O_X(p)|
 \le N_*(p)\{D_f(g+nQ)+(4\pi g+2\pi n)I_f+nJ_fQ^2\}.
\end{equation}
Thus, if also $g\ge nQ^2$, then
\begin{equation}\label{poly:normalized-apriori}
 \big|g^{-1}\E_{g,n}\mathcal O_X(p)\big|
 \le A_fN_*(p),\qquad A_f=2D_f+6\pi I_f+J_f.
\end{equation}
This a priori bound uses only the leading orbital coefficient.
The next two lemmas turn it into a polynomial approximation estimate:
the first controls monomial coefficients, and the second controls a
polynomial from its values at reciprocal integers.

\begin{lemma}[Coefficient bound]\label{poly:coefficients}
If $p(x)=\sum_{t=1}^Q a_tx^t$, then
$\sum_{t=1}^Q|a_t|\le7^QN_*(p)$.
\end{lemma}
\begin{proof}
Expand
$\widetilde p(x)=\sum_{k=0}^{Q-1}\gamma_kT_k(2x/a-1)$ in shifted Chebyshev polynomials. The cosine coefficient formulas give
$|\gamma_0|\le N_*(p)$ and $|\gamma_k|\le2N_*(p)$ for $k\ge1$.
For the sum $\|\cdot\|_{\ell^1}$ of absolute monomial coefficients,
one has $\|2x/a-1\|_{\ell^1}\le3$ because $a>1$.
The recurrence for $T_k$ therefore gives
\[
 \|T_{k+1}(2x/a-1)\|_{\ell^1}
 \le6\|T_k(2x/a-1)\|_{\ell^1}
    +\|T_{k-1}(2x/a-1)\|_{\ell^1}\le7^{k+1},
\]
by induction from degrees zero and one. It follows that
\[
 \|\widetilde p\|_{\ell^1}
 \le\left(1+2\sum_{k=1}^{Q-1}7^k\right)N_*(p)
 =\frac{7^Q-4}{3}N_*(p)\le7^QN_*(p).
\]
Multiplication by $x$ leaves this coefficient sum unchanged.
\end{proof}

\begin{lemma}[Reciprocal-integer grid]\label{poly:grid}
Let $d_0\ge1$ and $G\ge4d_0^2$ be integers, let $M\ge0$, and
let $R$ be a real polynomial of degree at most $d_0$. If $|R(0)|\le M$ and
$|R(1/k)|\le M$ for every integer $k\ge G$, then
\begin{equation}\label{poly:grid-bounds}
 \|R\|_{[0,1/G]}\le2M,\qquad
 \|R''\|_{[0,1/G]}\le8d_0^4G^2M.
\end{equation}
\end{lemma}
\begin{proof}
Write $S=\|R\|_{[0,1/G]}$. Consecutive nodes have spacing
$1/[k(k+1)]\le G^{-2}$, so every point of the interval lies within
$G^{-2}$ of a node or their limiting point zero. The classical Markov inequality, rescaled from $[-1,1]$, is
$\|R'\|_{[0,1/G]}\le2d_0^2GS$.
Apply the mean value theorem at a point attaining $S$ to obtain
$S\le M+(2d_0^2/G)S\le M+S/2$.
Hence $S\le2M$. Applying Markov once more gives
$\|R''\|\le4d_0^2(d_0-1)^2G^2S\le8d_0^4G^2M$.
For degree at most one the second derivative vanishes.
\end{proof}

\begin{proposition}[Polynomial approximation]\label{poly:pm}
There is $C_{\rm PM}<\infty$, depending only on $f$ and $c$, such that for
all integers $n,Q\ge1$, $g\ge G_{n,Q}$, and real polynomials
$p=x\widetilde p$ of degree at most $Q$,
\begin{equation}\label{poly:pm-estimate}
 \big|\E_{g,n}\mathcal O_X(p)-\omega_1(p)\big|
 \le \frac{C_{\rm PM}(n+1)^\alpha Q^\beta}{g}N_*(p).
\end{equation}
Moreover,
\begin{equation}\label{poly:threshold-upper}
 G_{n,Q}\le5[3c^2Q(n+1)]^{c+1}.
\end{equation}
\end{proposition}
\begin{proof}
We keep the cusp count and test polynomial fixed while varying the
genus. This bounds a single polynomial in reciprocal genus; Taylor's
theorem then isolates its linear coefficient.
Fix $n,Q,p$ and put $q=2cQ+2$, $u=cq(n+1)$,
$G=G_{n,Q}$, and $d_0=q-1$. The coefficient functionals in
Proposition~\ref{orb:cgm} are independent of genus and truncation order, so
\[
 R_p(z)=\sum_{j=1}^{q-1}\omega_{j,n}(p)z^j
\]
is one fixed polynomial while the integer genus varies. It has
$R_p(0)=0$ and $R_p'(0)=\omega_1(p)$, by Proposition~\ref{trace:leading}.
Linearity, $Q<q$, and Lemma~\ref{poly:coefficients} yield
\begin{equation}\label{poly:expanded-polynomial}
 \big|g^{-1}\E_{g,n}\mathcal O_X(p)-R_p(1/g)\big|
 \le7^Qu^{cq}g^{-q}N_*(p),\qquad g\ge u^c.
\end{equation}
The definition of $G$ gives
\begin{equation}\label{poly:threshold-comparisons}
 G\ge u^c\ge[2c(n+1)]^c,\qquad
 G\ge nQ^2,\qquad G\ge4d_0^2,\qquad
 G^{q-2}\ge7^Qu^{cq}.
\end{equation}
For the first three comparisons use $u\ge q$, $u\ge Q(n+1)$,
and $c+1\ge3$. For the last one, the exact exponent calculation is
\[
 (c+1)(q-2)-cq=2c(Q-1),\qquad
 \frac{G^{q-2}}{u^{cq}}
 \ge4^{2cQ}u^{2c(Q-1)}\ge7^Q.
\]
Thus for every integer $g\ge G$, equations
\eqref{poly:normalized-apriori} and \eqref{poly:expanded-polynomial} give
\[
 \big|g^{-1}\E_{g,n}\mathcal O_X(p)-R_p(1/g)\big|
 \le g^{-2}N_*(p),\qquad
 |R_p(1/g)|\le(A_f+1)N_*(p).
\]
Lemma~\ref{poly:grid}, followed by Taylor's theorem at zero, now implies
\begin{equation}\label{poly:taylor}
 \big|\E_{g,n}\mathcal O_X(p)-\omega_1(p)\big|
 \le\frac{1+4(A_f+1)d_0^4G^2}{g}N_*(p).
\end{equation}
Since $q\le3cQ$, $u\le3c^2Q(n+1)$, and
$G\le5u^{c+1}$, equation~\eqref{poly:threshold-upper} follows, and
\[
 d_0^4G^2\le25(3c)^4(3c^2)^{2c+2}
 Q^{2c+6}(n+1)^{2c+2}.
\]
This proves \eqref{poly:pm-estimate}; for example one may take
$C_{\rm PM}=1+100(A_f+1)(3c)^4(3c^2)^{2c+2}$.
\end{proof}

\subsection{A positive detector}\label{sec:detector}

Fix $0<\varepsilon<1/4$. We seek a polynomial that is exponentially
small on $[0,1]$ and close to one on $[b_\varepsilon,a]$, while its
quotient by $x$ remains uniformly bounded on $[0,a]$.
The constants below specify the transition and the degree thresholds:
\begin{equation}\label{poly:detector-constants}
\begin{gathered}
 b_\varepsilon=f(i\sqrt\varepsilon),\qquad
 \theta=\frac{1+b_\varepsilon}{2a},\qquad
 d=\frac{b_\varepsilon-1}{2a},\qquad
 \kappa=\frac{d^2}{4},\qquad B_\varepsilon=\frac1{a\theta},\\
 Q_{\rm w}=\max\left\{2,\ 1+\left\lceil\frac{2(1-\theta)}d\right\rceil\right\},\qquad
 Q_\varepsilon=\max\left\{Q_{\rm w},
 \left\lceil\frac{\log2}{8\kappa}\right\rceil,
 \left\lceil\frac{\log2}{2\sqrt\varepsilon L}\right\rceil\right\}.
\end{gathered}
\end{equation}
Indeed,
$f(is)=\int_{\R}h(u)\cosh(su)\,du$ is strictly increasing for $s>0$:
its derivative is $\int u h(u)\sinh(su)\,du>0$, since a smooth probability density is not supported at zero. Therefore
$1<b_\varepsilon<a$, $\theta-d=1/a$, and
$\theta+d=b_\varepsilon/a<1$.

\begin{lemma}[Binomial detector]\label{poly:detector}
For every integer $Q\ge1$, put
\begin{equation}\label{poly:binomial-polynomial}
 k_Q=\lceil\theta Q\rceil,\qquad
 p_Q(x)=\sum_{j=k_Q}^Q\binom Qj(x/a)^j(1-x/a)^{Q-j}.
\end{equation}
Then $p_Q(0)=0$, $p_Q(a)=1$, $0\le p_Q\le1$ on $[0,a]$, and
\begin{align}
 N_*(p_Q)&\le B_\varepsilon,\label{poly:full-detector-norm}\\
 N_{\mathbb R}(p_Q)&\le B_\varepsilon e^{-\kappa Q}
       &&(Q\ge Q_{\rm w}),\label{poly:real-detector-norm}\\
 \inf_{[b_\varepsilon,a]}p_Q&\ge1-e^{-8\kappa Q}.
       \label{poly:detector-lower}
\end{align}
\end{lemma}
\begin{proof}
For $t=x/a\in[0,1]$, $p_Q(x)$ is
$\PP\{\operatorname{Bin}(Q,t)\ge k_Q\}$, giving the endpoint values and positivity. We first establish the concentration bound used below.
For a Bernoulli variable of parameter $t$, let
$F_t(s)=\log(1-t+te^s)-ts$.
Then $F_t(0)=F_t'(0)=0$ and
\[
 F_t''(s)=\frac{t(1-t)e^s}{(1-t+te^s)^2}\le\frac14,
 \qquad F_t(s)\le s^2/8.
\]
Independence and the exponential Markov inequality at $s=\pm4\eta$ give,
for $B\sim\operatorname{Bin}(N,t)$, $N\ge1$, and $\eta>0$,
\begin{equation}\label{poly:concentration}
 \PP(B/N-t\ge\eta)\le e^{-2N\eta^2},\qquad
 \PP(B/N-t\le-\eta)\le e^{-2N\eta^2}.
\end{equation}
If $t\ge\theta+d$, failure of the detector means
$B\le k_Q-1<\theta Q$; hence \eqref{poly:detector-lower} follows from
\eqref{poly:concentration} with $\eta=d$.

To control the quotient at zero, the identity
$j\binom Qj=Q\binom{Q-1}{j-1}$ gives
\begin{equation}\label{poly:size-bias}
 \PP\{\operatorname{Bin}(Q,t)\ge k_Q\}
 \le\frac{Qt}{k_Q}
 \PP\{\operatorname{Bin}(Q-1,t)\ge k_Q-1\}.
\end{equation}
Here the binomial variable with zero trials is identically zero.
Dividing by $x=at>0$ and using $k_Q\ge\theta Q$ proves
\eqref{poly:full-detector-norm}, including $x=0$ by polynomial continuity.
If $Q\ge Q_{\rm w}$ and $0<x\le1$, then $t\le\theta-d$ and
\[
 \frac{k_Q-1}{Q-1}
 \ge\theta-\frac{1-\theta}{Q-1}
 \ge\theta-d/2\ge t+d/2.
\]
Applying \eqref{poly:concentration} to the last probability in
\eqref{poly:size-bias} gives
$N_{\mathbb R}(p_Q)\le B_\varepsilon e^{-(Q-1)d^2/2}
\le B_\varepsilon e^{-\kappa Q}$, since $Q\ge2$.
The displayed threshold also implies $k_Q\ge2$, so the quotient vanishes at zero in this range.
\end{proof}

The exponentially small norm $N_{\mathbb R}(p_Q)$ controls the
identity and parabolic terms through \eqref{poly:deterministic-bounds}.
The uniformly bounded norm $N_*(p_Q)$ controls the orbital error.

\subsection{Retained mass after cusp truncation}

Hide \cite[Lemma~2.7]{Hide23}, following Gamburd, gives a uniform
lower bound on the mass retained at a fixed cusp height. We prove
the following explicit estimate for the truncation $Y_Q$.

\begin{lemma}[Uniform retained mass]\label{poly:mass}
If $\psi$ is a normalized $L^2$ eigenfunction with
$0<\lambda\le1/4-\varepsilon$, then
\begin{equation}\label{poly:mass-bound}
 \|\chi_{Y_Q}\psi\|_2^2\ge1-e^{-2\sqrt\varepsilon LQ}.
\end{equation}
This estimate includes residual eigenfunctions and is uniform in the surface and its cusp count.
\end{lemma}
\begin{proof}
Write $\lambda=1/4-s^2$, so $\sqrt\varepsilon\le s<1/2$.
In a width-one cusp, Fourier expansion and the eigenfunction equation give
\begin{equation}\label{poly:cusp-fourier}
 \psi(x,y)=a_0y^{1/2-s}
  +\sum_{m\ne0}a_m\sqrt y\,K_s(2\pi|m|y)e^{2\pi imx}.
\end{equation}
To see the allowed terms, the $m$th mode satisfies
$-y^2\psi_m''+(2\pi m)^2y^2\psi_m=(1/4-s^2)\psi_m$.
For $m=0$, square integrability excludes $y^{1/2+s}$.
For $m\ne0$, writing $\psi_m(y)=\sqrt y\,v(2\pi|m|y)$ gives the modified Bessel equation. Its positive-axis asymptotics exclude the growing $I_s$ solution and leave $K_s$; see \cite[equations~10.40.1--10.40.2]{DLMF}. The zero Fourier coefficient $a_0$ need not vanish.

The positive integral representation and order symmetry
\cite[equations~10.32.10 and~10.27.3]{DLMF}, applied at order $-s$, give
\begin{equation}\label{poly:bessel-integral}
 z^sK_s(z)=2^{s-1}\int_0^\infty t^{s-1}
       \exp\left(-t-\frac{z^2}{4t}\right)\,dt,
 \qquad z>0.
\end{equation}
The integrand is positive and decreases with $z$. Hence
$0<K_s(rz)\le r^{-s}K_s(z)$ for $r\ge1$.
If $M_b(Y)$ is the squared mass in cusp $b$ above height $Y$, Parseval and Tonelli yield
\begin{equation}\label{poly:cusp-mass}
 M_b(Y)=\frac{|a_0|^2}{2s}Y^{-2s}
  +\sum_{m\ne0}|a_m|^2
       \int_Y^\infty K_s(2\pi|m|y)^2\,\frac{dy}{y}.
\end{equation}
The zero mode scales exactly by $r^{-2s}$ when $Y$ is replaced by $rY$.
In the other integrals substitute $y=rv$ and use the same Bessel bound. It follows that
$M_b(3r)\le r^{-2s}M_b(3)$.
The height-three cusp neighborhoods are disjoint, so
$\sum_bM_b(3)\le1$ and therefore $\sum_bM_b(3r)\le r^{-2s}$.
Taking $r=e^{LQ}$ proves \eqref{poly:mass-bound} without a factor depending on the number of cusps.
\end{proof}

\subsection{Exact cancellation and a preliminary spectral count}

Let $P_0$ denote orthogonal projection onto constants. The normalized constant eigenfunction has retained squared mass
$1-n/(\mathcal A_{g,n}Y_Q)$. Since $p_Q(a)=1$, define
\begin{equation}\label{poly:nonconstant-trace}
\begin{split}
 D_{X,Q}
 &=\tau_{X,Q}(p_Q)-\left(1-\frac{n}{\mathcal A_{g,n}Y_Q}\right)\\
 &=\Tr\bigl(\chi_{Y_Q}p_Q(A_X)(I-P_0)\chi_{Y_Q}\bigr)\ge0.
\end{split}
\end{equation}
The compressed operator differs from
$\chi_{Y_Q}p_Q(A_X)\chi_{Y_Q}$ by a rank-one operator, so it is
trace class by Theorem~\ref{trace:exact}.
It is positive because $P_0$ commutes with $p_Q(A_X)$ and
$p_Q\ge0$ on $[0,a]$; no commutation of $P_0$ with the cutoff is needed.
The positive spectral expansion contains all nonconstant discrete eigenfunctions and the continuous contribution. Every eigenvalue counted by $N_\varepsilon(X)$ has multiplier in $[b_\varepsilon,a]$.
Lemmas~\ref{poly:detector} and~\ref{poly:mass} therefore give
\begin{equation}\label{poly:counting}
 D_{X,Q}\ge
 (1-e^{-8\kappa Q})(1-e^{-2\sqrt\varepsilon LQ})N_\varepsilon(X)
 \ge\tfrac14N_\varepsilon(X),\qquad Q\ge Q_\varepsilon.
\end{equation}
This follows first for any finite collection of orthonormal low eigenfunctions. Finiteness of the positive trace then also proves finiteness of the count, with multiplicities.

\begin{proposition}[Finite-parameter probability bound]\label{poly:probability}
For every integer $n\ge1$, $Q\ge Q_\varepsilon$, and
$g\ge G_{n,Q}$,
\begin{equation}\label{poly:probability-bound}
 \PP_{g,n}(N_\varepsilon\ge1)\le\E_{g,n}N_\varepsilon
 \le C_{f,\varepsilon,c}\left\{
 (g+nQ^2)e^{-\kappa Q}
 +\frac{(n+1)^\alpha Q^\beta}{g}
 +\frac{n}{\mathcal A_{g,n}Y_Q}\right\}.
\end{equation}
The constant is independent of $g,n,Q$.
\end{proposition}
\begin{proof}
Proposition~\ref{trace:leading} identifies
$\omega_1(p)=p(a)+\mu_0(p)$, with
$\supp\mu_0\subset[0,1]$ and $\|\mu_0\|_{\TV}<5$.
Combining this identity with \eqref{poly:trace-interface},
Proposition~\ref{poly:pm}, and \eqref{poly:nonconstant-trace} gives the exact cancellation
\begin{equation}\label{poly:cancelled-expectation}
\begin{split}
 \E_{g,n}D_{X,Q}
 ={}&(\mathcal A_{g,n}-n/Y_Q)\mathcal I(p_Q)
    +n\mathcal P_{Y_Q}(p_Q)+\mu_0(p_Q)\\
    &+\frac{n}{\mathcal A_{g,n}Y_Q}+e_{g,n,Q},\qquad
 |e_{g,n,Q}|\le
 \frac{C_{\rm PM}(n+1)^\alpha Q^\beta}{g}N_*(p_Q).
\end{split}
\end{equation}
All expectations are absolutely convergent by
\eqref{poly:first-moment}--\eqref{poly:apriori}.
Since $|\mu_0(p)|\le5N_{\mathbb R}(p)$, the detector and deterministic bounds imply
\begin{equation}\label{poly:explicit-expectation}
\begin{split}
 0\le\E_{g,n}D_{X,Q}\le{}&
 B_\varepsilon(\mathcal A_{g,n}I_f+nJ_fQ^2+5)e^{-\kappa Q}\\
 &+B_\varepsilon C_{\rm PM}\frac{(n+1)^\alpha Q^\beta}{g}
 +\frac{n}{\mathcal A_{g,n}Y_Q}.
\end{split}
\end{equation}
Using $\mathcal A_{g,n}\le4\pi g+2\pi n$, $Q\ge1$, and
\eqref{poly:counting} proves the expected-count estimate.
The probability estimate follows from
$\ind_{\{N_\varepsilon\ge1\}}\le N_\varepsilon$.
\end{proof}

The unit coefficient of the atom at $a$ leaves exactly the displayed
truncation correction when the constant mode is removed. No interval
of nonconstant eigenvalues near zero has to be excluded.

\subsection{A common logarithmic degree}\label{sec:logdegree}

The detector degree must grow to suppress the identity and parabolic
terms. The following logarithmic choice will suffice in the final
spectral estimate:
\begin{equation}\label{poly:log-degree}
 K=\max\{1,4/\kappa,2/L\},\qquad Q_g=\lceil K\log g\rceil.
\end{equation}
For all sufficiently large $g$, independently of $n$,
$Q_g\ge Q_\varepsilon$, $\log g\ge1$, and $Q_g\le2K\log g$.
Whenever $1\le n\le g$, the following bounds hold uniformly:
\begin{align}
 (g+nQ_g^2)e^{-\kappa Q_g}
 &\le g^{-3}\bigl(1+4K^2(\log g)^2\bigr)
   =O_{f,\varepsilon}(g^{-2}),\label{poly:exp-decay}\\
 \frac{n}{\mathcal A_{g,n}Y_{Q_g}}
 &\le\frac1{6\pi}e^{-LQ_g}\le\frac1{6\pi}g^{-2},
 \label{poly:area-decay}
\end{align}
where the second line uses $\mathcal A_{g,n}\ge2\pi n$.

\subsection{Uniform norms of the cusp-polynomial coefficients}
\label{sec:sharppm}

We now use reciprocal-genus interpolation to bound each canonical
coefficient at finitely many cusp counts. Cusp polynomiality then
extends the bound to all $n$, with an explicit power of $n+1$.
Throughout the remainder of this section write $s=n+1$.

\begin{lemma}[Uniform norms of the cusp-polynomial coefficients]
\label{sharppm:coefficients}
Fix the integer $D_0=2c+4$. There is a constant $C_0\ge2$,
depending only on $f$, such that
for all integers $j\ge2$, $Q,n\ge1$, and all real polynomials
$p=x\widetilde p$ of degree at most $Q$,
\begin{equation}\label{sharppm:coefficientbound}
 |\omega_{j,n}(p)|
 \le s^j[C_0(Q+j)]^{D_0j}N_*(p).
\end{equation}
No relation between $j$ and $Q$ is assumed.
\end{lemma}
\begin{proof}
We first bound the coefficient at the fixed cusp counts
$1,\ldots,j+1$ by varying an auxiliary integer genus. Polynomiality
in the cusp count then gives the estimate at arbitrary $n$; the
auxiliary genus threshold is not imposed on the surfaces used later.
Put $A=A_f+1$, with $A_f$ from \eqref{poly:normalized-apriori}, and
\begin{equation}\label{sharppm:gridchoices}
 M=Q+j,\quad q=2cM+2,\quad d_0=q-1,\quad
 u=cq(j+2),\quad G=\lceil4u^{c+1}\rceil.
\end{equation}
At each fixed cusp count $r\in\{1,\ldots,j+1\}$, use the
polynomial $R_{r,p}(z)=\sum_{k=1}^{q-1}\omega_{k,r}(p)z^k$.
Proposition~\ref{orb:cgm}, applied only at these cusp counts, and
Lemma~\ref{poly:coefficients} give
\begin{equation}\label{sharppm:nodeexpansion}
 |g^{-1}\E_{g,r}\mathcal O_X(p)-R_{r,p}(1/g)|
 \le7^Qu^{cq}g^{-q}N_*(p),\qquad g\ge u^c.
\end{equation}
Here $q>Q,j$ and the coefficients are independent of the varying
integer genus.  The common grid threshold satisfies
\begin{equation}\label{sharppm:gridcomparisons}
 G\ge u^c\ge[2c(r+1)]^c,\qquad
 G\ge rQ^2,\qquad G\ge4d_0^2,\qquad
 G^{q-2}\ge7^Qu^{cq}.
\end{equation}
Indeed, $u\ge q$, $u\ge Q(r+1)$ and $c+1\ge3$ prove the first
three bounds. For the last one,
\begin{equation}\label{sharppm:absorbmonomials}
 (c+1)(q-2)-cq=2c(M-1),\qquad
 G^{q-2}/u^{cq}\ge4^{2cM}u^{2c(M-1)}\ge7^Q.
\end{equation}
The first two comparisons allow \eqref{poly:normalized-apriori}
at every node, while the last bounds the expansion error by
$g^{-2}N_*(p)$. Hence
$|R_{r,p}(1/g)|\le AN_*(p)$ for all integers $g\ge G$.
Also $R_{r,p}(0)=0$.
Lemma~\ref{poly:grid} gives
$\|R_{r,p}\|_{[0,1/G]}\le2AN_*(p)$.
Iterating the rescaled Markov inequality $j$ times gives
\begin{equation}\label{sharppm:nodecoefficient}
 |\omega_{j,r}(p)|=
 \frac{|R_{r,p}^{(j)}(0)|}{j!}
 \le\frac{2A}{j!}(2Gd_0^2)^jN_*(p)=:B_j(p).
\end{equation}

Recall that $\omega_{j,n}(p)=\kappa_{j-1,n}(H_p)$.
By Theorem~\ref{poly:orbital}, $W(n)=\omega_{j,n}(p)$ is therefore
a polynomial of degree at most $j$. Its shifted Newton formula is
\begin{equation}\label{sharppm:shiftednewton}
 W(n)=\sum_{r=0}^j\binom{n-1}{r}\Delta^rW(1),\qquad
 \Delta^rW(1)=\sum_{v=0}^r(-1)^{r-v}\binom rvW(1+v).
\end{equation}
For $n\ge1$, all binomial factors in the first sum are
nonnegative and bounded by $s^j$; those with $r>n-1$ vanish.
The second sum has absolute value at most $2^rB_j(p)$.
Consequently
\begin{equation}\label{sharppm:newtonbound}
 |W(n)|\le2^{j+1}s^jB_j(p)
 \le\frac{4A}{j!}s^j(4Gd_0^2)^jN_*(p)
 \le s^j(16AGd_0^2)^jN_*(p).
\end{equation}
Finally $q\le3cM$, $j+2\le2M$, $u\le6c^2M^2$ and
$G\le5u^{c+1}$, so
\begin{equation}\label{sharppm:explicitconstants}
 16AGd_0^2\le720Ac^2(6c^2)^{c+1}M^{2c+4}
 \le[C_0M]^{D_0},
\end{equation}
where $D_0=2c+4$ and one may take
$C_0=\max\{2,[720Ac^2(6c^2)^{c+1}]^{1/D_0}\}$.
This proves \eqref{sharppm:coefficientbound}, including $j>Q$.
We interpolate only the coefficient polynomial in $n$; no finite-genus
remainder is extrapolated to the actual cusp count.
\end{proof}

Fix the constants $C_0,D_0$ from this lemma and put $B=2D_0$.
They depend only on the fixed multiplier $f$.

\subsection{Polynomial transfer at the additive threshold}
\label{sec:additivepolynomial}

The additive genus hypothesis of Theorem~\ref{addgeo:orbital} is
preserved when we pass from monomials to polynomial observables.
We first specialize the orbital expansion, then estimate its terminal
remainder and retained coefficients separately. We keep the first
correction explicit because the next subsection will bound its
detector response independently of degree.

\begin{lemma}[Moments and integrability at the additive threshold]
\label{addspec:moments}
There is a fixed integer $c_3\ge2$, depending only on $f$, such
that, for integers $q>t\ge1$ and
$g\ge c_3s^2+[c_3q]^{c_3}s$,
\begin{equation}\label{addspec:momentbound}
 \left|g^{-1}\E_{g,n}\mathcal O_X(x^t)
       -\sum_{j=1}^{q-1}g^{-j}\omega_{j,n}(x^t)\right|
 \le g^{-q}s^q[c_3q]^{c_3q}.
\end{equation}
Here $\omega_{j,n}$ are the coefficients defined in \eqref{orb:polynomialtest}. Simultaneously for
every $Q\ge1$,
\begin{equation}\label{addspec:firstmoment}
 \E_{g,n}\tau_{X,Q}(x)\le D_f(g+nQ),
 \qquad g\ge c_3s^2+[2c_3]^{c_3}s,
\end{equation}
where $D_f$ is the constant in \eqref{poly:first-moment}.
At this threshold, the traces and orbital statistics associated with
polynomials of degree at most $Q$ and with zero constant term are
absolutely integrable.
\end{lemma}
\begin{proof}
Apply Theorem~\ref{addgeo:orbital} with $R=q-2$, $H=h^{*t}$,
and $T=tL$. Since $\|h\|_1=1$, convolution gives
$\|h^{*t}\|_\infty\le\|h\|_\infty$, and
$\rho=q-1+\lceil tL\rceil\le(1+L)q$.
Writing $C\ge2$ for its absolute constant, one may take
\begin{equation}\label{addspec:momentconstant}
 \begin{gathered}
 U=C(1+L),\qquad H_0=\max\{1,\|h\|_\infty\},\\
 c_3=\left\lceil3+U+U\log U+3L+\log H_0\right\rceil.
 \end{gathered}
\end{equation}
Indeed, after division by $g$, the remaining error factor has
logarithm at most
$q\{U\log q+U\log U+3L+\log H_0\}
\le c_3q\log(c_3q)$.
Also $c_3\ge C$ and $[C\rho]^C\le[c_3q]^{c_3}$, so the
two summands of the threshold in Theorem~\ref{addgeo:orbital} are
dominated separately.
The coefficient of $s^2$ is independent of $q$.

The case $(q,t)=(2,1)$ gives directly
\begin{equation}\label{addspec:basemoment}
 |\E_{g,n}\mathcal O_X(x)-\omega_1(x)|
 \le[2c_3]^{2c_3}s^2/g\le g,
 \qquad g\ge c_3s^2+[2c_3]^{c_3}s.
\end{equation}
For the last inequality use $g\ge[2c_3]^{c_3}s$.
The base orbital statistic is nonnegative and integrable. Thus
the exact trace identity and deterministic bounds give
\[
 \E\tau_{X,Q}(x)
 \le(4\pi g+2\pi n)I_f+nb_fQ+6a+g
 \le D_f(g+nQ).
\]
Trace domination and the exact identity prove the remaining
integrability statements. Thus the first-moment estimate holds under the additive threshold alone.
\end{proof}

\begin{proposition}[Polynomial transfer at the additive threshold]
\label{addspec:pm}
Put $D_3=\max\{D_0,3c_3\}$ and $A_Q=(3C_0Q)^{D_0}$.
There is a fixed $C_{\rm add}\ge2$, depending only on $f$, such
that, for all integers $n,Q\ge1$ and $g$ satisfying
\begin{equation}\label{addspec:pmthreshold}
 g\ge C_{\rm add}s^2+C_{\rm add}Q^{D_3}s,
\end{equation}
every real polynomial $p=x\widetilde p$ of degree at most $Q$ satisfies
\begin{equation}\label{addspec:correctedbound}
 \left|\E_{g,n}\mathcal O_X(p)-\omega_1(p)
              -g^{-1}\omega_{2,n}(p)\right|
 \le\left\{\frac{2s^3A_Q^3}{g^2}
              +\frac{s^2}{g}2^{-Q}\right\}N_*(p).
\end{equation}
At the same threshold one also has
\begin{equation}\label{addspec:pmbound}
 |\E_{g,n}\mathcal O_X(p)-\omega_1(p)|
 \le C_{\rm add}\frac{s^2Q^B}{g}N_*(p).
\end{equation}
\end{proposition}
\begin{proof}
Choose
\begin{equation}\label{addspec:pmconstant}
 C_{\rm add}=\max\{2,c_3,(3c_3)^{c_3},
       2(3C_0)^{D_0},14(3c_3)^{3c_3},3(3C_0)^{2D_0}\}.
\end{equation}
Set $q=Q+2$, so $q\le3Q$ and $q>t$ for every monomial in
$p$. The two summands of \eqref{addspec:pmthreshold} imply
\begin{equation}\label{addspec:comparisons}
 \begin{gathered}
 g\ge c_3s^2+[c_3q]^{c_3}s,\\
 g\ge2sA_Q,\qquad g\ge14(3c_3Q)^{3c_3}s.
 \end{gathered}
\end{equation}
The last two comparisons use only the summand linear in $s$.
They will respectively control the sum of retained coefficients
and absorb the monomial coefficient bound in the terminal remainder.
The same threshold also implies \eqref{addspec:firstmoment}.
Linearity, Lemma~\ref{poly:coefficients}, and Lemma~\ref{addspec:moments}
give the finite expansion
\begin{equation}\label{addspec:expansion}
 \begin{split}
 \E\mathcal O_X(p)
 &=\omega_1(p)+\sum_{j=2}^{Q+1}g^{1-j}\omega_{j,n}(p)
       +\mathcal R,\\
 |\mathcal R|&\le7^Qg^{1-q}s^q[c_3q]^{c_3q}N_*(p).
 \end{split}
\end{equation}
We first estimate the terminal remainder. For $N_*(p)>0$, divide
its bound by $s^2N_*(p)/g$:
\begin{equation}\label{addspec:terminal}
 7^Q[c_3q]^{c_3q}(s/g)^{q-2}
 \le[7(3c_3Q)^{3c_3}s/g]^Q\le2^{-Q}.
\end{equation}
Thus $|\mathcal R|\le2^{-Q}s^2N_*(p)/g$.
If $N_*(p)=0$, then $p=0$.

We next sum the retained corrections of order $j\ge3$.
The coefficient bound in Lemma~\ref{sharppm:coefficients}
holds without a condition on the genus.
Since $Q+j\le3Q$ for retained indices, it gives
\[
 \sum_{j=3}^{Q+1}g^{1-j}|\omega_{j,n}(p)|
 \le\frac{s^3A_Q^3}{g^2}N_*(p)
       \sum_{j=3}^{Q+1}(sA_Q/g)^{j-3}
 \le\frac{2s^3A_Q^3}{g^2}N_*(p).
\]
For $Q=1$ this sum is empty and $q=3$ remains valid.
This proves \eqref{addspec:correctedbound}.
Summing instead from $j=2$ bounds all retained corrections by
$2s^2A_Q^2g^{-1}N_*(p)$; together with
\eqref{addspec:terminal} and \eqref{addspec:pmconstant}, this
proves \eqref{addspec:pmbound}.
\end{proof}

\subsection{Positive exceptional measures and the first correction}
\label{sec:exceptionalmeasure}

The coefficient bound from Section~\ref{sec:sharppm} grows with
detector degree. We now remove this dependence for $\omega_{2,n}(p_Q)$
by using the positive trace. We first extend the weighted coefficient
to a distribution. After passing to its unweighted form, the additive
polynomial estimate and a genus limit at fixed $n$ show that its
restriction to $(1,\infty)$ is a positive measure supported in $(1,a]$.
A cutoff then separates the detector into two parts: one controlled
by the local mass of this measure, and one controlled by exponential
decay below the detector transition.

\begin{lemma}[Weighted distribution extensions]
\label{exc:distributions}
Put $m=B+2=2D_0+2$. For each $n\ge1$, the polynomial functional
\begin{equation}\label{exc:weighteddefinition}
 \mathscr K_n(q)=\omega_{2,n}(xq)
\end{equation}
has a unique extension to a distribution on $\mathbb R$,
supported in $[0,a]$, of order at most $m$.
There is a fixed $K_f$ such that
\begin{equation}\label{exc:distributionnorm}
 |\mathscr K_n(\varphi)|
 \le K_f(n+1)^2\|\varphi\|_{C^m([0,a])}.
\end{equation}
The coefficient $W_2(p)=\Lambda_2(H_p)$ of $\binom n2$ has the unique
weighted extension
\begin{equation}\label{exc:topdistribution}
 \mathscr T_2=\mathscr K_3-2\mathscr K_2+\mathscr K_1,
 \qquad \mathscr T_2(q)=W_2(xq),\qquad
 |\mathscr T_2(\varphi)|\le38K_f\|\varphi\|_{C^m([0,a])}.
\end{equation}
\end{lemma}
\begin{proof}
The construction uses the rapid decay of smooth Chebyshev coefficients
to absorb the degree growth in Lemma~\ref{sharppm:coefficients}
(cf.~\cite[Proposition~3.2]{HMT25}).
With $j=2$, that lemma gives
\[
 |\mathscr K_n(q)|
 \le s^2C_0^B(\deg q+3)^B\|q\|_{[0,a]}.
\]
For a smooth $\varphi$, set
$F(\theta)=\varphi(a(1+\cos\theta)/2)$ and expand
\[
 \varphi(x)=\sum_{k\ge0}\gamma_kT_k(2x/a-1).
\]
The function $F$ is even and periodic. Periodic integration by parts
$m$ times gives
\begin{equation}\label{exc:chebdecay}
 |\gamma_0|\le\|\varphi\|_\infty,
 \qquad |\gamma_k|\le2k^{-m}\|F^{(m)}\|_\infty
 \le2J_{a,m}k^{-m}\|\varphi\|_{C^m}\quad(k\ge1),
\end{equation}
where $\|\varphi\|_{C^m}=\max_{0\le j\le m}\|\varphi^{(j)}\|_\infty$
and one may take $J_{a,m}=(m\max\{1,a/2\})^m$.
Indeed, the chain rule in its partition form has at most $m^m$
terms, and every derivative of $a(1+\cos\theta)/2$ is bounded
by $a/2$.
Apply $\mathscr K_n$ termwise to the series. Since $m=B+2$ and
$(k+3)^B\le4^Bk^B$ for $k\ge1$, the resulting scalar series
converges absolutely and satisfies \eqref{exc:distributionnorm} with
\[
 K_f=C_0^B\{3^B+2J_{a,m}4^B\zeta(2)\}.
\]
For polynomial tests the series is finite, so it agrees with
\eqref{exc:weighteddefinition}. Restriction of smooth tests to $[0,a]$
defines the asserted compactly supported distribution.

For uniqueness, a compactly supported finite-order distribution is
determined by its polynomial values. To verify this directly, let $r$
be its order and approximate the $r$th derivative of a smooth test
uniformly by polynomials on a neighborhood of its support.
Integrating $r$ times and matching the lower derivatives at one endpoint
gives polynomial approximation in $C^r$. A fixed cutoff equal to one
near the support then proves that a distribution vanishing on all
polynomials vanishes on every test.
Finally, the second difference in $n$ of
\eqref{first:quadraticidentity} gives
$W_2=\omega_{2,3}-2\omega_{2,2}+\omega_{2,1}$.
Its norm factor is $4^2+2\cdot3^2+2^2=38$, proving
\eqref{exc:topdistribution}.
\end{proof}

We use the canonical unweighted extensions
\begin{equation}\label{exc:unweighted}
 \mathscr V_n(\psi)=\mathscr K_n\left(
                 \frac{\psi(x)-\psi(0)}x\right).
\end{equation}
The quotient has value $\psi'(0)$ at zero and equals
$\int_0^1\psi'(tx)\,dt$. Thus $\mathscr V_n$ has order at most
$m+1$, support in $[0,a]$, and total mass zero in the distributional
sense. It agrees with $\omega_{2,n}$ on polynomials vanishing at zero.
Only the unweighted extension requires this normalization: without it,
such polynomial data leave a multiple of $\delta_0$ undetermined.

\begin{theorem}[Positive exceptional measures]
\label{exc:positive}
For every $n\ge1$, the restriction of $\mathscr V_n$ to $(1,\infty)$
is a positive Radon measure $\mu_n$ supported in $(1,a]$.
For every fixed $1<b\le a$,
\begin{equation}\label{exc:localmass}
 \mu_n([b,a])\le C_{f,b}(n+1)^2.
\end{equation}
There is also a positive Radon measure $\eta_2$ on $(1,a]$ such that
\begin{equation}\label{exc:topmeasure}
 \mathscr T_2(\psi/x)=\int\psi\,d\eta_2,
 \qquad \eta_2([b,a])\le C_{f,b},
\end{equation}
for tests $\psi\in C_c^\infty((1,\infty))$.
These tests may be nonzero at $a$. In particular the statement includes
possible endpoint mass. Moreover $2n^{-2}\mu_n$ converges vaguely to
$\eta_2$ on $(1,a]$ as $n\to\infty$ through the integers.
\end{theorem}
\begin{proof}
Fix $n$ and let $0\le\psi\in C_c^\infty((1,\infty))$. Put
$\xi(x)=\psi(x)/x$, extended by zero near zero. Choose
\begin{equation}\label{exc:fixednscale}
 0<\sigma<\min\{1/D_3,1/(3D_0)\},
 \qquad Q=\lfloor g^\sigma\rfloor.
\end{equation}
The restriction $\sigma<1/D_3$ ensures that the threshold in
Proposition~\ref{addspec:pm} holds eventually at this fixed $n$,
since $Q^{D_3}/g\to0$. The restriction $\sigma<1/(3D_0)$ will
make the remainder vanish after multiplication by $g$.
Let $r_Q$ be the shifted-Chebyshev truncation of $\xi$ of
degree at most $Q-1$, and set $v_Q(x)=xr_Q(x)$.
Smoothness and \eqref{exc:chebdecay}, used with arbitrarily high order,
give for every fixed $M$
\begin{equation}\label{exc:smoothapprox}
 \delta_Q:=\|r_Q-\xi\|_{[0,a]}\le C_{\psi,M}Q^{-M},
 \qquad r_Q\longrightarrow\xi\text{ in }C^m([0,a]).
\end{equation}
For the second assertion, the iterated Markov inequality bounds the first
$m$ derivatives of the $k$th Chebyshev term by a fixed multiple of
$k^{2m}$, which is absorbed by the rapid coefficient decay.
Consequently $N_*(v_Q)$ is bounded,
$N_{\mathbb R}(v_Q)\le\delta_Q$, and
$\omega_{2,n}(v_Q)\to\mathscr V_n(\psi)$.

Let $P_0$ project onto constants. The smooth compressed trace
\begin{equation}\label{exc:smoothpositive}
 \mathcal D_{X,Q}(\psi)=
 \Tr\bigl(\chi_{Y_Q}\psi(A_X)(I-P_0)\chi_{Y_Q}\bigr)
 \ge0
\end{equation}
is well defined and integrable: on $[0,a]$ one has
$0\le\psi(x)\le C_\psi x$, so the first moment of
Lemma~\ref{addspec:moments} applies.
Spectral domination also gives
\begin{equation}\label{exc:traceapprox}
 |\E\mathcal D_{X,Q}(\psi)-\E\mathcal D_{X,Q}(v_Q)|
 \le\delta_Q\E\tau_{X,Q}(x)
 \le D_f\delta_Q(g+nQ).
\end{equation}
This comparison does not require $v_Q\ge0$ on $[0,a]$.
For each fixed cutoff, polynomial compressions are trace class by
\eqref{trace:factorization} and measurable by integration of their
periodized kernels over the truncated surface. Weighted polynomial
approximation makes the smooth trace their pointwise limit, proving
its measurability. The exact geometric trace identity is used below
only for $v_Q$, whose degree is at most $Q$.

The exact polynomial trace identity and the unit exceptional atom in
$\omega_1$ give
\begin{equation}\label{exc:smoothidentity}
 \begin{split}
 g\E\mathcal D_{X,Q}(v_Q)
 ={}&g(\mathcal A_{g,n}-n/Y_Q)\mathcal I(v_Q)
       +gn\mathcal P_{Y_Q}(v_Q)+g\mu_0(v_Q)\\
    &+\frac{gnv_Q(a)}{\mathcal A_{g,n}Y_Q}
       +\omega_{2,n}(v_Q)+g\mathcal R_{g,n,Q}.
 \end{split}
\end{equation}
By \eqref{poly:deterministic-bounds}, the first three terms have
absolute value at most $C_{f,n}(g^2+gnQ^2+g)\delta_Q$.
The fourth tends to zero because $v_Q(a)$ and
$g/\mathcal A_{g,n}$ are bounded and $Y_Q=3e^{LQ}$.
Equation~\eqref{addspec:correctedbound} gives
\[
 |g\mathcal R_{g,n,Q}|
 \le C_{f,n,\psi}\{Q^{3D_0}/g+2^{-Q}\}\longrightarrow0.
\]
Choose $M$ in \eqref{exc:smoothapprox} large enough that
$g^2\delta_Q\to0$. Since $\sigma<1/D_3\le1/2$ and $n$ is fixed,
all the other approximation errors, including
\eqref{exc:traceapprox} multiplied by $g$, tend to zero as well.
We conclude that
\begin{equation}\label{exc:positivelimit}
 g\E_{g,n}\mathcal D_{X,Q}(\psi)
       \longrightarrow\mathscr V_n(\psi)\ge0.
\end{equation}
The constant mode was removed before taking the limit, so this argument
does not require $\psi(a)=0$.

A positive distribution is a locally finite positive measure. Locally,
choose a nonnegative smooth cutoff $\zeta=1$ on a compact set $K$;
positivity bounds its action on tests supported in $K$ by
$\|\varphi\|_\infty\mathscr V_n(\zeta)$. Uniform approximation
and the representation theorem for positive functionals give the
measure. This proves the assertion for $\mu_n$ on $(1,\infty)$;
its support follows from that of $\mathscr V_n$.
To bound its mass, choose $\zeta_b\in C_c^\infty((1,\infty))$
with $0\le\zeta_b\le1$, zero below $(1+b)/2$ and one on $[b,a]$. Then
\begin{equation}\label{exc:massfromnorm}
 \mu_n([b,a])\le\mathscr V_n(\zeta_b)
 =\mathscr K_n(\zeta_b/x)
 \le K_fs^2\|\zeta_b/x\|_{C^m([0,a])}.
\end{equation}
This proves \eqref{exc:localmass}, uniformly in $n$.

The quadratic identity \eqref{first:quadraticidentity} extends from polynomials to smooth tests through
Lemma~\ref{exc:distributions}. Divide that identity by $\binom n2$
and let $n\to\infty$ through the integers. Since
$\mathscr V_n(\psi)\ge0$, the limit
$\mathscr T_2(\psi/x)$ is nonnegative.
The same positive-functional argument proves \eqref{exc:topmeasure};
its local norm bound follows either from
\eqref{exc:topdistribution} or from \eqref{exc:massfromnorm}
after multiplication by $2/n^2$.
Cusp polynomiality gives convergence on smooth compactly supported tests,
and the uniform local mass bound extends it to continuous compactly supported
tests, proving vague convergence. In \eqref{exc:positivelimit}, the genus limit is taken at fixed $n$,
before the subsequent limit $n\to\infty$.
\end{proof}

\begin{proposition}[Degree-independent detector response]
\label{exc:detector}
For every fixed $0<\varepsilon<1/4$, there is $C_{f,\varepsilon}$
such that, for all $n\ge1$ and $Q\ge Q_\varepsilon$,
\begin{equation}\label{exc:detectorbound}
 |\omega_{2,n}(p_Q)|\le C_{f,\varepsilon}(n+1)^2,
 \qquad |W_2(p_Q)|\le C_{f,\varepsilon}.
\end{equation}
\end{proposition}
\begin{proof}
Recall the detector transition point $a\theta=(1+b_\varepsilon)/2>1$.
Fix $1<r_0<r_1<a\theta$ and a smooth $0\le\chi\le1$ which is zero
below $r_0$ and one above $r_1$.
For a weighted distribution $\mathscr T$ supported in $[0,a]$,
of order at most $m$, whose unweighted exceptional restriction is a locally finite measure
$\nu$, one has
\begin{equation}\label{exc:cutoffsplit}
 \mathscr T(p_Q/x)=\int\chi p_Q\,d\nu
                 +\mathscr T((1-\chi)p_Q/x).
\end{equation}
Since $0\le p_Q\le1$,
the first term has absolute value at most
$|\nu|([r_0,a])$.

Set $\delta=\theta-r_1/a>0$ and
$Q_0=\max\{Q_\varepsilon,2,1+\lceil2(1-\theta)/\delta\rceil\}$.
For $Q\ge Q_0$ and $0<x\le r_1$, the size-bias bound
\eqref{poly:size-bias} and the binomial concentration bound give
\begin{equation}\label{exc:lowdetector}
 \frac{p_Q(x)}x
 \le\frac1{a\theta}
   \PP\{\operatorname{Bin}(Q-1,x/a)\ge k_Q-1\}
 \le B_\varepsilon e^{-\delta^2Q/4}.
\end{equation}
Indeed, $(k_Q-1)/(Q-1)\ge\theta-(1-\theta)/(Q-1)
\ge x/a+\delta/2$. The quotient at zero is defined by continuity.
It is a polynomial of degree at most $Q-1$, so the iterated Markov
inequality on $[0,r_1]$ gives
\[
 \|(p_Q/x)^{(j)}\|_{[0,r_1]}
 \le(2Q^2/r_1)^j B_\varepsilon e^{-\delta^2Q/4},
 \qquad 0\le j\le m.
\]
Derivatives above the actual degree vanish. Leibniz's rule and the
fixed cutoff therefore show that, if $\mathscr T$ has $C^m$ norm at
most $K$, the second term in \eqref{exc:cutoffsplit} is bounded by
$C_{f,\varepsilon}KQ^{2m}e^{-\delta^2Q/4}$.
Only values and derivatives on $[0,r_1]$ are used because $1-\chi$
vanishes above $r_1$.

Apply this with $(\mathscr T,\nu)=(\mathscr K_n,\mu_n)$ and then
$(\mathscr T_2,\eta_2)$, using the distribution and local mass bounds.
The function $Q^{2m}e^{-\delta^2Q/4}$ is bounded, independently of $Q$.
For the finitely many degrees $Q_\varepsilon\le Q<Q_0$, use
\eqref{sharppm:coefficientbound}, \eqref{exc:topdistribution}, and
$N_*(p_Q)\le B_\varepsilon$. This proves both estimates.
No finite total mass near the endpoint $1$ is assumed.
\end{proof}

\subsection{The spectral count in the full sub-square-root range}
\label{sec:additivespectral}

We now combine the additive polynomial remainder with the
degree-independent first-correction bound. At logarithmic detector
degree, the remaining errors are bounded by a constant times
$(n+1)^2/g$ for all sufficiently large $g$ in the additive domain.
We then verify that domain uniformly under a quadratic condition
on the cusp count.

\begin{theorem}[Spectral count throughout the sub-square-root range]
\label{addspec:spectral}
For each fixed $0<\varepsilon<1/4$ there is a constant
$C_{f,\varepsilon}$ such that, for integers $g\ge3$, $n\ge1$,
and $Q\ge Q_\varepsilon$ satisfying \eqref{addspec:pmthreshold},
\begin{equation}\label{addspec:finitecount}
 \begin{split}
 \PP_{g,n}(N_\varepsilon\ge1)\le\E_{g,n}N_\varepsilon
 \le C_{f,\varepsilon}\bigg\{
 &(g+nQ^2)e^{-\kappa Q}+\frac{s^2}{g}
       +\frac{s^3Q^{3D_0}}{g^2}\\
 &+\frac{s^2}{g}2^{-Q}
       +\frac{n}{\mathcal A_{g,n}Y_Q}\bigg\}.
 \end{split}
\end{equation}
After increasing $C_{f,\varepsilon}$, there is a fixed
$g_\varepsilon$, independent of $n$, such that
\begin{equation}\label{addspec:quadraticdomain}
 g\ge g_\varepsilon,\qquad
 g\ge C_{f,\varepsilon}(n+1)^2
\end{equation}
imply
\begin{equation}\label{addspec:cleanbound}
 \PP_{g,n}\{\SpG(X)\le1/4-\varepsilon\}
 \le\E_{g,n}N_\varepsilon
 \le C_{f,\varepsilon}\frac{(n+1)^2}{g}.
\end{equation}
Consequently, for every integer-valued sequence $n(g)\ge1$ with
$n(g)=o(\sqrt g)$, the expected count and this probability tend
to zero. For any envelope $n_{\max}(g)=o(\sqrt g)$,
\begin{equation}\label{addspec:uniformlittleoh}
 \sup_{1\le n\le n_{\max}(g)}\E_{g,n}N_\varepsilon
 \le C_{f,\varepsilon}\frac{(n_{\max}(g)+1)^2}{g}\longrightarrow0
\end{equation}
whenever the range is nonempty. All counts include residual
eigenvalues with multiplicity.
\end{theorem}
\begin{proof}
We first prove the estimate at fixed degree.
Proposition~\ref{exc:detector} bounds the canonical first correction
without a genus hypothesis, so it applies throughout the additive
domain. Thus
$|\omega_{2,n}(p_Q)|\le C_{f,\varepsilon}s^2$.
Every expectation at the additive threshold is justified by
Lemma~\ref{addspec:moments}.

The exact trace identity, the unit exceptional atom of $\omega_1$,
and \eqref{addspec:correctedbound} give
\begin{equation}\label{addspec:cancellation}
 \begin{split}
 \E D_{X,Q}
 ={}&(\mathcal A_{g,n}-n/Y_Q)\mathcal I(p_Q)
       +n\mathcal P_{Y_Q}(p_Q)+\mu_0(p_Q)\\
    &+\frac{n}{\mathcal A_{g,n}Y_Q}
       +g^{-1}\omega_{2,n}(p_Q)+\mathcal R.
 \end{split}
\end{equation}
The detector estimates, deterministic trace bounds, and
$N_\varepsilon\le4D_{X,Q}$ prove
\eqref{addspec:finitecount}. The residual modes and their
multiplicities are included by Lemma~\ref{poly:mass}.
The trace $D_{X,Q}$ remains nonnegative; the first correction enters
only when its expectation is expanded.

Use $K$ and $Q_g$ from \eqref{poly:log-degree}. For all sufficiently
large $g$, independently of $n$, the bounds in
Section~\ref{sec:logdegree} give $Q_g\ge Q_\varepsilon$,
$\log g\ge1$, and $Q_g\le2K\log g$. The threshold \eqref{addspec:pmthreshold} holds if
\begin{equation}\label{addspec:logthreshold}
 g\ge C_{\rm add}s^2+C_{\rm add}(2K\log g)^{D_3}s.
\end{equation}
In this domain $s\le\sqrt{g/C_{\rm add}}$, so the higher-order
term relative to $s^2/g$ satisfies
\begin{equation}\label{addspec:higherabsorption}
 \frac{sQ_g^{3D_0}}g
 \le\frac{Q_g^{3D_0}}{\sqrt{C_{\rm add}g}}
 \longrightarrow0
\end{equation}
uniformly over the allowed cusp counts. The exponential term and the truncation correction
$n/(\mathcal A_{g,n}Y_{Q_g})$ are $O_{f,\varepsilon}(g^{-2})$:
use $n\le g$, $\mathcal A_{g,n}\ge2\pi n$, and
\eqref{poly:exp-decay}--\eqref{poly:area-decay}.
Since $s\ge2$, these terms are bounded by a constant times $s^2/g$;
the term with $2^{-Q_g}\le1$ satisfies the same bound.
This proves \eqref{addspec:cleanbound}
under \eqref{addspec:logthreshold}, after imposing a fixed lower genus bound.

It remains to show that an eventual quadratic condition implies
\eqref{addspec:logthreshold}. Choose $g_\varepsilon$ large enough
that all the preceding requirements hold and
\begin{equation}\label{addspec:halflinear}
 (2K\log g)^{D_3}\le\sqrt{g/(2C_{\rm add})}
 \qquad(g\ge g_\varepsilon).
\end{equation}
If $g\ge2C_{\rm add}s^2$, then $C_{\rm add}s^2\le g/2$ and
\[
 C_{\rm add}(2K\log g)^{D_3}s
 \le C_{\rm add}\sqrt{g/(2C_{\rm add})}
                    \sqrt{g/(2C_{\rm add})}=g/2.
\]
Thus \eqref{addspec:logthreshold} follows, uniformly in $n$.
Increase $C_{f,\varepsilon}$ to cover both $2C_{\rm add}$ and
the count constant. This proves
\eqref{addspec:quadraticdomain}--\eqref{addspec:cleanbound}.
The probability bound uses the spectral-gap event identity and
$\ind_{\{N_\varepsilon\ge1\}}\le N_\varepsilon$.

Finally, $(n_{\max}(g)+1)^2/g\to0$ makes the quadratic domain hold
eventually for every $1\le n\le n_{\max}(g)$. Substituting this
upper bound for $s^2$ proves \eqref{addspec:uniformlittleoh}
and the asserted sequence statement.
\end{proof}

\subsection{Sharpness of the count bound and cusp range}
\label{sec:sharpness}

Short separating curves give a matching lower bound. We use
Mirzakhani's integration formula and uniform volume comparisons.
The count below for pairs of pants with two cusps also appears in
\cite[\S4.2 and the proof of Proposition~35]{SW25}.
The torus-cap argument for $n=1$ is the once-cusped analogue of the
closed-surface lower-bound argument of Anantharaman--Monk
\cite[\S9.1.1]{AMI}.

\begin{proposition}[Uniform lower bound]\label{sharp:lower}
For every fixed $0<\varepsilon<1/4$ there are constants
$c_\varepsilon>0$ and $g_\varepsilon^{\rm low}\ge3$, independent
of $n$, such that, for all integers $g\ge g_\varepsilon^{\rm low}$
and $n\ge1$,
\begin{equation}\label{sharp:global-lower}
 \PP_{g,n}(N_\varepsilon\ge1)
 \ge c_\varepsilon
       \min\left\{\frac{(n+1)^2}{g},1\right\}.
\end{equation}
\end{proposition}

\begin{proof}
Put $\tau_\varepsilon=1/4-\varepsilon$ and choose
\[
 \ell_\varepsilon=
 \min\left\{\frac12,\frac{\pi\tau_\varepsilon}{4\sinh1}\right\}.
\]
All curve counts in this proof are unoriented and count each
geodesic once.

\emph{A short cap gives a small eigenvalue.}
Suppose a simple closed geodesic of length
$\ell\le\ell_\varepsilon$ bounds either a pair of pants with two
cusps or a one-holed torus with no cusps. The cap has area $2\pi$.
The collar lemma gives an embedded collar of half-width
$\operatorname{arsinh}(1/\sinh(\ell/2))>1$, with metric
$dr^2+\ell^2\cosh^2r\,d\theta^2$, where $\theta\in\R/\Z$;
see \cite[proof of Lemma~34]{SW25}.
Let $u$ equal one on the cap, $1-r$ on its exterior collar
$0\le r\le1$, and zero elsewhere. This is an $H^1(X)$ function,
also when the cap contains cusps, since the cap has finite area.
Direct integration gives
\begin{equation}\label{sharp:cap-energy}
 \int_X|\nabla u|^2=\ell\sinh1,\qquad
 \int_Xu^2\ge2\pi,\qquad
 \int_Xu\le2\pi+\ell\sinh1\le3\pi.
\end{equation}
For $g\ge3$ and $n\ge1$, $\mathcal A_{g,n}\ge10\pi$. With
$\bar u=\mathcal A_{g,n}^{-1}\int_Xu$ we therefore have
\[
 \|u-\bar u\|_2^2
 =\int_Xu^2-\frac{(\int_Xu)^2}{\mathcal A_{g,n}}
 \ge\pi,\qquad
 \frac{\int_X|\nabla(u-\bar u)|^2}{\|u-\bar u\|_2^2}
 \le\frac{\ell_\varepsilon\sinh1}{\pi}
 \le\frac{\tau_\varepsilon}{4}.
\]
The variational principle on the orthogonal complement of
constants, together with the essential spectral threshold $1/4$,
gives a nonzero discrete $L^2$ eigenvalue below
$\tau_\varepsilon$; see \cite[Theorem~7]{SW25}.
Thus the presence of either cap implies $N_\varepsilon\ge1$.

\emph{Two-cusp pants when $n\ge2$.}
Let $P(X)$ count simple closed geodesics of length at most
$\ell_\varepsilon$ cutting off exactly two cusps and no genus.
Distinct such geodesics are disjoint by the collar lemma.
Their caps cannot be nested, since the region between them would
be an annulus bounded by distinct closed geodesics. Nor can two
caps cover $X$, since $g\ge3$. Hence the caps, and their cusp
pairs, are disjoint.

Mirzakhani's marked integration formula, with the convention of
Lemma~\ref{orb:simpleintegration}, gives
\begin{equation}\label{sharp:pants-first}
 \E_{g,n}P
 =\frac{\binom n2}{V_{g,n}}
   \int_0^{\ell_\varepsilon}
       xV_{g,n-1}(x,\mathbf0_{n-2})\,dx.
\end{equation}
For $n\ge4$, the ordered count of two distinct caps gives
\begin{equation}\label{sharp:pants-second}
 \E_{g,n}[P(P-1)]
 =\frac{\binom n2\binom{n-2}2}{V_{g,n}}
  \int_{[0,\ell_\varepsilon]^2}
       xyV_{g,n-2}(x,y,\mathbf0_{n-4})\,dx\,dy.
\end{equation}
For $n=2,3$ this factorial moment is zero.
The second formula has no exchange factor because the pair is
ordered. Both formulas also follow from
\cite[Theorem~6]{SW25}, with $V_{0,3}=1$.

For every stable type, the volume comparisons
\cite[Lemma~14]{SW25} give
\begin{equation}\label{sharp:volume-comparisons}
 V_{g,m}\le V_{g,m}(\boldsymbol x)
       \le e^{\sum_i x_i/2}V_{g,m},\qquad
 \frac{V_{g,m+1}}{(2g-2+m)V_{g,m}}\asymp1,
\end{equation}
with absolute comparison constants. These bounds have no
restriction on the growth of $m$.
Consequently, with
$q_{g,n}=n(n-1)/(g+n)$, there are positive constants
$a_\varepsilon,A_\varepsilon,B_\varepsilon$, independent of $g,n$,
such that
\begin{equation}\label{sharp:pants-moments}
 a_\varepsilon q_{g,n}\le\E_{g,n}P
       \le A_\varepsilon q_{g,n},\qquad
 \E_{g,n}[P(P-1)]\le B_\varepsilon q_{g,n}^2.
\end{equation}
Indeed, products of one or two adjacent-volume ratios are comparable to
$(g+n)^{-1}$ or $(g+n)^{-2}$, respectively, and the bounded
boundary lengths contribute only constants depending on
$\ell_\varepsilon$.
Cauchy--Schwarz now yields
\begin{equation}\label{sharp:second-moment}
 \PP_{g,n}(P\ge1)
 \ge\frac{(\E_{g,n}P)^2}{\E_{g,n}P^2}
 \ge c'_\varepsilon\frac{q_{g,n}}{1+q_{g,n}}.
\end{equation}
For $n\ge2$, putting $t=(n+1)^2/g$, we have
\[
 \frac{q_{g,n}}{1+q_{g,n}}
 =\frac{n(n-1)}{g+n^2}
 \ge\frac29\frac{t}{1+t}
 \ge\frac19\min\{t,1\}.
\]
Together with the cap test, this proves
\eqref{sharp:global-lower} for $n\ge2$.

\emph{Torus caps when $n=1$.}
Let $T(X)$ count simple closed geodesics of length at most
$\ell_\varepsilon$ cutting off a one-holed torus containing no
cusp. Distinct counted boundaries are disjoint. The same annulus
argument excludes nesting of the caps, and two genus-one caps
cannot cover a surface of genus at least three. The complementary
types for one cap and an ordered pair are $(g-1,2)$ and $(g-2,3)$.
Marked integration therefore gives
\begin{align}
 \E_{g,1}T
 &=\frac1{V_{g,1}}\int_0^{\ell_\varepsilon}
      xV_{1,1}(x)V_{g-1,2}(x,0)\,dx,
       \label{sharp:torus-first}\\
 \E_{g,1}[T(T-1)]
 &=\frac1{V_{g,1}}\int_{[0,\ell_\varepsilon]^2}
      xyV_{1,1}(x)V_{1,1}(y)V_{g-2,3}(x,y,0)\,dx\,dy.
       \label{sharp:torus-second}
\end{align}
Here $V_{1,1}(x)=(x^2+4\pi^2)/48$ already includes the torus
involution factor, as in Lemma~\ref{orb:simpleintegration}.
The ordered pair again has no exchange factor.
The boundary comparisons in \eqref{sharp:volume-comparisons}
and the asymptotic for a fixed number of marked points
\cite[Theorem~1.2]{MZ15} imply, for all sufficiently large $g$,
\begin{equation}\label{sharp:torus-moments}
 \E_{g,1}T\asymp_\varepsilon
       \frac{V_{g-1,2}}{V_{g,1}}\asymp_\varepsilon g^{-1},
 \qquad
 \E_{g,1}[T(T-1)]\ll_\varepsilon
       \frac{V_{g-2,3}}{V_{g,1}}\asymp g^{-2}.
\end{equation}
Applying the same second-moment inequality gives
$\PP_{g,1}(T\ge1)\ge c''_\varepsilon/g$.
Since $(n+1)^2=4$ when $n=1$, the cap test proves
\eqref{sharp:global-lower} after decreasing $c_\varepsilon$
and choosing $g_\varepsilon^{\rm low}$ sufficiently large, with
$g_\varepsilon^{\rm low}\ge4$.
\end{proof}

\begin{corollary}[Sharp count and cusp range]\label{sharp:consequences}
For every fixed $0<\varepsilon<1/4$ there are constants
$c_\varepsilon,C_\varepsilon>0$, $K_\varepsilon\ge1$ and
$g_\varepsilon\ge3$, independent of $n$, such that all integers
$g\ge3$ and $n\ge1$ with $g\ge g_\varepsilon$ and
$g\ge K_\varepsilon(n+1)^2$ satisfy
\begin{equation}\label{sharp:two-sided}
 c_\varepsilon\frac{(n+1)^2}{g}
 \le\PP_{g,n}(N_\varepsilon\ge1)
 \le\E_{g,n}N_\varepsilon
 \le C_\varepsilon\frac{(n+1)^2}{g}.
\end{equation}
Moreover, for any positive integer sequence $n(g)$ and
$X_g$ with law $\PP_{g,n(g)}$, as $g\to\infty$,
\begin{equation}\label{sharp:criterion}
 \SpG(X_g)\xrightarrow{\;\mathrm{prob.}\;}\frac14
 \quad\Longleftrightarrow\quad n(g)=o(\sqrt g).
\end{equation}
\end{corollary}

\begin{proof}
Combine Proposition~\ref{sharp:lower} with
Theorem~\ref{thm:subcritical}, taking $K_\varepsilon$ at least one
and at least the quadratic-threshold constant in that theorem,
and increasing the lower genus bound. Then $(n+1)^2/g\le1$,
which proves \eqref{sharp:two-sided}.
For \eqref{sharp:criterion}, sufficiency follows from
Theorem~\ref{thm:subcritical} and $\SpG(X_g)\le1/4$.
Conversely, convergence to $1/4$ implies
$\PP_{g,n(g)}(N_{1/8}\ge1)=
\PP_{g,n(g)}(\SpG\le1/8)\to0$.
Proposition~\ref{sharp:lower} with $\varepsilon=1/8$ then forces
$(n(g)+1)^2/g\to0$, proving necessity.
\end{proof}

\par
\end{document}